%% file: 1105.4307.English.tex
\documentclass{amsart}
\input{1105.4307}

\begin{document}
\title{Algebra with Conjugation}

\begin{abstract}
In the paper, I consider properties and mappings of free algebra
with unit. I consider also conjugation of free algebra
with unit.
\end{abstract}

\ShowEq{contents}
\end{document}

%% file: 1105.4307.tex
\input{Commands}
\def\BookNumber{1105.4307}
\xRefDef{0909.0855}
\xRefDef{8433-5163}

\DefEq
{
\maketitle
\tableofcontents
\input{\FilePrefix Preface.1105.4307.English}
\input{\FilePrefix Convention.English}

\input{\FilePrefix Algebra.Conjugate.English}
\input{\FilePrefix Biblio.English}
\input{\FilePrefix Index.English}
\input{\FilePrefix Symbol.English}
}
{contents}

%% file: Commands.tex
\def\Defined{}
\scrollmode
\ifx\FilePrefix\undefined
\newcommand{\FilePrefix}{}
\fi
\ifx\UseRussian\Defined
	\usepackage{cmap}
		\usepackage[T2A,T2B]{fontenc}
	\usepackage[cp1251]{inputenc}
	\usepackage[english,russian]{babel}
\fi
\ifx\Presentation\Defined
	\usepackage[margin=1cm]{geometry}
\fi
\ifx\PublishBook\Defined
	\def\PrintPaper{}
	\usepackage{setspace}
	\ifx\UseRussian\undefined
		\usepackage{pslatex}
	\fi
	\usepackage[margin=2cm]{geometry}
\fi
\usepackage{footmisc}
\usepackage[all]{xy}
\usepackage{color}
\ifx\PrintPaper\undefined
	\definecolor{UrlColor}{rgb}{.9,0,.3}
	\definecolor{SymbColor}{rgb}{.4,0,.9}
	\definecolor{IndexColor}{rgb}{1,.3,.6}
	\newcommand\BlueText[1]{\textcolor{blue}{#1}}
	
\else
	\definecolor{UrlColor}{rgb}{.1,.1,.1}
	\definecolor{SymbColor}{rgb}{.1,.1,.1}
	\definecolor{IndexColor}{rgb}{.1,.1,.1}
	\newcommand\BlueText[1]{#1}
	
\fi

\usepackage{xr-hyper}
	\usepackage[unicode]{hyperref}
	\hypersetup{pdfdisplaydoctitle=true}
	\hypersetup{colorlinks}
	\hypersetup{citecolor=UrlColor}
	\hypersetup{urlcolor=UrlColor}
	\hypersetup{linkcolor=UrlColor}
	\hypersetup{pdffitwindow=true}
	\hypersetup{pdfnewwindow=true}
	\hypersetup{pdfstartview={FitH}}

\def\hyph{\penalty0\hskip0pt\relax-\penalty0\hskip0pt\relax}
\def\Hyph{-\penalty0\hskip0pt\relax}

\def\ValueOff{off}
\def\ValueOn{on}
\def\Items#1{\ItemList#1,LastItem,}%
\def\LastItem{LastItem}%
\def\ItemList#1,{\def\ViewBook{#1}%
	\ifx\ViewBook\LastItem%
	\else%
		\ifx\ViewBook\BookNumber%
			\def\Semafor{on}%
		\fi%
		\expandafter\ItemList%
	\fi%
}%

\newcommand{\ePrints}[1]
{%
	\def\Semafor{off}%
	\Items{#1}%
}%

\newcommand{\Basis}[1]{\overline{\overline{#1}}{}}
\newcommand{\Vector}[1]{\overline{#1}{}}
\ifx\PrintPaper\undefined
	\newcommand{\gi}[1]{\boldsymbol{\textcolor{IndexColor}{#1}}}
\else
	\newcommand{\gi}[1]{\boldsymbol{#1}}
\fi

\newcommand{\VX}[1]{\Vector{#1}_{[1]}}
\makeatletter
\newcommand{\NameDef}[1]{%
	\expandafter\gdef\csname #1\endcsname%
}%
\newcommand{\ShowSymbol}[1]{%
	\@nameuse{ViewSymbol#1}%
}%
\ifx\PrintPaper\undefined
	\newcommand{\symb}[3]{%
		\@ifundefined{ViewSymbol#3}{%
			\NameDef{ViewSymbol#3}{\textcolor{SymbColor}{#1}}%
			\NameDef{RefSymbol#3}{\pageref{symbol: #3}}%
			\@namedef{LabelSymbol#3}{\label{symbol: #3}}%
		}{%
			\NameDef{RefSymbol#3}{\pageref{symbol: #3}, \pageref{symbol 1: #3}}%
			\@namedef{LabelSymbol#3}{\label{symbol 1: #3}}%
		}%
		\ifcase#2
		\or
			$\@nameuse{ViewSymbol#3}$%
		\or
			\[\@nameuse{ViewSymbol#3}\]%
		\else%
		\fi%
		\@nameuse{LabelSymbol#3}%
	}%
\else
	\newcommand{\symb}[3]{%
		\@ifundefined{ViewSymbol#3}{%
			\NameDef{ViewSymbol#3}{#1}%
			\NameDef{RefSymbol#3}{\pageref{symbol: #3}}%
			\@namedef{LabelSymbol#3}{\label{symbol: #3}}%
		}{%
			\message {error: extra symb #3}
		}%
		\ifcase#2
		\or
			$\@nameuse{ViewSymbol#3}$%
		\or
			\[\@nameuse{ViewSymbol#3}\]%
		\else%
		\fi%
		\@nameuse{LabelSymbol#3}%
	}%
\fi
\newcommand{\DefEq}[2]{%
	\@ifundefined{ViewEq#2}{%
		\NameDef{ViewEq#2}{#1}%
	}{%
	}%
}%
\newcommand{\DefEquation}[2]{%
	\DefEq
	{
	\begin{equation}
	#1
	\EqLabel{#2}
	\end{equation}
	}
	{#2}
}%
\newcommand\EqRef[1]{\eqref{eq: #1}}%
\newcommand\EqLabel[1]{\label{eq: #1}}%
\newcommand\ShowEq[1]{%
	\@ifundefined{ViewEq#1}{%
		\message {error: missed ShowEq #1}
  }{%
	\@nameuse{ViewEq#1}%
	}%
}%
\newcommand\DrawEq[2]{%
	\@ifundefined{ViewEq#1}{%
		\message {error: missed ShowEq #1}
  }{%
	\def\Temp{}
	\def\Tempa{#2}
	\ifx\Tempa\Temp
		\[%
		\@nameuse{ViewEq#1}%
		\]%
	\else
		\begin{equation}%
		\@nameuse{ViewEq#1}%
		\EqLabel{#1, #2}%
		\end{equation}%
	\fi%
	}%
}%
\makeatother
 
\DeclareMathOperator{\id}{\mathrm{id}} 
 
\newcommand{\subs}{${}_*$\Hyph}
\newcommand{\sups}{${}^*$\Hyph}

\newcommand{\CRstar}{{}^*{}_*}
\newcommand{\RCstar}{{}_*{}^*}
\newcommand{\CRcirc}{{}^{\circ}{}_{\circ}}
\newcommand{\RCcirc}{{}_{\circ}{}^{\circ}}

\newcommand{\RC}{$\RCstar$\Hyph}
\newcommand{\CR}{$\CRstar$\Hyph}
\newcommand{\drc}{$D\RCstar$\Hyph}
\newcommand{\Drc}{$\mathcal D\RCstar$\Hyph}
\newcommand{\dcr}{$D\CRstar$\hyph}
\newcommand{\rcd}{$\RCstar D$\Hyph}
\newcommand{\crd}{$\CRstar D$\Hyph}

\newcommand{\Acr}{$A\CRcirc$\Hyph}

\newcommand\sT{$\star T$\Hyph}%
\newcommand\Ts{$T\star$\Hyph}%
\newcommand\sD{$\star D$\Hyph}%
\newcommand\Ds{$D\star$\Hyph}%

\newcommand\pC[2]{{}_{#1\cdot #2}}%
\newcommand\DrcPartial[1]%
{%
	\def\tempa{}%
	\def\tempb{#1}%
	\ifx\tempa\tempc%
		(\partial\RCstar)%
	\else%
		({}_{#1}\partial\RCstar)%
	\fi%
}%
\newcommand\crDPartial[1]%
{%
	\def\tempa{}%
	\def\tempb{#1}%
	\ifx\tempa\tempc%
		(\CRstar\partial)%
	\else%
		(\CRstar\partial_{#1})%
	\fi%
}%
\newcommand\StandPartial[3]%
{%
	\frac{\partial^{\gi{#3}} #1}{\partial #2}%
}%

	\renewcommand{\uppercasenonmath}[1]{}

\makeatletter
\newcommand\@dotsep{4.5}
\def\@tocline#1#2#3#4#5#6#7
{\relax
		\par \addpenalty\@secpenalty\addvspace{#2}%
		\begingroup 
		\@ifempty{#4}{%
			\@tempdima\csname r@tocindent\number#1\endcsname\relax
		}{%
			\@tempdima#4\relax
		}%
		\parindent\z@ \leftskip#3\relax \advance\leftskip\@tempdima\relax
		\rightskip\@pnumwidth plus1em \parfillskip-\@pnumwidth
		#5\leavevmode\hskip-\@tempdima #6\relax
		\leaders\hbox{$\m@th
		\mkern \@dotsep mu\hbox{.}\mkern \@dotsep mu$}\hfill
		\hbox to\@pnumwidth{\@tocpagenum{#7}}\par
		\nobreak
		\endgroup
}
\makeatother 

\ifx\PrintBook\undefined

	\makeatletter
		\renewcommand{\@indextitlestyle}{%
			\twocolumn[\section{\indexname}]%
			\def\IndexSpace{off}%
		}
	\makeatother 
		\ifx\PrintPaper\undefined
			\thanks{\href{mailto:Aleks\_Kleyn@MailAPS.org}{Aleks\_Kleyn@MailAPS.org}}
			\ePrints{1102.1776,Finsler}
			\ifx\Semafor\ValueOff
					\thanks{\ \ \ \url{http://sites.google.com/site/AleksKleyn/}}
					\thanks{\ \ \ \url{http://arxiv.org/a/kleyn\_a\_1}}
					\thanks{\ \ \ \url{http://AleksKleyn.blogspot.com/}}
			\fi
		\fi
\else

\makeatletter
\def\@maketitle{%
  \cleardoublepage \thispagestyle{empty}%
  \begingroup \topskip\z@skip
  \null\vfil
  \begingroup
  \LARGE\bfseries \centering
  \openup\medskipamount
  \@title\par\vspace{24pt}%
  \def\and{\par\medskip}\centering
  \mdseries\authors\par\bigskip
  \endgroup
  \vfil\vspace{24pt}
  \ifx\@empty\addresses \else \@setaddresses \fi
  \vfil
  \ifx\@empty\@dedicatory
  \else \begingroup
    \centering{\footnotesize\itshape\@dedicatory\@@par}%
    \endgroup
  \fi
  \vfill
  \newpage\thispagestyle{empty}
  \@setabstract
  \begin{center}
    \ifx\@empty\@subjclass\else\@setsubjclass\fi
    \ifx\@empty\@keywords\else\@setkeywords\fi
    \ifx\@empty\@translators\else\vfil\@settranslators\fi
    \ifx\@empty\thankses\else\vfil\@setthanks\fi
  \end{center}
  \vfil
  \endgroup}
\makeatother 

	\makeatletter
	\renewcommand{\@indextitlestyle}{%
		\twocolumn[\chapter{\indexname}]%
		\def\IndexSpace{off}%
		\let\@secnumber\@empty
		\chaptermark{\indexname}%
	}
	\makeatother 
	\email{\href{mailto:Aleks\_Kleyn@MailAPS.org}{Aleks\_Kleyn@MailAPS.org}}
	\ePrints{1102.1776,Finsler}
	\ifx\Semafor\ValueOff
		\urladdr{\url{http://sites.google.com/site/alekskleyn/}}
		\urladdr{\url{http://arxiv.org/a/kleyn\_a\_1}}
		\urladdr{\url{http://AleksKleyn.blogspot.com/}}
	\fi
\fi

\ifx\SelectlEnglish\undefined
	\ifx\UseRussian\undefined
		\def\SelectlEnglish{}
	\fi
\fi

\newcommand\arXivRef{http://arxiv.org/PS_cache/}
\newcommand\wRefDef[2]
	{
	\def\Tempa{#1}
	\def\Tempb{0405.027}
	\ifx\Tempa\Tempb
	\def\wRef{\arXivRef gr-qc/pdf/0405/0405027v3.pdf}
	\fi
	\def\Tempb{0405.028}
	\ifx\Tempa\Tempb
	\def\wRef{\arXivRef gr-qc/pdf/0405/0405028v5.pdf}
	\fi
	\def\Tempb{0412.391}
	\ifx\Tempa\Tempb
	\def\wRef{\arXivRef math/pdf/0412/0412391v4.pdf}
	\fi
	\def\Tempb{0612.111}
	\ifx\Tempa\Tempb
	\def\wRef{\arXivRef math/pdf/0612/0612111v2.pdf}
	\fi
	\def\Tempb{0701.238}
	\ifx\Tempa\Tempb
	\def\wRef{\arXivRef math/pdf/0701/0701238v4.pdf}
	\fi
	\def\Tempb{0702.561}
	\ifx\Tempa\Tempb
	\def\wRef{\arXivRef math/pdf/0702/0702561v3.pdf}
	\fi
	\def\Tempb{0707.2246}
	\ifx\Tempa\Tempb
	\def\wRef{\arXivRef arxiv/pdf/0707/0707.2246v2.pdf}
	\fi
	\def\Tempb{0803.3276}
	\ifx\Tempa\Tempb
	\def\wRef{\arXivRef arxiv/pdf/0803/0803.3276v3.pdf}
	\fi
	\def\Tempb{0812.4763}
	\ifx\Tempa\Tempb
	\def\wRef{\arXivRef arxiv/pdf/0812/0812.4763v6.pdf}
	\fi
 	\def\Tempb{0906.0135}
	\ifx\Tempa\Tempb
	\def\wRef{\arXivRef arxiv/pdf/0906/0906.0135v3.pdf}
	\fi
 	\def\Tempb{0909.0855}
	\ifx\Tempa\Tempb
	\def\wRef{\arXivRef arxiv/pdf/0909/0909.0855v5.pdf}
	\fi
 	\def\Tempb{0912.3315}
	\ifx\Tempa\Tempb
	\def\wRef{\arXivRef arxiv/pdf/0912/0912.3315v2.pdf}
	\fi
 	\def\Tempb{0912.4061}
	\ifx\Tempa\Tempb
	\def\wRef{\arXivRef arxiv/pdf/0912/0912.4061v2.pdf}
	\fi
 	\def\Tempb{1003.1544}
	\ifx\Tempa\Tempb
	\def\wRef{\arXivRef arxiv/pdf/1003/1003.1544v2.pdf}
	\fi
 	\def\Tempb{1104.5197}
	\ifx\Tempa\Tempb
	\def\wRef{\arXivRef arxiv/pdf/1104/1104.5197.pdf}
	\fi
 	\def\Tempb{1107.1139}
	\ifx\Tempa\Tempb
	\def\wRef{\arXivRef arxiv/pdf/1104/1107.1139.pdf}
	\fi
 	\def\Tempb{1107.5037}
	\ifx\Tempa\Tempb
	\def\wRef{\arXivRef arxiv/pdf/1107/1107.5037.pdf}
	\fi
 	\def\Tempb{1202.6021}
	\ifx\Tempa\Tempb
	\def\wRef{\arXivRef arxiv/pdf/1202/1202.6021.pdf}
	\fi
 	\def\Tempb{8433-5163}
	\ifx\Tempa\Tempb
	\def\wRef{http://www.amazon.com/}
	\fi
 	\def\Tempb{8443-0072}
	\ifx\Tempa\Tempb
	\def\wRef{http://www.amazon.com/}
	\fi
	\externaldocument[#1-#2-]{\FilePrefix #1.#2}[\wRef]
	}
\newcommand\LanguagePrefix{}%
\makeatletter
\newcommand\StartLabelItem[1]
{
\def\theenumi{\ref*{#1}.\@arabic\c@enumi}
\def\labelenumi{\theenumi:}
}
\newcommand\StopLabelItem
{
\def\theenumi{\@arabic\c@enumi}
\def\labelenumi{(\theenumi)}
}
\makeatother

\ifx\SelectlEnglish\undefined
	\newcommand\CurrentLanguage{Russian.}%
	\author{Александр Клейн}
		\newtheorem{theorem}{Теорема}[section]
		\newtheorem{corollary}[theorem]{Следствие}
		\newtheorem{convention}[theorem]{Соглашение}
		\theoremstyle{definition}
		\newtheorem{definition}[theorem]{Определение}
		\newtheorem{example}[theorem]{Пример}
		\newtheorem{xca}[theorem]{Exercise}
		\theoremstyle{remark}
		\newtheorem{remark}[theorem]{Замечание}
		\newtheorem{lemma}[theorem]{Лемма}
		
	\newcommand\xRefDef[1]
		{
		\wRefDef{#1}{Russian}
		\NameDef{xRefDef#1}{}%
		}
	\makeatletter
	\newcommand\xRef[2]%
	{%
		\@ifundefined{xRefDef#1}{%
		\ref{#2}%
		}{%
		\citeBib{#1}-\ref{#1-Russian-#2}%
		}%
	}%
	\newcommand\xEqRef[2]%
	{%
		\@ifundefined{xRefDef#1}{%
		\eqref{eq: #2}%
		}{%
		\citeBib{#1}-\eqref{#1-Russian-eq: #2}%
		}%
	}%
	\makeatother
	\ifx\PrintBook\undefined
		\newcommand{\BibTitle}{%
			\section{Список литературы}%
		}
	\else
		\newcommand{\BibTitle}{%
			\chapter{Список литературы}%
		}
	\fi
\else
	\newcommand\CurrentLanguage{English.}%
	\author{Aleks Kleyn}
		\newtheorem{theorem}{Theorem}[section]
		\newtheorem{corollary}[theorem]{Corollary}
		\newtheorem{convention}[theorem]{Convention}
		\theoremstyle{definition}
		\newtheorem{definition}[theorem]{Definition}
		\newtheorem{example}[theorem]{Example}
		
		\theoremstyle{remark}
		\newtheorem{remark}[theorem]{Remark}
		
	\newcommand\xRefDef[1]
		{
		\wRefDef{#1}{English}
		\NameDef{xRefDef#1}{}%
		}
	\makeatletter
	\newcommand\xRef[2]%
	{%
		\@ifundefined{xRefDef#1}{%
		\ref{#2}%
		}{%
		\citeBib{#1}-\ref{#1-English-#2}%
		}%
	}%
	\newcommand\xEqRef[2]%
	{%
		\@ifundefined{xRefDef#1}{%
		\eqref{eq: #2}%
		}{%
		\citeBib{#1}-\eqref{#1-English-eq: #2}%
		}%
	}%
	\makeatother
	\ifx\PrintBook\undefined
		\newcommand{\BibTitle}{%
			\section{References}%
		}
	\else
		\newcommand{\BibTitle}{%
			\chapter{References}%
		}
	\fi
\fi

\ifx\PrintBook\undefined
		\numberwithin{Hfootnote}{section}
\else
	\numberwithin{section}{chapter}
	\numberwithin{footnote}{chapter}
	\numberwithin{Hfootnote}{chapter}
\fi

\ifx\Presentation\undefined
	\numberwithin{equation}{section}
	\numberwithin{figure}{section}
	\numberwithin{table}{section}
	\numberwithin{Item}{section}
\fi

\makeatletter
\newcommand\org@maketitle{}
\let\org@maketitle\maketitle
\def\maketitle{%
	\hypersetup{pdftitle={\@title}}%
	\hypersetup{pdfauthor={\authors}}%
	\hypersetup{pdfsubject=\@keywords}%
	\ifx\UseRussian\Defined
		\pdfbookmark[1]{\@title}{TitleRussian}
	\else
		\pdfbookmark[1]{\@title}{TitleEnglish}
	\fi
	\org@maketitle
}
\def\make@stripped@name#1{%
	\begingroup
		\escapechar\m@ne
		\global\let\newname\@empty
		\protected@edef\Hy@tempa{\CurrentLanguage #1}%
		\edef\@tempb{%
			\noexpand\@tfor\noexpand\Hy@tempa:=%
			\expandafter\strip@prefix\meaning\Hy@tempa
		}%
		\@tempb\do{%
			\if\Hy@tempa\else
				\if\Hy@tempa\else
					\xdef\newname{\newname\Hy@tempa}%
				\fi
			\fi
		}%
	\endgroup
}%
\newenvironment{enumBib}{%
	\BibTitle
	\advance\@enumdepth \@ne
	\edef\@enumctr{enum\romannumeral\the\@enumdepth}\list
	{\csname biblabel\@enumctr\endcsname}{\usecounter
	{\@enumctr}\def\makelabel##1{\hss\llap{\upshape##1}}}
}{%
	\endlist
}

\makeatletter

\newcommand{\BiblioItem}[2]
{
	\def\Semafor{off}
	\@ifundefined{\LanguagePrefix ViewCite#1}{}{%
		\def\Semafor{on}%
	}%
	\ifx\Semafor\ValueOff
		\@ifundefined{xRefDef#1}{}{%
		\def\Semafor{on}%
		}%
	\fi
	\ifx\Semafor\ValueOn
		\ifx\IndexState\ValueOff
			\begin{enumBib}
			\def\IndexState{on}
		\fi
		\item \label{\LanguagePrefix bibitem: #1}#2%
	\fi
}
\makeatother
\newcommand{\OpenBiblio}
{
	\def\IndexState{off}
}
\newcommand{\CloseBiblio}
{
	\ifx\IndexState\ValueOn
		\end{enumBib}
		\def\IndexState{off}
	\fi
}

\makeatletter
\def\StartCite{[}%
\def\citeBib#1{\protect\showCiteBib#1,endCite,}%
\def\endCite{endCite}%
\def\showCiteBib#1,{\def\temp{#1}%
\ifx\temp\endCite
]%
\def\StartCite{[}%
\else
	\StartCite\LanguagePrefix \ref{\LanguagePrefix bibitem: #1}%
	\@ifundefined{\LanguagePrefix ViewCite#1}{%
		\NameDef{\LanguagePrefix ViewCite#1}{}%
	}{%
	}%
	\def\StartCite{, }%
\expandafter\showCiteBib%
\fi}%
\makeatother

\newcommand{\arp}{\ar @{-->}}
\newcommand{\ars}{\ar @{.>}}
\newcommand\Bundle[1]{{\mathbb #1}}
\newcommand{\bundle}[4]%
{%
	\def\tempa{}%
	\def\tempb{#3}%
	\def\tempc{#1}%
	\ifx\tempa\tempb%
		\ifx\tempa\tempc%
			#2%
		\else%
			\xymatrix{#2:#1\arp[r]&#4}%
		\fi%
	\else%
		\ifx\tempa\tempc%
			#2[#3]%
		\else%
			\xymatrix{#2[#3]:#1\arp[r]&#4}%
		\fi%
	\fi%
}%
\makeatletter
\newcommand{\AddIndex}[2]%
{%
	\@ifundefined{RefIndex#2}{%
		\NameDef{RefIndex#2}{\pageref{index: #2}}%
		\label{index: #2}%
	}{%
		\NameDef{RefIndex#2}{\pageref{index: #2}, \pageref{index 1: #2}}%
		\label{index 1: #2}%
	}%
	{\bf #1}%
}%
\newcommand{\Index}[2]%
{%
	\def\Semafor{off}%
	\@ifundefined{RefIndex#2}{%
	}{%
		\def\Semafor{on}
	}%
	\ifx\Semafor\ValueOn%
		\def\tempa{}%
		\def\tempb{#2}%
		\ifx\IndexState\ValueOff%
			\begin{theindex}%
			\def\IndexState{on}%
		\fi%
		\ifx\IndexSpace\ValueOn%
			\indexspace%
			\def\IndexSpace{off}%
		\fi%
		\item #1%
		\ifx\tempa\tempb%
		\else%
			\ \@nameuse{RefIndex#2}%
		\fi%
	\fi%
}%

\newcommand{\Symb}[2]
{
	\def\Semafor{off}
	\@ifundefined{ViewSymbol#2}{%
	}{%
		\def\Semafor{on}
	}%
	\ifx\Semafor\ValueOn
		\ifx\IndexState\ValueOff
			\begin{theindex}
			\def\IndexState{on}
		\fi
		\ifx\IndexSpace\ValueOn
			\indexspace
			\def\IndexSpace{off}
		\fi
		\item $\displaystyle\@nameuse{ViewSymbol#2}$\ \ #1
		\@nameuse{RefSymbol#2}%
	\fi
}

\makeatother

\newcommand{\SetIndexSpace}%
{%
	\def\IndexSpace{on}%
}%

\newcommand{\OpenIndex}
{
	\def\IndexState{off}
}
\newcommand{\CloseIndex}
{
	\ifx\IndexState\ValueOn
		\end{theindex}
		\def\IndexState{off}
	\fi
}

\def\LastMemo{LastMemo}%
\def\MemoList#1//{\def\temp{#1}%
	\ifx\temp\LastMemo
	\else%
	\begin{flushright}
\begin{minipage}{200pt}
\setlength{\parindent}{5mm}
		\BlueText{#1}%
\end{minipage}
	\end{flushright}
		\expandafter\MemoList%
	\fi%
}%

%% file: Preface.1105.4307.English.tex

\section{Preface}

When I started to write the book \citeBib{0812.4763},
I realized that classical definition of a linear mapping
\[
y=ax\ \ \ y=xa
\]
restricts our possibilities to study the derivative of mapping
into division ring.
Intuitively, it was clear that we can represent
the linear mapping into division ring as a sum of terms like
\[axb\]
Following study of the structure of linear mapping into division ring
was finished in the book \citeBib{8433-5163}.
I expanded the circle of considered algebras
and started to study free algebra over commutative ring.

If the set $A$ is $D$\Hyph algebra,
then we can define few algebraic structures on the set $A$.
According to the definition
\xRef{8433-5163}{definition: algebra over ring},
$A$ is $D$\Hyph module.
If we consider only sum and product
in $D$\Hyph algebra $A$, then we can consider $D$\Hyph algebra
$A$ as ring (as option, the ring maybe nonassociative).
According to this,
we can consider $D$\Hyph algebra $A$
as $A\star$\Hyph module of dimension $1$.
This is why it is important to study mappings that preserve one or another
algebraic structure on the set $A$.

In this paper, I considered  $A\star$\Hyph linear mapping
of $D\star$\Hyph algebra $A$.

The concept of algebra with conjugation fits naturally into the theory
of $D$\Hyph algebra with unit. However I believe that we can define
conjugation even in case if algebra does not
have unit. However the answer on this question requires further
research.

%% file: Convention.English.tex
\input{\FilePrefix Convention.Eq}

\section{Conventions}

\ePrints{0812.4763,0906.0135,0908.3307,0909.0855,0912.3315,1003.1544}
\Items{1006.2597}
\ifx\Semafor\ValueOn
\begin{convention}
Function and mapping are synonyms. However according to
tradition, correspondence between either rings or vector
spaces is called mapping and a mapping of
either real field or quaternion algebra is called function.
\qed
\end{convention}
\fi

\ePrints{0701.238,0812.4763,0908.3307,0912.4061,1001.4852}
\Items{1003.1544}
\ifx\Semafor\ValueOn
\begin{convention}
In any expression where we use index I assume
that this index may have internal structure.
For instance, considering the algebra $A$ we enumerate coordinates of
$a\in A$ relative to basis $\Basis e$ by an index $i$.
This means that $a$ is a vector. However, if $a$
is matrix, then we need two indexes, one enumerates
rows, another enumerates columns. In the case, when index has
structure, we begin the index from symbol $\cdot$ in
the corresponding position. 
For instance, if I consider the matrix $a^i_j$ as an element of a vector
space, then I can write the element of matrix as $a^{\cdot}{}^i_j$.
\qed
\end{convention}
\fi

\ePrints{0701.238,0812.4763,0908.3307,0912.4061,1006.2597,1011.3102}
\Items{Calculus.Paper,jCalc}
\ifx\Semafor\ValueOn
\begin{convention}
I assume sum over index $s$ in expression like
\ShowEq{Sum over repeated index}
\qed
\end{convention}
\fi

\ePrints{0701.238,0812.4763,0906.0135,0908.3307,0909.0855}
\ifx\Semafor\ValueOn
\begin{convention}
We can consider division ring $D$ as $D$\Hyph vector space
of dimension $1$. According to this statement, we can explore not only
homomorphisms of division ring $D_1$ into division ring $D_2$,
but also linear maps of division rings.
This means that map is multiplicative over
maximum possible field. In particular, linear map
of division ring $D$ is multiplicative over center $Z(D)$. This statement
does not contradict with
definition of linear map of field because for field $F$ is true
$Z(F)=F$.
When field $F$ is different from
maximum possible, I explicit tell about this in text.
\qed
\end{convention}
\fi

\ePrints{0912.4061}
\ifx\Semafor\ValueOn
\begin{convention}
For given field $F$, unless otherwise stated,
we consider finite dimensional $F$\Hyph algebra.
\qed
\end{convention}
\fi

\ePrints{0701.238,0812.4763,0906.0135,0908.3307}
\ifx\Semafor\ValueOn
\begin{convention}
In spite of noncommutativity of product a lot of statements
remain to be true if we substitute, for instance, right representation by
left representation or right vector space by left
vector space.
To keep this symmetry in statements of theorems
I use symmetric notation.
For instance, I consider \Ds vector space
and \sD vector space.
We can read notation \Ds vector space
as either D\Hyph star\Hyph vector space or
left vector space.
We can read notation \Ds linear dependent vectors
as either D\Hyph star\Hyph linear dependent vectors or
vectors that are linearly dependent from left.
\qed
\end{convention}
\fi

\ePrints{0701.238,0812.4763,0906.0135,0908.3307,0909.0855,0912.4061}
\Items{1001.4852,1003.1544,1006.2597,1104.5197,1105.4307,1107.1139}
\Items{1202.6021,MQuater,jCalc}
\ifx\Semafor\ValueOn
\begin{convention}
Let $A$ be free finite
dimensional algebra.
Considering expansion of element of algebra $A$ relative basis $\Basis e$
we use the same root letter to denote this element and its coordinates.
However we do not use vector notation in algebra.
In expression $a^2$, it is not clear whether this is component
of expansion of element
$a$ relative basis, or this is operation $a^2=aa$.
To make text clearer we use separate color for index of element
of algebra. For instance,
\ShowEq{Expansion relative basis in algebra}
\qed
\end{convention}

\begin{convention}
If free finite dimensional algebra has unit, then we identify
the vector of basis $\Vector e_{\gi 0}$ with unit of algebra.
\qed
\end{convention}
\fi

\ePrints{1104.5197,1105.4307}
\ifx\Semafor\ValueOn
\begin{convention}
Although the algebra is a free module over some
ring, we do not use the vector notation
to write elements of algebra. In the case when I consider the
matrix of coordinates of element of algebra, I will use vector
notation to write corresponding element.
In order to avoid ambiguity when I use conjugation,
I denote $a^*$ element conjugated to element $a$.
\qed
\end{convention}
\fi

\ePrints{0906.0135,0912.3315,8443-0072,1111.6035,1102.5168}
\ifx\Semafor\ValueOn
\begin{convention}
In \citeBib{Cohn: Universal Algebra},
an arbitrary operation of algebra is denoted by letter $\omega$,
and $\Omega$ is the set of operations of some universal algebra.
Correspondingly, the universal algebra with the set of operations
$\Omega$ is denoted as $\Omega$\Hyph algebra.
Similar notations we see in
\citeBib{Burris Sankappanavar} with small difference
that an operation in the algebra is denoted by letter $f$
and $\mathcal F$ is the set of operations.
I preferred first case of notations because in this case it is
easier to see where I use operation.
\qed
\end{convention}
\fi

\ePrints{0906.0135,0912.3315,8443-0072}
\ifx\Semafor\ValueOn
\begin{convention}
Since the number of universal algebras
in the tower of representations is varying,
then we use vector notation for a tower of
representations. We denote the set
$(A_1,...,A_n)$ of $\Omega_i$\Hyph algebras $A_i$, $i=1$, ..., $n$
as $\Vector A$. We denote the set of representations
$(f_{1,2},...,f_{n-1,n})$ of these algebras as $\Vector f$.
Since different algebras have different type, we also
talk about the set of $\Vector{\Omega}$\Hyph algebras.
\ePrints{8443-0072}
\ifx\Semafor\ValueOn
We
\else
In relation to the set $\Vector A$,
we also use matrix notations 
that we discussed
in section \xRef{0701.238}{section: Concept of Generalized Index}.
For instance, we
\fi
use the symbol $\Vector A_{[1]}$ to denote the
set of $\Vector{\Omega}$\Hyph algebras $(A_2,...,A_n)$.
In the corresponding notation $(\VX A,\Vector f)$ of tower
of representation, we assume that $\Vector f=(f_{2,3},...,f_{n-1,n})$.
\qed
\end{convention}

\begin{convention}
Since we use vector notation for elements of the
tower of representations, we need convention about notation of operation.
We assume that we get result of operation componentwise. For instance,
\ShowEq{vector notation in tower of representations}
\qed
\end{convention}
\fi

\ePrints{8443-0072,1111.6035,0906.0135,NewAffine,1102.5168}
\ifx\Semafor\ValueOn
\begin{convention}
Let $A$ be $\Omega_1$\Hyph algebra.
Let $B$ be $\Omega_2$\Hyph algebra.
Notation
\ShowEq{A->*B}
means that there is representation of $\Omega_1$\Hyph algebra $A$
in $\Omega_2$\Hyph algebra $B$.
\qed
\end{convention}
\fi

\ePrints{0702.561,0707.2246,0803.2620}
\ifx\Semafor\ValueOn
\begin{convention}
I use arrow $\xymatrix{\arp[r]&}$ to represent
projection of bundle on diagram.
I use arrow $\xymatrix{\ars[r]&}$ to represent
section of bundle on diagram.
\qed
\end{convention}
\fi

\ePrints{0912.3315}
\ifx\Semafor\ValueOn
\begin{remark}
I believe that diagrams of maps are an important tool.
However, sometimes I want
to see the diagram as three dimensional figure
and I expect that this would increase its expressive
power. Who knows what surprises the future holds.
In 1992, at a conference in Kazan, I have described to my colleagues
what advantages the computer preparation of papers has.
8 years later I learned from the letter from Kazan that now we can
prepare paper using LaTeX.
\qed
\end{remark}
\fi

\ePrints{1001.4852,1003.1544,1006.2597,1011.3102}
\Items{Calculus.Paper,jCalc}
\ifx\Semafor\ValueOn
\begin{convention}
If, in a certain expression, we use several operations
which include the operation $\circ$, then
it is assumed that the operation $\circ$ is executed first.
Below is an example of equivalent expressions.
\ShowEq{list circ expressions}
\qed
\end{convention}
\fi


\ePrints{1107.1139}
\ifx\Semafor\ValueOn
\begin{convention}
For given $D$\Hyph algebra $A$
we define left shift
\ShowEq{left shift, D algebra}
by the equation
\ShowEq{left shift 1, D algebra}
and right shift
\ShowEq{right shift, D algebra}
by the equation
\ShowEq{right shift 1, D algebra}
\qed
\end{convention}
\fi

\ifx\PrintPaper\undefined
Without a doubt, the reader may have questions,
comments, objections. I will appreciate any response.
\fi

%% file: Convention.Eq.tex

\DefEq
{
\[
\Vector r(\Vector a)=(r_1(a_1),...,r_n(a_n))
\]
}
{vector notation in tower of representations}

\DefEq
{
\[
\xymatrix
{
A\ar[r]|{*}&B
}
\]
}
{A->*B}

\DefEq
{
\[
\begin{array}{r@{\ }lr@{\ }l}
f\circ xy&\equiv f(x)y
&
f\circ(xy)&\equiv f(xy)
\\
f\circ x+y&\equiv f(x)+y
&
f\circ (x+y)&\equiv f(x+y)
\end{array}
\]
}
{list circ expressions}

\DefEq
{
\[
a\pC s0xa\pC s1
\]
}
{Sum over repeated index}

\DefEq
{
\[
a=a^{\gi i}\Vector e_{\gi i}
\]
}
{Expansion relative basis in algebra}

\DefEq
{
\symb{a\circ}1{left shift, D algebra}
}
{left shift, D algebra}

\DefEq
{
\symb{a\star}1{right shift, D algebra}
}
{right shift, D algebra}

\DefEq
{
\[
\ShowSymbol{left shift, D algebra}x=ax
\]
}
{left shift 1, D algebra}

\DefEq
{
\[
\ShowSymbol{right shift, D algebra}x=xa
\]
}
{right shift 1, D algebra}

%% file: Algebra.Conjugate.English.tex
\input{Algebra.Conjugate.Eq}

\ifx\PrintBook\Defined
\chapter{Algebra with Conjugation}
\fi

\section{Algebra with unit}

Let $D$ be commutative ring.
Let $A$ be $D$\hyph algebra with unit $1$.
Therefore there exists effective representation
\ShowEq{representation D in A}
of the ring $D$ in the algebra $A$.

\begin{theorem}
\label{theorem: m+n-1 linear map algebra}
Let $D$ be commutative ring.
Let $A$ be $D$\hyph algebra.
Let \ShowEq{n linear map algebra} be
$n$\hyph linear mapping into $D$\hyph algebra $A$.
Let \ShowEq{m linear map algebra} be
$m$\hyph linear mapping into $D$\hyph algebra $A$.
Then if for given $i$ in the expression
\ShowEq{n linear map algebra, value}
we make the substitution
\ShowEq{m linear map algebra, value}
then the resulting expression
\ShowEq{m+n-1 linear map algebra, value}
is $(m+n-1)$\hyph linear mapping.
\end{theorem}
\begin{proof}
Since polylinear mapping is linear mapping with respect to each
argument, then the mapping $h$ is linear with respect to $x^j$, $j\ne i$.
Since the mapping $f$ is linear mapping with respect to $x^i$ and for each $k$
the mapping $g$ is linear mapping with respect to $y^k$, then according
to the theorem
\xRef{8433-5163}{theorem: product of linear mapping, algebra},
the mapping $h$ is linear mapping with respect to $y^k$.
\end{proof}

\begin{theorem}
Let $D$ be commutative ring.
Let $A$ be $D$\hyph algebra.
Then
\ShowEq{representation commutative with product}
\end{theorem}
\begin{proof}
Let $f$ be bilinear mapping corresponding to
product in $D$\hyph algebra $A$.\footnote{I consider the definition of the product
in $D$\hyph algebra according to the definition
\xRef{8433-5163}{definition: algebra over ring}.}
Then
\ShowEq{representation commutative with product, 1}
The equation \EqRef{representation commutative with product}
follows from the equation \EqRef{representation commutative with product, 1}.
\end{proof}

\begin{theorem}
\label{theorem: commutator associator maps}
Let $D$ be commutative ring.\footnote{You can see the definition of commutator
in the definition
\xRef{8433-5163}{definition: commutator of algebra}.
You can see the definition of associator
in the definition
\xRef{8433-5163}{definition: associator of algebra}.
You can see also \citeBib{Richard D. Schafer}, p. 13.}
The commutator of $D$\hyph algebra $A$ is bilinear mapping.
The associator of $D$\hyph algebra $A$ is $3$\hyph linear mapping.
\end{theorem}
\begin{proof}
According to definitions
\xRef{8433-5163}{definition: algebra over ring},
\xRef{8433-5163}{definition: commutator of algebra},
the commutator is bilinear mapping.
According to definitions
\xRef{8433-5163}{definition: algebra over ring},
\xRef{8433-5163}{definition: associator of algebra}
and the theorem \ref{theorem: m+n-1 linear map algebra},
the associator is $3$\hyph linear mapping.
\end{proof}

\begin{theorem}
\label{theorem: ring subset of center of algebra}
Let $D$ be commutative ring.
Let $A$ be $D$\hyph algebra with unit $e$.
Then
we can identify $d\in D$ and $d\,1\in A$.
In such case, $D\subseteq Z(A)$.\footnote{You can see
the definition of nucleus $N(A)$ of algebra $A$ in the definition
\xRef{8433-5163}{definition: nucleus of algebra};
you can see the definition of center $Z(A)$ of algebra $A$
in the definition
\xRef{8433-5163}{definition: center of algebra};
see also \citeBib{Richard D. Schafer}, p. 13, 14.}
\end{theorem}
\begin{proof}
From the theorem \ref{theorem: commutator associator maps}, it follows that
\ShowEq{representation commutative with product, 4}
Therefore,
\ShowEq{representation commutative with product, 5}.
If we assume $b=1$ in the equation
\EqRef{representation commutative with product}, then we get
\ShowEq{representation commutative with product, 2}
From the equation \EqRef{representation commutative with product, 2},
it follows that the commutator
has form
\ShowEq{representation commutative with product, 3}
\end{proof}

\begin{theorem}
Structural constants of
$D$\hyph algebra with unit $e$ satisfy condition
\ShowEq{structural constant re algebra, 1}
\end{theorem}
\begin{proof}
The equation
\EqRef{structural constant re algebra, 1}
follows from the equation $1d=d1=d$.
\end{proof}

\begin{theorem}
$A\star$\Hyph linear mapping of
$D$\Hyph algebra $A$ with unit has form
\ShowEq{A linear mapping of algebra A}
Let $\Basis e$ be $(\re A)$\Hyph basis of the algebra $A$.
Coordinates of the mapping $f$
\ShowEq{Re A linear mapping of algebra A}
relative to the basis $\Basis e$ satisfy to the equation
\ShowEq{Re A linear mapping of algebra A 2}
\ShowEq{Re A linear mapping of algebra A 5}
\end{theorem}
\begin{proof}
$A\star$\Hyph linear mapping satisfies to equations
\ShowEq{A linear mapping of algebra A 1}
\ShowEq{A linear mapping of algebra A 2}
From equations
\EqRef{A linear mapping of algebra A 1},
\EqRef{A linear mapping of algebra A 2},
it follows that
\ShowEq{A linear mapping of algebra A 3}
The equation
\ShowEq{Re A linear mapping of algebra A 1}
follows from the equation
\EqRef{A linear mapping of algebra A 3}.

From the equation
\EqRef{Re A linear mapping of algebra A 1}
it follows that\footnote{Let $\gi l=\gi 0$.}
\ShowEq{Re A linear mapping of algebra A 3}
From equations
\EqRef{structural constant re algebra, 1},
\EqRef{Re A linear mapping of algebra A 3},
it follows that
\ShowEq{Re A linear mapping of algebra A 4}
The equation
\EqRef{Re A linear mapping of algebra A 2}
follows from the equation
\EqRef{Re A linear mapping of algebra A 4}.

From equations
\EqRef{Re A linear mapping of algebra A},
\EqRef{Re A linear mapping of algebra A 2},
it follows that
\ShowEq{A linear mapping of algebra A 6}
The equation
\EqRef{Re A linear mapping of algebra A 5}
follows from the equations
\EqRef{A linear mapping of algebra A},
\EqRef{A linear mapping of algebra A 6}.
\end{proof}

\section{Algebra with Conjugation}

Let $D$ be commutative ring.
Let $A$ be $D$\hyph algebra with unit $e$, $A\ne D$.
According to the theorem
\ref{theorem: ring subset of center of algebra},
we can identify the ring $D$ and subalgebra of $D$\hyph algebra $A$.
Therefore, $D$\Hyph algebra $A$ has a nontrivial center $Z(A)$.

Let there exist subalgebra $F$ of algebra $A$ such that
$F\ne A$, $D\subseteq F\subseteq Z(A)$, and algebra $A$ is a free
module over the ring $F$.
Let $\Basis e$ be the basis of free module $A$ over ring $F$.
We assume that $\Vector e_{\gi 0}=1$.

Consider mappings
\ShowEq{Re Im A->A}
defined by equation
\ShowEq{Re Im A->F 0}
\ShowEq{Re Im A->F 1}
The expression
\ShowEq{Re Im A->F 2}
is called
\AddIndex{scalar of element}{scalar of algebra} $d$.
The expression
\ShowEq{Re Im A->F 3}
is called
\AddIndex{vector of element}{vector of algebra} $d$.\footnote{In the section
\xRef{0909.0855}{section: algebra CCH},
I consider an example of algebra in which there are two conjugation.
We assume that we chose conjugation.}

According to
\EqRef{Re Im A->F 1}
\ShowEq{Re Im A->F 4}
We will use notation
\ShowEq{scalar algebra of algebra}
to denote
\AddIndex{scalar algebra of algebra}{scalar algebra of algebra} $A$.

\begin{theorem}
\label{theorem: Re Im A->F}
The set
\ShowEq{vector module of algebra}
is $(\re A)$\Hyph module
which is called
\AddIndex{vector module of algebra}{vector module of algebra} $A$.
\ShowEq{A=re+im}
\end{theorem}
\begin{proof}
Let $c$, $d\in\im A$. Then $c_{\gi 0}=d_{\gi 0}=0$.
Therefore,
\[
(c+d)_{\gi 0}=c_{\gi 0}+d_{\gi 0}=0
\]
If $a\in\re A$, then
\[
(ad)_{\gi 0}=ad_{\gi 0}=0
\]
Therefore, $\im A$ is $(\re A)$\Hyph module.

Sequence of modules
\ShowEq{re->A->im}
is exact sequence.
According to the definition
\EqRef{Re Im A->F 1}
of the mapping $\re$,
following diagram is commutative
\ShowEq{re->A->im 1}
According to the statement 2 of the proposition \citeBib{Serge Lang}-III.3.2,
\ShowEq{re->A->im 2}
According to the definition
\EqRef{vector module of algebra}
\ShowEq{re->A->im 3}
The equation
\EqRef{A=re+im}
follows from the equations
\EqRef{re->A->im 2},
\EqRef{re->A->im 3}.
\end{proof}

According to the theorem \ref{theorem: Re Im A->F},
there is unique defined representation
\ShowEq{d Re Im}

\begin{definition}
\label{definition: conjugation in algebra}
The mapping
\ShowEq{conjugation in algebra}
\ShowEq{conjugation in algebra, 0}
is called
\AddIndex{conjugation in algebra}{conjugation in algebra}
provided that this mapping satisfies
\ShowEq{conjugation in algebra, 1}
$(\re A)$\Hyph algebra $A$
equipped with conjugation is called
\AddIndex{algebra with conjugation}{algebra with conjugation}.
\qed
\end{definition}

\begin{corollary}
\ShowEq{double conjugation in algebra}
\qed
\end{corollary}

\begin{theorem}
\label{theorem: d +- dstar}
\ShowEq{d+dstar}
\ShowEq{d-dstar}
\end{theorem}
\begin{proof}
The theorem follows from equations
\EqRef{d Re Im},
\EqRef{conjugation in algebra, 0}.
\end{proof}

\begin{theorem}
The $(\re A)$\Hyph algebra $A$ is
algebra with conjugation
if structural constants of
$(\re A)$\Hyph algebra $A$ satisfy condition
\ShowEq{conjugation in algebra, 5}
\ShowEq{conjugation in algebra, 4}
\end{theorem}
\begin{proof}
From equations
\EqRef{d Re Im},
\EqRef{conjugation in algebra, 0}
it follows that
\ShowEq{algebra with conjugation 1}
\ShowEq{algebra with conjugation 2}
From equations
\EqRef{conjugation in algebra, 1},
\EqRef{algebra with conjugation 1},
\EqRef{algebra with conjugation 2},
it follows that
\ShowEq{conjugation in algebra, 2}
Let $\Basis e$ be the basis of the $(\re D)$\Hyph algebra $D$.
From the equation
\EqRef{conjugation in algebra, 2},
it follows that
\ShowEq{conjugation in algebra, 3}
\ShowEq{conjugation in algebra, 4}
The equation
\EqRef{conjugation in algebra, 5}
follows from the equation
\EqRef{conjugation in algebra, 3}.
\end{proof}

\begin{corollary}
\label{corollary: conjugation in algebra}
\ShowEq{conjugation in algebra, 6}
\qed
\end{corollary}

\begin{theorem}
\ShowEq{kl>0}
\ShowEq{ek el}
\end{theorem}
\begin{proof}
From the equation
\EqRef{conjugation in algebra, 5},
it follows that
\ShowEq{ek el 2}
\ShowEq{ek el 1}
The equation
\EqRef{ek el}
follows from the equation
\EqRef{ek el 1}.
\end{proof}

\begin{theorem}
\ShowEq{kl>0}
\ShowEq{ek+el}
\ShowEq{ek-el}
\end{theorem}
\begin{proof}
The theorem follows from the equation
\EqRef{ek el}
and the theorem
\ref{theorem: d +- dstar}.
\end{proof}

\begin{theorem}
\ShowEq{d dstar}
\end{theorem}
\begin{proof}
From the equations
\EqRef{conjugation in algebra, 0},
\EqRef{conjugation in algebra, 5},
it follows that
\ShowEq{ek el 2}
\ShowEq{d dstar 1}
First $3$ items in the expression
\EqRef{d dstar 1}
have form
\ShowEq{d dstar 2}
The last item in the expression
\EqRef{d dstar 1}
has form
\ShowEq{d dstar 3}
According to the corollary
\ref{corollary: conjugation in algebra}
and statement
\EqRef{ek+el}
\ShowEq{d dstar 4}
The theorem
follows from the equation
\EqRef{d dstar 1}
and statements
\EqRef{d dstar 2},
\EqRef{d dstar 4}.
\end{proof}

We can represent the conjugation using the matrix $I$
\ShowEq{conjugation in algebra, 7}

\begin{example}
The product in the set of complex numbers $C$ is commutative.
However, complex field has subfield $R$.
The vector space $\im C$ has dimension $1$ and the basis
$\Vector e_{\gi 1}=i$. Accordingly
\[
c^*=c^{\gi 0}-c^{\gi 1}i
\]
\qed
\end{example}

\begin{example}
The division ring of quaternions $H$ has subfield $R$.
The vector space $\im H$ has dimension $3$ and the basis
\[
\begin{matrix}
\Vector e_{\gi 1}=i&
\Vector e_{\gi 2}=j&
\Vector e_{\gi 3}=k
\end{matrix}
\]
Accordingly
\[
d^*=d^{\gi 0}-d^{\gi 1}i-d^{\gi 2}j-d^{\gi 3}k
\]
\qed
\end{example}

\section{Antilinear Mapping of Algebra with Conjugation}

\begin{definition}
\label{definition: antilinear mapping of algebra with conjugation}
Let $A$ be the algebra with conjugation.
$(\re A)$\Hyph linear mapping
\ShowEq{f:A->A}
is called
\AddIndex{$A\star$\Hyph antilinear}
{antilinear mapping of algebra with conjugation},
if the mapping $f$ satisfies to the equation
\ShowEq{antilinear mapping of algebra with conjugation}
\qed
\end{definition}

\begin{theorem}
$A\star$\Hyph antilinear mapping of the algebra with conjugation $A$ has form
\ShowEq{A antilinear mapping of algebra A}
Let $\Basis e$ be $(\re A)$\Hyph basis of the algebra $A$.
Coordinates of the mapping $f$
\ShowEq{Re A antilinear mapping of algebra A}
relative to the basis $\Basis e$ satisfy to the equation
\ShowEq{Re A antilinear mapping of algebra A 2}
\ShowEq{Re A antilinear mapping of algebra A 5}
\end{theorem}
\begin{proof}
According to the definition
\ref{definition: antilinear mapping of algebra with conjugation},
$A\star$\Hyph antilinear mapping satisfies to equations
\ShowEq{A antilinear mapping of algebra A 1}
\ShowEq{A antilinear mapping of algebra A 2}
From equations
\EqRef{A antilinear mapping of algebra A 1},
\EqRef{A antilinear mapping of algebra A 2},
it follows that
\ShowEq{A antilinear mapping of algebra A 3}
The equation
\ShowEq{Re A antilinear mapping of algebra A 1}
follows from the equation
\EqRef{A antilinear mapping of algebra A 3}.

From the equation
\EqRef{Re A antilinear mapping of algebra A 1}
it follows that\footnote{Let $\gi q=\gi 0$.}
\ShowEq{Re A antilinear mapping of algebra A 3}
From equations
\EqRef{structural constant re algebra, 1},
\EqRef{conjugation in algebra, 7},
\EqRef{Re A antilinear mapping of algebra A 3},
it follows that
\ShowEq{Re A antilinear mapping of algebra A 4}
The equation
\EqRef{Re A antilinear mapping of algebra A 2}
follows from the equation
\EqRef{Re A antilinear mapping of algebra A 4}.

From equations
\EqRef{Re A antilinear mapping of algebra A},
\EqRef{Re A antilinear mapping of algebra A 2},
it follows that
\ShowEq{A antilinear mapping of algebra A 6}
The equation
\EqRef{Re A antilinear mapping of algebra A 5}
follows from the equations
\EqRef{A linear mapping of algebra A},
\EqRef{A antilinear mapping of algebra A 6}.
\end{proof}

%% file: Algebra.Conjugate.Eq.tex

\newcommand\re{\mathrm{Re}\,} 
\newcommand\im{\mathrm{Im}\,}

\DefEquation
{
d+d^*\in\re A
}
{d+dstar}

\DefEquation
{
d-d^*\in\im A
}
{d-dstar}

\DefEquation
{
\Vector e_{\gi k}\,\Vector e_{\gi l}+\Vector e_{\gi l}\,\Vector e_{\gi k}
\in\re A
}
{ek+el}

\DefEq
{
($\gi k>0$, $\gi l>0$)
}
{kl>0}

\DefEquation
{
\Vector e_{\gi k}\,\Vector e_{\gi l}-\Vector e_{\gi l}\,\Vector e_{\gi k}
\in\im A
}
{ek-el}

\DefEquation
{
\Vector r_2(\Vector C_1(\Vector v,\Vector w))
=\Vector C_2(\Vector r_2(\Vector v),\Vector r_2(\Vector w))
}
{algebra, morphism of representation 23, 1}

\DefEquation
{
\Vector e_{\gi k}\,\Vector e_{\gi l}
=(\Vector e_{\gi l}\,\Vector e_{\gi k})^*
}
{ek el}

\DefEquation
{
\begin{split}
\Vector e_{\gi k}\,\Vector e_{\gi l}
&=C^{\gi 0}_{\gi{ab}}e^{\gi a}_{\gi k}e^{\gi b}_{\gi l}\Vector e_{\gi 0}
+C^{\gi p}_{\gi{ab}}e^{\gi a}_{\gi k}e^{\gi b}_{\gi l}\Vector e_{\gi p}
\\
&=C^{\gi 0}_{\gi{ba}}e^{\gi b}_{\gi l}e^{\gi a}_{\gi k}\Vector e_{\gi 0}
-C^{\gi p}_{\gi{ba}}e^{\gi b}_{\gi l}e^{\gi a}_{\gi k}\Vector e_{\gi p}
\\
&=(C^{\gi 0}_{\gi{ba}}e^{\gi b}_{\gi l}e^{\gi a}_{\gi k}\Vector e_{\gi 0}
+C^{\gi p}_{\gi{ba}}e^{\gi b}_{\gi l}e^{\gi a}_{\gi k}\Vector e_{\gi p})^*
\\
&=(\Vector e_{\gi l}\,\Vector e_{\gi k})^*
\end{split}
}
{ek el 1}

\DefEq
{
$\Vector e_{\gi k}\,\Vector e_{\gi k}\in\re A$
}
{conjugation in algebra, 6}

\DefEquation
{
\begin{array}{l@{\ \ \ }r@{\,}c@{\,}l@{\ \ \ }r@{\,}c@{\,}l@{\ \ \ }r@{\,}l}
d^*=I\circ d
&I^{\gi 0}_{\gi k}=&&\delta^{\gi 0}_{\gi k}
&I^{\gi k}_{\gi 0}=&&\delta^{\gi k}_{\gi 0}
&\gi k&=\gi 0,...,\gi n
\\
&I^{\gi m}_{\gi k}=&-&\delta^{\gi m}_{\gi k}
&I^{\gi k}_{\gi m}=&-&\delta^{\gi k}_{\gi m}
&\gi m&=\gi 1,...,\gi n
\end{array}
}
{conjugation in algebra, 7}

\DefEquation
{
\begin{matrix}
f(c)=cd&d\in A
\end{matrix}
}
{A linear mapping of algebra A}

\DefEquation
{
f(c)
=c^{\gi l}f^{\gi i}_{\gi l}\,\Vector e_{\gi i}
}
{Re A linear mapping of algebra A}

\DefEquation
{
f^{\gi j}_{\gi k}
=f^{\gi i}_{\gi 0}C^{\gi j}_{\gi{ki}}
}
{Re A linear mapping of algebra A 2}

\DefEquation
{
d=f^{\gi i}_{\gi 0}\,\Vector e_{\gi i}
}
{Re A linear mapping of algebra A 5}

\DefEquation
{
f(ac)
=(ac)^{\gi i}f^{\gi j}_{\gi i}\Vector e_{\gi j}
=a^{\gi k}c^{\gi l}C^{\gi i}_{\gi{kl}}f^{\gi j}_{\gi i}\Vector e_{\gi j}
}
{A linear mapping of algebra A 1}

\DefEquation
{
af(c)
=a^{\gi k}(f(c))^{\gi i}C^{\gi j}_{\gi{ki}}\Vector e_{\gi j}
=a^{\gi k}c^{\gi l}f^{\gi i}_{\gi l}C^{\gi j}_{\gi{ki}}\Vector e_{\gi j}
}
{A linear mapping of algebra A 2}

\DefEquation
{
a^{\gi k}c^{\gi l}C^{\gi i}_{\gi{kl}}f^{\gi j}_{\gi i}\Vector e_{\gi j}
=a^{\gi k}c^{\gi l}f^{\gi i}_{\gi l}C^{\gi j}_{\gi{ki}}\Vector e_{\gi j}
}
{A linear mapping of algebra A 3}

\DefEquation
{
C^{\gi i}_{\gi{kl}}f^{\gi j}_{\gi i}
=f^{\gi i}_{\gi l}C^{\gi j}_{\gi{ki}}
}
{Re A linear mapping of algebra A 1}

\DefEquation
{
C^{\gi i}_{\gi{k0}}f^{\gi j}_{\gi i}
=f^{\gi i}_{\gi 0}C^{\gi j}_{\gi{ki}}
}
{Re A linear mapping of algebra A 3}

\DefEq
{
\symb{\re d}0{scalar of algebra}
\symb{\im d}0{vector of algebra}
}
{Re Im A->F 0}

\DefEquation
{
\begin{matrix}
\ShowSymbol{scalar of algebra}=d^{\gi 0}
&\ShowSymbol{vector of algebra}=d-d^{\gi 0}
&d\in D
&d=d^{\gi i}\Vector e_{\gi i}
\end{matrix}
}
{Re Im A->F 1}

\DefEq
{
$\ShowSymbol{scalar of algebra}$
}
{Re Im A->F 2}

\DefEq
{
\[
F=\{d\in A:\re d=d\}
\]
}
{Re Im A->F 4}

\DefEq
{
\symb{\im A}0{vector module of algebra}
\begin{equation}
\EqLabel{vector module of algebra}
\ShowSymbol{vector module of algebra}=\{d\in A:\re d=0\}
\end{equation}
}
{vector module of algebra}

\DefEquation
{
A=\re A\oplus\im A
}
{A=re+im}

\DefEq
{
\[
\xymatrix
{
0\ar[r]&\re A\ar[r]^{\mathrm{id}}
&A\ar[r]^{\im}&\im A\ar[r]&0
}
\]
}
{re->A->im}

\DefEq
{
\[
\xymatrix
{
\re A\ar[r]^{\id}\ar[d]_{\id}&A\ar[ld]^{\re}
\\
\re A
}
\]
}
{re->A->im 1}

\DefEquation
{
A=\id(\re A)\oplus\mathrm{ker}\,\re
}
{re->A->im 2}

\DefEquation
{
\mathrm{ker}\,\re=\{d\in A:\re d=0\}=\im A
}
{re->A->im 3}

\DefEquation
{
d=\re d+\im d
}
{d Re Im}

\DefEq
{
\symb{d^*}0{conjugation in algebra}
}
{conjugation in algebra}

\DefEq
{
$dd^*\in\re A$
}
{d dstar}

\DefEq
{
($\gi p=\gi 1$, ..., $\gi n$)
}
{ek el 2}

\DefEquation
{
\begin{matrix}
C^{\gi 0}_{\gi{kl}}=C^{\gi 0}_{\gi{lk}}
&
C^{\gi p}_{\gi{kl}}=-C^{\gi p}_{\gi{lk}}
\end{matrix}
}
{conjugation in algebra, 5}

\DefEq
{
\[
\begin{matrix}
\gi 1<\gi k<\gi n&\gi 1<\gi l<\gi n&\gi 1<\gi p<\gi n
\end{matrix}
\]
}
{conjugation in algebra, 4}

\DefEquation
{
C_{\gi{0k}}^{\gi l}=C^{\gi l}_{\gi{k0}}=\delta^{\gi k}_{\gi l}
}
{structural constant re algebra, 1}

\DefEquation
{
\begin{array}{r@{\,}l}
(cd)^*&=(\re c\ \re d+\re c\ \im d+\im c\ \re d+\im c\ \im d)^*
\\
&=\re c\ \re d-\re c\ \im d-\im c\ \re d+(\im c\ \im d)^*
\end{array}
}
{algebra with conjugation 1}

\DefEquation
{
\begin{array}{r@{\,}l}
d^*\,c^*&
=(\re d-\im d)(\re c-\im c)
\\
&=\re d\ \re c-\re d\ \im c-\im d\ \re c+\im d\ \im c
\end{array}
}
{algebra with conjugation 2}

\DefEquation
{
\begin{array}{r@{\,}l}
(C^{\gi 0}_{\gi{kl}}c^{\gi k}d^{\gi l}
+C^{\gi p}_{\gi{kl}}c^{\gi k}d^{\gi l}\Vector e_{\gi p})^*
&=C^{\gi 0}_{\gi{kl}}c^{\gi k}d^{\gi l}
-C^{\gi p}_{\gi{kl}}c^{\gi k}d^{\gi l}\Vector e_{\gi p}
\\&=
C^{\gi 0}_{\gi{kl}}d^{\gi k}c^{\gi l}
+C^{\gi p}_{\gi{kl}}d^{\gi k}c^{\gi l}\Vector e_{\gi p}
\end{array}
}
{conjugation in algebra, 3}

\DefEquation
{
(\im c\ \im d)^*=\im d\ \im c
}
{conjugation in algebra, 2}

\DefEq
{
$(d^*)^*=d$
}
{double conjugation in algebra}

\DefEquation
{
(cd)^*=d^*\,c^*
}
{conjugation in algebra, 1}

\DefEquation
{
\begin{split}
dd^*
&=(d^{\gi 0}\Vector e_{\gi 0}+d^{\gi p}\Vector e_{\gi p})
(d^{\gi 0}\Vector e_{\gi 0}+d^{\gi q}\Vector e_{\gi q})^*
\\
&=(d^{\gi 0}\Vector e_{\gi 0}+d^{\gi p}\Vector e_{\gi p})
(d^{\gi 0}\Vector e_{\gi 0}-d^{\gi q}\Vector e_{\gi q})
\\
&=d^{\gi 0}\Vector e_{\gi 0}\,d^{\gi 0}\Vector e_{\gi 0}
-d^{\gi 0}\Vector e_{\gi 0}\,d^{\gi q}\Vector e_{\gi q}
+d^{\gi p}\Vector e_{\gi p}\,d^{\gi 0}\Vector e_{\gi 0}
-d^{\gi p}\Vector e_{\gi p}\,d^{\gi q}\Vector e_{\gi q}
\\
&=(d^{\gi 0})^2\,\Vector e_{\gi 0}
-d^{\gi 0}d^{\gi p}\Vector e_{\gi p}
+d^{\gi 0}d^{\gi p}\Vector e_{\gi p}
-d^{\gi p}d^{\gi q}\Vector e_{\gi p}\,\Vector e_{\gi q}
\end{split}
}
{d dstar 1}

\DefEquation
{
-d^{\gi p}d^{\gi q}\Vector e_{\gi p}\,\Vector e_{\gi q}
=
-(d^{\gi p})^2\,\Vector e_{\gi p}\,\Vector e_{\gi p}
-\sum_{q>p}d^{\gi p}d^{\gi q}(\Vector e_{\gi p}\,\Vector e_{\gi q}
+\Vector e_{\gi q}\,\Vector e_{\gi p})
}
{d dstar 3}

\DefEquation
{
-d^{\gi p}d^{\gi q}\Vector e_{\gi p}\,\Vector e_{\gi q}
\in\re A
}
{d dstar 4}

\DefEquation
{
(d^{\gi 0})^2\,\Vector e_{\gi 0}
-d^{\gi 0}d^{\gi p}\Vector e_{\gi p}
+d^{\gi 0}d^{\gi p}\Vector e_{\gi p}
=
(d^{\gi 0})^2\,\Vector e_{\gi 0}\in\re A
}
{d dstar 2}

\DefEquation
{
\ShowSymbol{conjugation in algebra}=\re d-\im d
}
{conjugation in algebra, 0}

\DefEq
{
\symb{\re A}1{scalar algebra of algebra}
}
{scalar algebra of algebra}

\DefEq
{
$\ShowSymbol{vector of algebra}$
}
{Re Im A->F 3}

\DefEq
{
\[
f:A\rightarrow A
\]
}
{f:A->A}

\DefEquation
{
f(da)=f(a)d^*
}
{antilinear mapping of algebra with conjugation}

\DefEq
{
\[
\begin{matrix}
f(c)=dc^*=dI\circ c&d\in A
\end{matrix}
\]
}
{A antilinear mapping of algebra A}

\DefEquation
{
f(c)
=c^{\gi l}f^{\gi i}_{\gi l}\Vector e_{\gi i}
}
{Re A antilinear mapping of algebra A}

\DefEquation
{
f^{\gi j}_{\gi p}
=f^{\gi k}_{\gi 0}I^{\gi i}_{\gi p}
C^{\gi j}_{\gi{ki}}
}
{Re A antilinear mapping of algebra A 2}

\DefEquation
{
d=f^{\gi i}_{\gi 0}\,\Vector e_{\gi i}
}
{Re A antilinear mapping of algebra A 5}

\DefEquation
{
f(ac)
=(ac)^{\gi i}f^{\gi j}_{\gi i}\Vector e_{\gi j}
=a^{\gi p}c^{\gi q}C^{\gi i}_{\gi{pq}}
f^{\gi j}_{\gi i}\Vector e_{\gi j}
}
{A antilinear mapping of algebra A 1}

\DefEquation
{
f(c)a^*
=(f(c))^{\gi k}(a^*)^{\gi i}C^{\gi j}_{\gi{ki}}\Vector e_{\gi j}
=c^{\gi l}f^{\gi k}_{\gi l}I^{\gi i}_{\gi q}a^{\gi q}C^{\gi j}_{\gi{ki}}
\Vector e_{\gi j}
}
{A antilinear mapping of algebra A 2}

\DefEquation
{
a^{\gi p}c^{\gi q}C^{\gi i}_{\gi{pq}}
f^{\gi j}_{\gi i}\Vector e_{\gi j}
=c^{\gi q}f^{\gi k}_{\gi q}I^{\gi i}_{\gi p}a^{\gi p}
C^{\gi j}_{\gi{ki}}\Vector e_{\gi j}
}
{A antilinear mapping of algebra A 3}

\DefEquation
{
C^{\gi i}_{\gi{pq}}
f^{\gi j}_{\gi i}
=f^{\gi k}_{\gi q}I^{\gi i}_{\gi p}
C^{\gi j}_{\gi{ki}}
}
{Re A antilinear mapping of algebra A 1}

\DefEquation
{
\delta^{\gi i}_{\gi p}
f^{\gi j}_{\gi i}
=f^{\gi k}_{\gi 0}I^{\gi i}_{\gi p}
C^{\gi j}_{\gi{ki}}
}
{Re A antilinear mapping of algebra A 4}

\DefEquation
{
f(c)
=c^{\gi l}f^{\gi k}_{\gi 0}I^{\gi i}_{\gi l}
C^{\gi j}_{\gi{ki}}\Vector e_{\gi j}
=f^{\gi k}_{\gi 0}c^{\gi l}I^{\gi i}_{\gi l}
C^{\gi j}_{\gi{ki}}\Vector e_{\gi j}
=f^{\gi k}_{\gi 0}(c^*)^{\gi i}
C^{\gi j}_{\gi{ki}}\Vector e_{\gi j}
}
{A antilinear mapping of algebra A 6}

\DefEquation
{
C^{\gi i}_{\gi{p0}}
f^{\gi j}_{\gi i}
=f^{\gi k}_{\gi 0}I^{\gi i}_{\gi p}
C^{\gi j}_{\gi{ki}}
}
{Re A antilinear mapping of algebra A 3}

\DefEquation
{
\delta^{\gi i}_{\gi k}f^{\gi j}_{\gi i}
=f^{\gi i}_{\gi 0}C^{\gi j}_{\gi{ki}}
}
{Re A linear mapping of algebra A 4}

\DefEq
{
\[
g:A^m\rightarrow A
\]
}
{m linear map algebra}

\DefEquation
{
f(c)
=c^{\gi l}f^{\gi j}_{\gi 0}C^{\gi i}_{\gi{lj}}\Vector e_{\gi i}
}
{A linear mapping of algebra A 6}

\DefEq
{
\[
[d,a]=0
\]
}
{representation commutative with product, 3}

\DefEquation
{
\begin{matrix}
da=ad&d\in D&a\in A
\end{matrix}
}
{representation commutative with product, 2}

\DefEq
{
$D\subseteq N(A)$
}
{representation commutative with product, 5}

\DefEq
{
\begin{align*}
(d,a,b)&=d(1,a,b)=d((1a)b-1(ab))=d(ab-ab)=0
\\
(a,d,b)&=d(a,1,b)=d((a1)b-a(1b))=d(ab-ab)=0
\\
(a,b,d)&=d(a,b,1)=d((ab)1-a(b1))=d(ab-ab)=0
\end{align*}
}
{representation commutative with product, 4}

\DefEquation
{
\begin{matrix}
df(ab)=f(da,b)=f(a,db)&d\in D&a,b\in A
\end{matrix}
}
{representation commutative with product, 1}

\DefEquation
{
\begin{matrix}
d(ab)=(da)b=a(db)&d\in D&a,b\in A
\end{matrix}
}
{representation commutative with product}

\DefEq
{
\[
h(x_1,...,\widehat{x_i},...,x_n,y_1,...,y_m)=
f\circ(x_1,...,g\circ(y_1,...,y_m),...,x_n)
\]
}
{m+n-1 linear map algebra, value}

\DefEq
{
\[
f\circ(x_1,...,x_n)
\]
}
{n linear map algebra, value}

\DefEq
{
\[
x_i=g\circ(y_1,...,y_m)
\]
}
{m linear map algebra, value}

\DefEq
{
\[
f:A^n\rightarrow A
\]
}
{n linear map algebra}

\DefEquation
{
\begin{matrix}
f_{1,2}:d\circ a=da
&d\in D&a\in A
\end{matrix}
}
{representation D in A}

\DefEq
{
\begin{align*}
\re&:A\rightarrow A
\\
\im&:A\rightarrow A
\end{align*}
}
{Re Im A->A}

%% file: Biblio.English.tex
\OpenBiblio


\BiblioItem{Einstein: Electrodynamics of Moving Bodies}
{
Albert Einstein,
On the Electrodynamics of Moving Bodies, 1905,
\\
The Principle of Relativity: A Collection of Original
Memoirs on the Special and General Theory of Relativity , 37 - 65,
\\
Courier Dover Publications, 1952; ISBN-13: 978-0486600819
\\
Zur Elektrodynamik der bewegter K\"orper. Ann. Phys., 1905, 17, 891-921. 
}%

\BiblioItem{Einstein: On the Relativity Principle}
{
Albert Einstein,
On the Relativity Principle and the Conclusions Drawn from It, 1907,
\\
The Collected Papers of Albert Einstein, Volume 2:
The Swiss Years: Writings, 1900-1909. English translation. 252 - 311.
\\
Anna Beck, translator; Peter Havas, consultant.
Princeton University Press, 1989; ISBN-13: 9780691085494
\\
\"Uber das Relativit\"atsprinzip und die aus demselben gezogenen Folgerungen. 
Jahrb. d. Radioaktivit\"at u. Elektronik, 1907, 4, 411-462. 
}%

\BiblioItem{Einstein: Foundations of general relativity}
{
Albert Einstein,
Die Grundlage der allgemeinen Relativit\"atstheorie,
Ann. Phys., 1916, {\bf 49}, 769 - 822,\\
Einstein's Annalen Papers: The Complete Collection 1901-1922,
edited by J\"urgen Renn, 517 - 571,\\
Wiley-VCH Verlag GmbH \& Co. KGaA, 2005
}%

\BiblioItem{Einstein: Geometry and Experience}
{
Albert Einstein, Geometry and Experience, (1921)\\
Albert Einstein, Sidelights on Relativity, 25 - 56,\\
Courier Dover Publications, 1983
}%

\BiblioItem{Einstein: Main problems of general relativity}
{
Albert Einstein,
Grundgedanken und Probleme der Relativit\"atstheorie, (1923),\\
Nobelstiftelsen, Les Prix Nobel en 1921 - 1922,
Imprimerie Royale, Stockholm, 1923
}%

\BiblioItem{Einstein: Noneuclidean Geometry and Physics}
{
Albert Einstein,
Nichtenklidische Geometrie in der Physik Neue Rundschan, (1925)
Berlin, S. 16 - 20
}%

\BiblioItem{Einstein: Isaak Newton}
{
Albert Einstein,
Isaak Newton, 1927,
Out of My Later Years, 
Citadel Press, 1995, 219 - 223
}%

\BiblioItem{Einstein: On Science}
{
Albert Einstein,
On Science, 
Cosmic Religion, with Other Opinions and Aphorisms,142 - 146,
New York, 1931, 97 - 103
}%

\BiblioItem{Einstein: Autobiographical Notes}
{
Albert Einstein,
Autobiographical Notes, 1949,\\
Paul A. Schilpp, editor, Albert Einstein: Philosopher-Scientist,
Evanston, 
Illinois, The Library of Living Philosophers, 1949, 1 - 95
}%

\BiblioItem{Cite: 104}
{
Cite 104, Source unknown
}%

\BiblioItem{Ghez}
{
Ghez et al.,
The First Measurement of Spectral Lines in a Short-Period Star Bound to the Galaxy's Central Black Hole: A Paradox of Youth,
\href{http://www.journals.uchicago.edu/ApJ/journal/issues/ApJL/v586n2/16990/brief/16990.abstract.html}{ApJL, 586, L127} (2003),
eprint \href{http://arxiv.org/abs/astro-ph/0302299}{arXiv:astro-ph/0302299} (2003)
}%

\BiblioItem{Schodel}
{
R. Sch\"odel et al.,
A star in a 15.2-year orbit around the supermassive black hole at the centre of the Milky Way,
\href{http://www.nature.com/cgi-taf/DynaPage.taf?file=/nature/journal/v419/n6908/abs/nature01121_fs.html}{Nature 419, 694} (2002)
}%

\BiblioItem{Mielke}
{
Eckehard W. Mielke, Affine generalization of the Komar complex of general relativity,
\href{http://prola.aps.org/searchabstract/PRD/v63/i4/e044018}{Phys. Rev. D 63, 044018} (2001)
}%

\BiblioItem{Obukhov}
{
Yu. N. Obukhov and J. G. Pereira, Metric\hyph affine approach to teleparallel gravity,
\href{http://scitation.aip.org/getabs/servlet/GetabsServlet?prog=normal&id=PRVDAQ000067000004044016000001&idtype=cvips&gifs=Yes}
{Phys. Rev. D 67, 044016} (2003),
eprint \href{http://arxiv.org/abs/gr-qc/0212080}{arXiv:gr-qc/0212080} (2002)
}%

\BiblioItem{Sardanashvily}
{
Giovanni Giachetta, Gennadi Sardanashvily, Dirac Equation in Gauge and Affine-Metric Gravitation Theories,
eprint \href{http://arxiv.org/abs/gr-qc/9511035}{arXiv:gr-qc/9511035} (1995)
}%

\BiblioItem{Gauge}
{
Frank Gronwald and Friedrich W. Hehl, On the Gauge Aspects of Gravity, eprint
\href{http://arxiv.org/abs/gr-qc/9602013}{arXiv:gr-qc/9602013} (1996)
}%

\BiblioItem{Neeman}
{
Yuval Neeman, Friedrich W. Hehl, Test Matter in a Spacetime with Nonmetricity, eprint
\href{http://arxiv.org/abs/gr-qc/9604047}{arXiv:gr-qc/9604047} (1996)
}%

\BiblioItem{torsion}
{
F. W. Hehl, P. von der Heyde, G. D. Kerlick, and J. M. Nester,
General relativity with spin and torsion: Foundations and prospects,\\
\href{http://prola.aps.org/abstract/RMP/v48/i3/p393_1}{Rev. Mod. Phys. 48, 393} (1976)
}%

\BiblioItem{Megged}
{
O. Megged, Post-Riemannian Merger of Yang-Mills Interactions with Gravity,
eprint \href{http://arxiv.org/abs/hep-th/0008135}{arXiv:hep-th/0008135} (2001)
}%


\BiblioItem{gr-qc-9604027}
{
Yu.N. Obukhov, E.J. Vlachynsky, W. Esser, R. Tresguerres and F.W. Hehl,
An exact solution of the metric\hyph affine gauge theory with dilation, shear, and spin charges,
eprint \href{http://arxiv.org/abs/gr-qc/9604027}{arXiv:gr-qc/9604027} (1996)
}%

\BiblioItem{4419-7514}
{
Mari\'an Fabian, Petr Habala, Petr H\'ajek, Vicente Montesinos, V\'aclav Zizler.
Banach Space Theory: The Basis for Linear and Nonlinear Analysis.
\\
Springer; New York, 2010; ISBN-13: 978-1441975140
}%

\BiblioItem{Weinberg I}
{
Steven Weinberg.
The Quantum Theory of Fields. Volume I. Foundations.
Cambridge university press, 1995
}%

\BiblioItem{Weinberg II}
{
Steven Weinberg.
The Quantum Theory of Fields. Volume II. Modern applications.
Cambridge university press, 1996
}%

\BiblioItem{Reinhardt}
{
Walter Greiner, Joachim Reinhardt. Field Quantization. Springer.
}%

\BiblioItem{978-3540875604}
{
Walter Greiner, Joachim Reinhardt. Quantum Electrodynamics. Springer, 2009.
}%

\BiblioItem{978-1898563020}
{
H. Robert Mills. Practical Astronomy. Woodhead Publishing, 1994. ISBN-13: 978-1898563020.
}%

\BiblioItem{Landau}
{
L. D. Landau, E. M. Lifshich, The classical theory of fields.
Oxford, New York, Pergamon Press
}%

\BiblioItem{Wheeler}
{
Ignazio Ciufolini, John Wheeler. Gravitation and Inertia.
Princeton university press.
}%

\BiblioItem{Anderson02}
{
J. D. Anderson, P. A. Laing, E. L. Lau, A. S. Liu, M. M. Nieto, and S. G. Turyshev,
Study of the anomalous acceleration of Pioneer 10 and 11,
\href{http://prola.aps.org/searchabstract/PRD/v65/i8/e082004}{Phys. Rev. D 65, 082004, 50 pp.}, (2002),
eprint \href{http://arxiv.org/abs/gr-qc/0104064}{arXiv:gr-qc/0104064} (2001)
}%

\BiblioItem{Anderson98}
{
J. D. Anderson, P. A. Laing, E. L. Lau, A. S. Liu, M. M. Nieto, and S. G. Turyshev,
Indication, from Pioneer 10/11, Galileo, and Ulysses Data, of an Apparent Anomalous, Weak, Long-Range Acceleration,
\href{http://prola.aps.org/abstract/PRL/v81/i14/p2858_1}{Phys. Rev. Lett. 81, 2858}, (1998),
eprint \href{http://arxiv.org/abs/gr-qc/9808081}{arXiv:gr-qc/9808081} (1998)
}%


\BiblioItem{H. Aslaksen}
{
H. Aslaksen.  Quaternionic determinants \textit{Math.
Intelligencer} {\bf 18}(3), pp.57-65, (1996).
}%

\BiblioItem{L. Chen: Definition of determinant}
{
L. Chen, Definition of determinant and Cramer solutions over
quaternion field, \textit{Acta Math. Sinica (N.S.)} {\bf 7},
pp.171-180, (1991).
}%

\BiblioItem{L. Chen: Inverse matrix}
{
L. Chen,
Inverse matrix and properties of double determinant over quaternion
field, \textit{Sci. China, Ser. A} {\bf 34}, pp.528-540, (1991).
}%

\BiblioItem{N. Cohen S. De Leo}
{
N. Cohen, S. De Leo, The quaternionic determinant, \textit{The Electronic Journal Linear
Algebra} {\bf 7}, pp.100-111, (2000).
}%

\BiblioItem{Dyson: Quaternion determinants}
{
F. J. Dyson, Quaternion determinants, \textit{Helvetica Phys.
Acta} {\bf 45}, pp. 289-302, (1972).
}%

\BiblioItem{Melvin Hausner}
{
Melvin Hausner,
A Vector Space Approach to Geometry,
Dover Publications, 1998
}%

\BiblioItem{Serge Lang}
{
Serge Lang,
Algebra, Springer, 2002
}%

\BiblioItem{Burris Sankappanavar}
{
S. Burris, H.P. Sankappanavar,
A Course in Universal Algebra, Springer-Verlag (March, 1982),
\\eprint
\href{http://www.math.uwaterloo.ca/~snburris/htdocs/ualg.html}
{http://www.math.uwaterloo.ca/~snburris/htdocs/ualg.html}
\\(The Millennium Edition)
}%

\BiblioItem{Shilov}
{
G. E. Shilov,
Calculus, Multivariable Functions,
Moscow, Nauka, 1972
}%

\BiblioItem{Kolmogorov Fomin}
{
A. N. Kolmogorov and S. V. Fomin,
Elements of the Theory of Functions and Functional Analysis,
Courier Dover Publication, 1999
}%

\BiblioItem{Lebedev Vorovich}
{
L. P. Lebedev, I. I. Vorovich,
Functional Analysis in Mechanics,
Springer, 2002
}%

\BiblioItem
{Rashevsky}
{
P. K. Rashevsky, Riemann Geometry and Tensor Calculus,\\
Moscow, Nauka, 1967
}%

\BiblioItem
{Kurosh: High Algebra}
{
A. G. Kurosh, High Algebra,
Moscow, Nauka, 1968
}%

\BiblioItem
{Kurosh: General Algebra}
{
A. G. Kurosh, Lectures on General Algebra,
Chelsea Pub Co, 1965 
}%

\BiblioItem
{Sabinin: Smooth Quasigroups}
{
Lev V. Sabinin, Smooth Quasigroups and Loops,
Kluwer Academic Publisher, 1999 
}%

\BiblioItem{Dubrovin Fomenko Novikov part 1}
{
B. A. Dubrovin, A. T. Fomenko, S. P. Novikov,
Modern Geometry - Methods and Applications,\\
Part 1, The Geometry of Surfaces, Transformation Groups, and Fields,\\
Translated by Robert G. Burns,\\
Springer - New York, 1992
}%

\BiblioItem{Korn}
{
Granino A. Korn, Theresa M. Korn,
Mathematical Handbook for Scientists and Engineer,
McGraw-Hill Book Company, New York, San Francisco,
Toronto, London, Sydney, 1968
}%

\BiblioItem{Hocking Young Topology}
{
John G. Hocking, Gail S. Young,
Topology,\\
Courier Dover Publications, 1988
}%

\BiblioItem{Olver: Lie groups to differential equations}
{
Peter J. Olver,
Applications of Lie groups to differential equations,\\
Springer, 2000
}%

\BiblioItem{Tartaglia}
{
Angelo Tartaglia and Matteo Luca Ruggiero,
Angular Momentum Effects in Michelson\Hyph Morley Type Experiments,
Gen.Rel.Grav. 34, 1371-1382 (2002),\\
eprint \href{http://arxiv.org/abs/gr-qc/0110015}{arXiv:gr-qc/0110015} (2001)
}%

\BiblioItem{Tomozawa}
{
Yukio Tomozawa, Speed of Light in Gravitational Fields, eprint
\href{http://arxiv.org/abs/astro-ph/0303047}{arXiv:astro-ph/0303047} (2004)
}%

\BiblioItem{Magueijo}
{
Joao Magueijo,
Covariant and locally Lorentz-invariant varying speed of light theories,
\href{http://prola.aps.org/abstract/PRD/v62/i10/e103521}{Phys. Rev. D 62, 103521} (2000),
eprint \href{http://arxiv.org/abs/gr-qc/0007036}{arXiv:gr-qc/0007036} (2000)
}%

\BiblioItem{Bassett}
{
Bruce A. Bassett, Stefano Liberati, Carmen Molina-Paris, and Matt Visser,
Geometrodynamics of variable-speed-of-light cosmologies,
\href{http://prola.aps.org/abstract/PRD/v62/i10/e103518}{Phys. Rev. D 62}, 103518 (2000),
eprint \href{http://arxiv.org/abs/astro-ph/0001441}{arXiv:astro-ph/0001441} (2000)
}%

\BiblioItem{C.A. Deavours The Quaternion Calculus}
{
C.A. Deavours, The Quaternion Calculus, 
American Mathematical Monthly, {\bf 80} (1973), pp. 995 - 1008
}%

\BiblioItem{Straumann}
{
Lochlainn O'Raifeartaigh and Norbert Straumann,
Gauge theory: Historical origins and some modern developments,
\href{http://prola.aps.org/abstract/RMP/v72/i1/p1_1}{Rev. Mod. Phys. 72, 1} (2000)
}%

\BiblioItem{Lammerzahl}
{
Claus L\"ammerzahl, Mark P. Haugan,
On the interpretation of Michelson\Hyph Morley experiments,
{Phys. Lett. A282 223-229} (2001),\\
eprint \href{http://arxiv.org/abs/gr-qc/0103052}{arXiv:gr-qc/0103052} (2001)
}%

\BiblioItem{0305117}
{
Holger Mueller, Sven Herrmann, Claus Braxmaier, Stephan Schiller, Achim Peters.
Modern Michelson-Morley Experiment using Cryogenic Optical Resonators.
eprint \href{http://arxiv.org/abs/physics/0305117}{arXiv:physics/0305117} (2003)
\\
Phys. Rev. Lett. 91:020401, 2003
}%

\BiblioItem{0706.2031}
{
Holger Mueller, Paul Louis Stanwix, Michael Edmund Tobar,
Eugene Ivanov, Peter Wolf, Sven Herrmann, Alexander Senger,
Evgeny Kovalchuk, Achim Peters.
Relativity tests by complementary rotating Michelson-Morley experiments.
eprint \href{http://arxiv.org/abs/0706.2031}{arXiv:0706.2031 [physics.class-ph]} (2006)
\\
Phys. Rev. Lett. 99:050401, 2007
}%

\BiblioItem{1008.1205}
{
M. Nagel, K. M\"ohle, K. D\"oringshoff, S. Herrmann, A. Senger, E.V. Kovalchuk, A. Peters.
Testing Lorentz Invariance by Comparing Light Propagation in Vacuum and Matter.
eprint \href{http://arxiv.org/abs/1008.1205}{arXiv:1008.1205 [physics.ins-det]} (2010)
}%

\BiblioItem{1109.4897}
{
The OPERA Collaboration.
Measurement of the neutrino velocity with the OPERA detector in the CNGS beam.
eprint \href{http://arxiv.org/abs/1109.4897}{arXiv:1109.4897 [hep-ex]} (2011)
}%

\BiblioItem{Ranada}
{
Antonio F. Ranada,
Pioneer acceleration and variation of light speed: experimental situation,
eprint \href{http://arxiv.org/abs/gr-qc/0402120}{arXiv:gr-qc/0402120} (2004)
}%

\BiblioItem{Gelfand Minlos: rotation and Lorentz groups}
{
Izrail Moiseevich Gelfand, Robert Adolfovich Minlos,
Representations of the rotation and Lorentz groups and their applications;\\
Engl. transl. ed. H. K. Farahat; Transl. by G. Cummins and T. Boddongton;\\
Pergamon Press, 1963
}%

\BiblioItem{math.QA-0208146}
{
I. Gelfand, S. Gelfand, V. Retakh, R. Wilson,
Quasideterminants,\\
eprint \href{http://arxiv.org/abs/math.QA/0208146}{arXiv:math.QA/0208146} (2002)
}%

\BiblioItem{q-alg-9705026}
{
I.Gelfand, V.Retakh,
Quasideterminants, I,\\
eprint \href{http://arxiv.org/abs/q-alg/9705026}{arXiv:q-alg/9705026} (1997)
}%

\BiblioItem{Gelfand Retakh 1991}
{
I. Gelfand and V. Retakh, Determinants of Matrices over Noncommutative Rings, Funct.
Anal. Appl. 25 (1991), no. 2, 91-102
}%

\BiblioItem{Gelfand Retakh 1992}
{
I. Gelfand and V. Retakh, A Theory of Noncommutative Determinants and Characteristic
Functions of Graphs, Funct. Anal. Appl. 26 (1992), no. 4, 1-20
}%

\BiblioItem{hep-th-9407124}
{
I. M. Gelfand, D. Krob, A. Lascoux, B. Leclerc, V.S. Retakh and J.-Y. Thibon,
Noncommutative symmetric functions,\\
eprint \href{http://arxiv.org/abs/hep-th/9407124}{arXiv:hep-th/9407124} (1994)
}%

\BiblioItem{Naimark Shtern: Theory of group representations}
{
Mark Aronovich Naimark, Aleksandr Isaakovich Shtern,
Theory of group representations;\\
Heidelberg, 1982
}%

\BiblioItem{Barut Raczka: Theory of group representations}
{
Asim Orhan Barut; Ryszard R\c{a}czka;
Theory of group representations and applications;\\
World Scientific Publishing Co. Pre. Ltd., 1986
}%

\BiblioItem{Mihalev Pilz: concise handbook of algebra}
{
Aleksandr Vasilevich Mikhalev; G\"{u}nter Pilz;
The concise handbook of algebra;\\
Kluwer Academic Publishers, 2002
}%

\BiblioItem{Shafarevich: Basic notions of algebra}
{
I. R. Shafarevich,
Basic notions of algebra,\\
Translated from the Russian by M. Reid,\\
Springer, 2005
}%

\BiblioItem{Elsgolts: Differential Equations}
{
Lev Elsgolts,
Differential Equations and the Calculus of Variations,\\
University Press of the Pacific, 2003 
}%

\BiblioItem{Baez Huerta: algebra of grand unified theories}
{
John Baez; John Huerta;
The algebra of grand unified theories;\\
Bull. Amer. Math. Soc. {\bf 47} (2010), 483-552
}%

\BiblioItem{J. Fan: Determinants}
{
J. Fan, Determinants and multiplicative functionals
on quaternion matrices, \textit{Linear Algebra and Its
Applications} {\bf 369}, pp. 193-201, (2003).
}%

\BiblioItem{Carl Faith 1}
{
Carl Faith, Algebra: Rings, Modules and Categories I,
Springer - Verlag, Berlin - Heidelberg - New York, 1973
}%

\BiblioItem{Gilson Nimmo Ohta}
{
 C.R.Gilson, J.J.C.Nimmo, Y.Ohta, Quasideterminant solutions of a non-Abelian Hirota-Miwa
 equation, \textit{Journal of Physics A: Mathematical and Theoretical} {\bf 40}(42), pp.
 12607-12617,(2007).
}%

\BiblioItem{Haider Hassan}
{
B. Haider, M. Hassan, Quasideterminant solutions of an integrable chiral model in two
 dimensions, \textit{Journal of Physics A: Mathematical and Theoretical} {\bf 42} (35), art. no.
 355211, (2009).
}%



\BiblioItem{0702447}
{
I.I. Kyrchei, Cramer's rule for quaternion systems of linear equations,
\textit{Journal of Mathematical Sciences} {\bf 155}(6), 839-858, (2008).
 Translated from  \textit{Fundamental and Appl. Math.}
 {\bf 13}(4), pp.67-94, (2007). (in Russian)\\
eprint
\href{http://arxiv.org/abs/math/0702447}{arXiv:math.RA/0702447}
(2007)
}%

\BiblioItem{1004.4380}
{
I.I. Kyrchei, Cramer's rule for some quaternion matrix
    equations,  \textit{Applied Mathematics and Computation} {\bf 217}(5), pp.2024-2030, (2010).\\eprint
\href{http://arxiv.org/abs/1004.4380
}{arXiv:math.RA/arXiv:1004.4380 } (2010)
}%

\BiblioItem{1005.0736}
{
I.I. Kyrchei,Determinantal representations of the Moore-Penrose inverse
 over the quaternion skew field and corresponding Cramer's rules,
 \\
eprint
\href{http://arxiv.org/abs/1005.0736}{arXiv:math.RA/1005.0736}
(2010)
}%

\BiblioItem{0412.391}
{
Aleks Kleyn,
Basis Manifold,
eprint \href{http://arxiv.org/abs/math.DG/0412391}{arXiv:math.DG/0412391} (2007)
}%

\BiblioItem{0405.027}
{
Aleks Kleyn,
Reference Frame in General Relativity,\\
eprint \href{http://arxiv.org/abs/gr-qc/0405027}{arXiv:gr-qc/0405027} (2008)
}%

\BiblioItem{0405.028}
{
Aleks Kleyn, Metric\hyph Affine Manifold,\\
eprint \href{http://arxiv.org/abs/gr-qc/0405028}{arXiv:gr-qc/0405028} (2008)
}%

\BiblioItem{0612.111}
{
Aleks Kleyn,
Biring of Matrices,\\
eprint \href{http://arxiv.org/abs/math.OA/0612111}{arXiv:math.OA/0612111} (2007)
}%

\BiblioItem{0701.238}
{
Aleks Kleyn,
Lectures on Linear Algebra over Division Ring,\\
eprint \href{http://arxiv.org/abs/math.GM/0701238}{arXiv:math.GM/0701238} (2010)
}%

\BiblioItem{0702.561}
{
Aleks Kleyn,
Fibered $\mathfrak{F}$\Hyph Algebra,\\
eprint \href{http://arxiv.org/abs/math.DG/0702561}{arXiv:math.DG/0702561} (2007)
}%

\BiblioItem{math.RA-0501237}
{
Aleks Kleyn,
Vector Space Over Division Ring,\\
eprint \href{http://arxiv.org/abs/math.RA/0412391}{arXiv:math.RA/0501237} (2007)
}%

\BiblioItem{math.RA-0501237v1}
{
Aleks Kleyn,
Module Over Division Ring, version 1,\\
eprint \href{http://arxiv.org/abs/math/0501237v1}{arXiv:math.RA/0501237v1} (2005)
}%

\BiblioItem{0707.2246}
{
Aleks Kleyn,
Fibered Correspondence,\\
eprint \href{http://arxiv.org/abs/0707.2246}{arXiv:0707.2246} (2007)
}%

\BiblioItem{0803.2620}
{
Aleks Kleyn,
Morphism of \Ts Representations,\\
eprint \href{http://arxiv.org/abs/0803.2620}{arXiv:0803.2620} (2008)
}%

\BiblioItem{0803.3276}
{
Aleks Kleyn,
Lorentz Transformation and General Covariance Principle,\\
eprint \href{http://arxiv.org/abs/0803.3276}{arXiv:0803.3276} (2009)
}%

\BiblioItem{0812.4763}
{
Aleks Kleyn,
Introduction into Calculus over Division Ring,\\
eprint \href{http://arxiv.org/abs/0812.4763}{arXiv:0812.4763} (2010)
}%

\BiblioItem{0906.0135}
{
Aleks Kleyn,
Introduction into Geometry over Division Ring,\\
eprint \href{http://arxiv.org/abs/0906.0135}{arXiv:0906.0135} (2010)
}%

\BiblioItem{0909.0855}
{
Aleks Kleyn,
Quaternion Rhapsody,\\
eprint \href{http://arxiv.org/abs/0909.0855}{arXiv:0909.0855} (2010)
}%

\BiblioItem{0912.3315}
{
Aleks Kleyn,
Representation of Universal Algebra,\\
eprint \href{http://arxiv.org/abs/0912.3315}{arXiv:0912.3315} (2009)
}%

\BiblioItem{0912.4061}
{
Aleks Kleyn,
Linear Equation in Finite Dimensional Algebra,\\
eprint \href{http://arxiv.org/abs/0912.4061}{arXiv:0912.4061} (2010)
}%

\BiblioItem{1001.4852}
{
Aleks Kleyn,
The Matrix of Linear Mappings,\\
eprint \href{http://arxiv.org/abs/1001.4852}{arXiv:1001.4852} (2010)
}%

\BiblioItem{1003.1544}
{
Aleks Kleyn,
Linear Mappings of Free Algebra,\\
eprint \href{http://arxiv.org/abs/1003.1544}{arXiv:1003.1544} (2010)
}%

\BiblioItem{1006.2597}
{
Aleks Kleyn,
The G\^ateaux Derivative and Integral over Banach Algebra,\\
eprint \href{http://arxiv.org/abs/1006.2597}{arXiv:1006.2597} (2010)
}%

\BiblioItem{1011.3102}
{
Aleks Kleyn,
Polylinear Mapping of Free Algebra,\\
eprint \href{http://arxiv.org/abs/1011.3102}{arXiv:1011.3102} (2010)
}%

\BiblioItem{1104.5197}
{
Aleks Kleyn,
$C^*$-Rhapsody,\\
eprint \href{http://arxiv.org/abs/1104.5197}{arXiv:1104.5197} (2011)
}%

\BiblioItem{1105.4307}
{
Aleks Kleyn,
Algebra with Conjugation,\\
eprint \href{http://arxiv.org/abs/1105.4307}{arXiv:1105.4307} (2011)
}%

\BiblioItem{1107.1139}
{
Aleks Kleyn,
Linear Mappings of Quaternion Algebra,\\
eprint \href{http://arxiv.org/abs/1107.1139}{arXiv:1107.1139} (2011)
}%

\BiblioItem{1107.5037}
{
Aleks Kleyn,
Orthogonal Basis and Motion in Finsler Geometry,\\
eprint \href{http://arxiv.org/abs/1107.5037}{arXiv:1107.5037} (2011)
}%

\BiblioItem{1202.6021}
{
Aleks Kleyn,
Mappings of Conjugation of Quaternion Algebra,\\
eprint \href{http://arxiv.org/abs/1202.6021}{arXiv:1202.6021} (2012)
}%

\BiblioItem{8433-5163}
{
Aleks Kleyn,
Linear Mappings of Free Algebra: First Steps in Noncommutative Linear Algebra,\\
Lambert Academic Publishing, 2010
}%

\BiblioItem{8443-0072}
{
Aleks Kleyn,
Representation Theory: Representation of Universal Algebra,\\
Lambert Academic Publishing, 2011
}%

\BiblioItem{Lauve: Quantum coordinates}
{
A. Lauve, Quantum- and quasi-Plucker coordinates,
\textit{Journal of Algebra} {\bf 296}(2), pp.440-461,
(2006).
}%

\BiblioItem{Lewis D. W. Quaternion algebras}
{
Lewis D. W. Quaternion algebras and the algebraic legacy
of Hamilton's quaternions, \textit{Irish Math. Soc. Bulletin} {\bf
57}, pp. 41-64, (2006).
}%

\BiblioItem{0812.2865}
{
Jos\'e Miguel Figueroa-O'Farrill,
Three lectures on 3-algebras,
eprint \href{http://arxiv.org/abs/0812.2865}{arXiv:0812.2865} (2008)
}%

\BiblioItem{1202.4546}
{
Ming-Liang Hu,
Disentanglement, Bell-nonlocality violation
and teleportation capacity of the decaying tripartite states,
eprint \href{http://arxiv.org/abs/1202.4546}{arXiv:1202.4546} (2012)
}%

\BiblioItem{1203.1629}
{
Borivoje Dakic, Yannick Ole Lipp, Xiaosong Ma, Martin Ringbauer,
Sebastian Kropatschek, Stefanie Barz, Tomasz Paterek, Vlatko Vedral,
Anton Zeilinger, Caslav Brukner, Philip Walther,
Quantum Discord as Optimal Resource for Quantum Communication,
eprint \href{http://arxiv.org/abs/1203.1629}{arXiv:1203.1629} (2012)
}%

\BiblioItem{Li Nimmo: Darboux transformations}
{
C.X.Li, J.J.C. Nimmo, Darboux transformations for a twisted
derivation and quasideterminant solutions to the super KdV
equation, \textit{Proceedings of the Royal Society A:
Mathematical, Physical and Engineering Sciences} {\bf 466} (2120),
pp. 2471-2493, (2010)
}%

\BiblioItem{Schiebold: Cauchy-type determinants}
{
C. Schiebold, Cauchy-type determinants and integrable
systems, \textit{Linear Algebra and Its Applications} {\bf 433}
(2), pp. 447-475, (2010)
}%

\BiblioItem{Suzuki: Noncommutative spectral decomposition}
{
T. Suzuki, Noncommutative
spectral decomposition with qua\-si\-de\-ter\-mi\-nant, \textit{Advances in
Mathematics} {\bf 217}(5), pp. 2141-2158, (2008)
}%

\BiblioItem{1105.3456}
{
C. W. F. Everitt, D. B. DeBra, B. W. Parkinson, J. P. Turneaure, J. W. Conklin,
M. I. Heifetz, G. M. Keiser, A. S. Silbergleit, T. Holmes, J. Kolodziejczak,
M. Al-Meshari, J. C. Mester, B. Muhlfelder, V. Solomonik, K. Stahl, P. Worden,
W. Bencze, S. Buchman, B. Clarke, A. Al-Jadaan, H. Al-Jibreen, J. Li, J. A. Lipa,
J. M. Lockhart, B. Al-Suwaidan, M. Taber, S. Wang,\\
Gravity Probe B: Final Results of a Space Experiment to Test General Relativity,\\
eprint \href{http://arxiv.org/abs/1105.3456}{arXiv:1105.3456[gr-qc]} (2011)
}%

\BiblioItem{0009305}
{
G. S. Asanov.
Can Neutrinos and High-Energy Particles Test Finsler Metric of Space-Time?\\
eprint \href{http://arxiv.org/abs/hep-ph/0009305}{arXiv:hep-ph/0009305} (2000)
}%

\BiblioItem{Asanov 2004}
{
G. S. Asanov.
Finsleroid - space supplemented by angle and scalar product.\\
Hypercomplex Numbers in Geometry and Physics, {\bf 1}, 2004, p. 40 - 62
}%

\BiblioItem{1004.3007}
{
Sergiu I. Vacaru,
Principles of Einstein-Finsler Gravity and Perspectives in Modern Cosmology,\\
eprint \href{http://arxiv.org/abs/1004.3007}{arXiv:1004.3007[math-ph]} (2010)
}%

\BiblioItem{1012.4148}
{
Sergiu I. Vacaru.
Principles of Einstein-Finsler Gravity and Cosmology.\\
eprint \href{http://arxiv.org/abs/1012.4148}{arXiv:1012.4148[physics.gen-ph]} (2010)
}%

\BiblioItem{1112.5641}
{
Christian Pfeifer, Mattias N.R. Wohlfarth.
Finsler geometric extension of Einstein gravity.\\
eprint \href{http://arxiv.org/abs/1112.5641}{arXiv:1112.5641[gr-qc]} (2011)
}%

\BiblioItem{0711.0056}
{
Zhe Chang, Xin Li.
Lorentz Invariance Violation and Symmetry in Randers\Hyph Finsler Spaces.\\
eprint \href{http://arxiv.org/abs/0711.0056}{arXiv:0711.0056[hep-th]} (2011)
}%

\BiblioItem{Rund Finsler geometry}
{
Hanno Rund,
The differential geometry of Finsler spaces.
\\
Springer - Verlag, Berlin - G\"ottingen - Heidelberg, 1959
}%

\BiblioItem{Smirnov vol 1}
{
V. I. Smirnov,
A Course of Higher Mathematics, volume I.
\\
Translated by D. E. Brown.
\\
Translation, edited and additions made by I. N. Sneddon.
\\
Pergamon Press, Addison-Wesley Publishing Company, 1964
}%

\BiblioItem{Beem Dostoglou Ehrlich}
{
John K. Beem, Stamatis A. Dostoglou, Paul E. Ehrlich,
Advances in differential geometry and general relativity.
\\
American Mathematical Society, 2004
}%

\BiblioItem{978-0719033414}
{
Malcolm Pemberton, Nicholas Rau,
Mathematics for economists: an introductory textbook.
\\
Manchester University Press, November 2001; ISBN-13: 978-0719033414
}%

\BiblioItem{0 521 59180 5}
{
Cyrus D. Cantrell,
Modern mathematical methods for physicists and engineers.
\\
Cambridge University Press, 2000
}%

\BiblioItem{Arveson spectral theory}
{
William Arveson,
A short course on spectral theory.
\\
Springer - Verlag, New York, 2002
}%

\BiblioItem{Robert Hermann}
{
Robert Hermann,
Topics in the mathematics of quantum mechanics.
\\
Math Sci Press, 1973
}%

\BiblioItem{9705.009}
{
John C. Baez,
An Introduction to n-Categories,\\
eprint \href{http://arxiv.org/abs/q-alg/9705009}{arXiv:q-alg/9705009} (1997)
}%

\BiblioItem{0105.155}
{
John C. Baez,
The Octonions,\\
eprint \href{http://arxiv.org/abs/math.RA/0105155}{arXiv:math.RA/0105155} (2002)
}%

\BiblioItem{John Baez: Math Blogs}
{
John C. Baez,
What do mathematicians need to know about blogging?,\\
Notices of the American Mathematical Society,
(2010), 3, {\bf 57}, 333,\\
\url{http://www.ams.org/notices/201003/rtx100300333p.pdf}
}%

\BiblioItem{Tolstoi about Anna Karenina}
{
Tolstoi about Anna Karenina,
in book A Karenina Companion, by C. J. G. Turner,
published by Wilfrid Laurier University Press (August 1993)
}%

\BiblioItem
{Cohn: Universal Algebra}
{
Paul M. Cohn,
Universal Algebra,
Springer, 1981
}%

\BiblioItem
{Maunder: Algebraic Topology}
{
C. R. F. Maunder,
Algebraic Topology,
Dover Publications, Inc, Mineola, New York, 1996
}%

\BiblioItem{Pommaret: Partial Differential Equations}
{
J.-F. Pommaret,
Partial Differential Equations and Group Theory,
Springer, 1994
}%

\BiblioItem{Bourbaki: Set Theory}
{
N. Bourbaki,
Theory of sets,
Springer, 2004
}%

\BiblioItem{Bourbaki: Algebra 1}
{
N. Bourbaki,
Algebra 1,
Springer, 2004
}%

\BiblioItem
{Bourbaki: General Topology 1}
{
N. Bourbaki,
General Topology, Chapters 1 - 4,
Springer, 1989
}

\BiblioItem{Bourbaki: General Topology: Chapter 5 - 10}
{
N. Bourbaki,
General Topology, Chapters 5 - 10,
Springer, 1989
}

\BiblioItem{Bourbaki: Topological Vector Space}
{
N. Bourbaki,
Topological Vector Spaces, Chapters 1 - 5,
Transl. by H. G. Eggleston $\&$ S. Madan,
Springer, 2003
}

\BiblioItem{Bourbaki: Coxeter Group Lie}
{
N. Bourbaki,
Lie Groups and Lie Algebras, Chapters 4 - 6,
Translator Andrew Pressley,
Springer, 2002
}

\BiblioItem{Bourbaki: Real Group Lie}
{
N. Bourbaki,
Lie Groups and Lie Algebras, Chapters 7 - 9,
Translator Andrew Pressley,
Springer, 2005
}

\BiblioItem{Shabat: Complex Analysis}
{
Shabat B. V.,
Introduction to Complex Analysis,
\\ \url{http://www.math.uchicago.edu/~ryzhik/shabat-all.pdf},
\\Translated from Russian by L.Ryzhik, 2003
(Moscow, Nauka, 1969)
}

\BiblioItem{Pontryagin: Topological Group}
{
L. S. Pontryagin,
Selected Works, Volume Two, Topological Groups,
Gordon and Breach Science Publishers, 1986
}

\BiblioItem
{Eisenhart: Riemannian Geometry}
{
Eisenhart,
Riemannian Geometry,
Princeton University Press, Princeton, 1949
}

\BiblioItem
{Eisenhart: Continuous Groups of Transformations}
{
Eisenhart,
Continuous Groups of Transformations,
Dover Publications, New York, 1961
}

\BiblioItem
{Condon Odabasi}
{
Edward Uhler Condon, Halis Odabasi,
Atomic Structure,
CUP Archive, 1980
}

\BiblioItem{Postnikov: Differential Geometry}
{
Postnikov M. M.,
Geometry IV: Differential geometry,
Moscow, Nauka, 1983
}

\BiblioItem{Fihtengolts: Calculus volume 1}
{
Fihtengolts G. M.,
Differential and Integral Calculus Course, volume 1,
Moscow, Nauka, 1969
}

\BiblioItem{Hatcher: Algebraic Topology}
{
Allen Hatcher,
Algebraic Topology,
Cambridge University Press, 2002
}

\BiblioItem{geometry of differential equations}
{
Vinogradov, A. M., Krasil'shchik, I. S., and Lychagin, V. V.,
Introduction to geometry of nonlinear differential equations,
Nauka, Moscow, 1986
}

\BiblioItem{cohomological analysis}
{
A. M. Vinogradov,
Cohomological Analysis of Partial Differential Equations
and Secondary Calculus,
American Mathematical Society, 2001
}

\BiblioItem{0801.1734}
{
Brandon S. DiNunno, Richard A. Matzner,
The Volume Inside a Black Hole,\\
eprint \href{http://arxiv.org/abs/0801.1734v1}{arXiv:0801.1734v1} (2008)
}

\BiblioItem{0702.447}
{
Ivan Kyrchei,
Cramer's rule for some quaternion matrix equations,\\
eprint \href{http://arxiv.org/abs/math/0702447}{arXiv:math.RA/0702447} (2007)
}

\BiblioItem{Izrail M. Gelfand: Quaternion Groups}
{
I. M. Gelfand, M. I. Graev,
Representation of Quaternion Groups over Localy Compact and
Functional Fields,\\
Funct. Anal. Appl. {\bf 2} (1968) 19 - 33;\\
Izrail Moiseevich Gelfand, Semen Grigorevich Gindikin,\\
Izrail M. Gelfand: Collected Papers, volume II, 435 - 449,\\
Springer, 1989
}

\BiblioItem{Richard D. Schafer}
{
Richard D. Schafer,
An Introduction to Nonassociative Algebras,
Dover Publications, Inc., New York, 1995
}

\BiblioItem{Bamberg Sternberg}
{
Paul Bamberg, Shlomo Sternberg,
A course in mathematics for students of physics,
Cambridge University Press, 1991
}

\BiblioItem{Conway Smith}
{
John Horton Conway, Derek Alan Smith,
On quaternions and octonions: their geometry, arithmetic, and symmetry,
A K Peters, Natick, Massachussets, 2003
}

\BiblioItem{Fueter}
{
Fueter, R.
Die Funktionentheorie der Differentialgleichungen $\Delta u = 0$ und
$\Delta \Delta u = 0$ mit vier reellen Variablen.
Comment. Math. Helv. {\bf 7} (1935), 307-330
}

\BiblioItem{Sudbery Quaternionic Analysis}
{
A. Sudbery,
Quaternionic Analysis,
Math. Proc. Camb. Phil. Soc. (1979), {\bf 85}, 199 - 225
}

\BiblioItem{0902.4771}
{
Fabrizio Colombo, Graziano Gentili, Irene Sabadini,
A Cauchy kernel for slice regular functions,\\
eprint \href{http://arxiv.org/abs/0902.4771v1}{arXiv:0902.4771v1} (2009)
}

\BiblioItem{Vadim Komkov}
{
Vadim Komkov,
Variational Principles of Continuum Mechanics with Engineering Applications: Critical Points Theory,\\
Springer, 1986
}

\BiblioItem{Alain Connes 1994}
{
Alain Connes,
Noncommutative Geometry,\\
Academic Press, 1994
}

\BiblioItem{Hamilton papers 3}
{
Sir William Rowan Hamilton,
The Mathematical Papers, Vol. III, Algebra,\\
Cambridge at the University Press, 1967
}

\BiblioItem{Hamilton Elements of Quaternions 1}
{
Sir William Rowan Hamilton,
Elements of Quaternions, Volume I,\\
Longmans, Green, and Co., London, New York, and Bombay, 1899
}

\BiblioItem{Cartan geometry in reper}
{
Elie Cartan, Vladislav V. Goldberg, Serge\u{i} Pavlovich Finikov,\\
Riemannian geometry in an orthogonal frame:
from lectures delivered by Elie Cartan at the Sorbonne in 1926-1927,\\
translated by Vladislav V. Goldberg,\\
World Scientific, 2001
}

\BiblioItem{Arnautov Glavatsky Mikhalev}
{
V. I. Arnautov, S. T. Glavatsky, A. V. Mikhalev,\\
Introduction to the theory of topological rings and modules,
Volume 1995,\\
Marcel Dekker, Inc, 1996
}

\BiblioItem{Moore Yaqub}
{
Hal G. Moore, Adil Yaqub,
A first course in linear algebra with applications,
Edition 3, Academic Press, 1998 
}

\BiblioItem{math.CV-0405471}
{
S. V. Ludkovsky,
Differentiable functions of Cayley-Dickson numbers,\\
eprint \href{http://arxiv.org/abs/math.CV/0405471}{arXiv:math.CV/0405471} (2004)
}%

\BiblioItem{W.Bertram H.Glockner K.Neeb}
{
W.Bertram, H.Glockner, K.Neeb,
Differential Calculus over General Base Fields and Rings,
Expositiones Mathematicae (2004), Volume 22, Issue 3, Pages 213-282
}

\CloseBiblio

%% file: Index.English.tex
\OpenIndex
\SetIndexSpace%
\Index
   {$1$-\rcd form}%
   {1-rcd form, vector spaces}%
\SetIndexSpace%
\Index
   {$2$\Hyph ary fibered relation}%
   {2 ary fibered relation}%
\SetIndexSpace%
\Index
   {$(^a_b)$\hyph \CR quasideterminant}%
   {a b cr-quasideterminant}%
\Index
   {$A\CRcirc$\Hyph basis for module}%
   {A CRcirc basis, module over algebra}%
\Index
   {$A\CRcirc$\Hyph linearly dependent set of vectors}%
   {CRcirc linearly dependent, Astar module over algebra}%
\Index
   {$A\CRcirc$\Hyph linearly independent set of vectors}%
   {CRcirc linearly independent, Astar module over algebra}%
\Index
   {$A$\Hyph linear mapping of modules}%
   {A linear map of modules}%
\Index
   {$\mathcal A(A)$\Hyph mapping}%
   {A(A) mapping}%
\Index
   {$A$\Hyph module}%
   {module over algebra}%
\Index
   {$A$\Hyph valued function}%
   {A valued function}%
\Index
   {Abelian $\Omega$\Hyph group}%
   {Abelian Omega group}%
\Index
   {absolute value on division ring}%
   {absolute value on division ring}%
\Index
   {absolute value on ring}%
   {absolute value on ring}%
\Index
   {\Acr linear mapping of modules}%
   {Acr linear map of modules}%
\Index
   {$A\CRcirc$\Hyph linear combination}%
   {ACRcirc linear combination}%
\Index
   {active representation}%
   {active representation}%
\Index
   {active representation of group $G(f)$ in basis manifold of representation}%
   {active representation in basis manifold}%
\Index
   {active representation of group $G(\Vector f)$ in basis manifold of tower of representations}%
   {active representation in basis manifold, tower of representations}%
\Index
   {active \sT representation}%
   {active representation, vector space}%
\Index
   {active transformation of basis manifold of representation}%
   {active transformation of basis, representation}%
\Index
   {active transformation of basis manifold of tower of representations}%
   {active transformation of basis, tower of representations}%
\Index
   {active transformation on basis manifold}%
   {active transformation}%
\Index
   {active transformation on the set of \rcd bases}%
   {active transformation, vector space}%
\Index
   {affine basis}%
   {Affine Basis}%
\Index
   {affine structure on set}%
   {affine structure on set}%
\Index
   {affine transformation group}%
   {drc affine transformation group}%
\Index
   {affine transformation group}%
   {affine transformation group}%
\Index
   {affine transformation on basis manifold}%
   {affine transformation}%
\Index
   {algebra over ring}%
   {algebra over ring}%
\Index
   {algebra with conjugation}%
   {algebra with conjugation}%
\Index
   {alternative representation of matrix}%
   {Alternative representation}%
\Index
   {anholonomic coordinate}%
   {anholonomic coordinate}%
\Index
   {anholonomic coordinates of connection}%
   {anholonomic coordinates of connection}%
\Index
   {anholonomic coordinates of vector}%
   {vector anholonomic coordinates}%
\Index
   {anholonomic coordinates on manifold}%
   {anholonomic coordinates on manifold}%
\Index
   {anholonomity object}%
   {anholonomity object}%
\Index
   {antihomomorphism of fibered groups}%
   {antihomomorphism of fibered groups}%
\Index
   {antisymmetric $2$\Hyph ary fibered relation}%
   {antisymmetric 2 ary fibered relation}%
\Index
   {$A\RCstar$\Hyph basis for vector space}%
   {Arc basis, vector space}%
\Index
   {$A\RCstar$\Hyph linearly dependent vectors}%
   {linearly dependent, A vector space}%
\Index
   {$A\RCstar$\Hyph linearly independent vectors}%
   {linearly independent, A vector space}%
\Index
   {arity of operation}%
   {arity of operation}%
\Index
   {associative $D$\Hyph algebra}%
   {associative D algebra}%
\Index
   {associative law for $A\star$\Hyph linear mappings of vector spaces}%
   {associative law for Astar linear mappings of vector spaces}%
\Index
   {associative law for $A\star$\Hyph module}%
   {associative law, Astar module over algebra}%
\Index
   {associative law for $A\star$\Hyph vector space}%
   {associative law, Astar vector space}%
\Index
   {associative law for covariant \sT representation}%
   {associative law for covariant starT representation}%
\Index
   {associative law for covariant \Ts representation}%
   {associative law for covariant Tstar representation}%
\Index
   {associative law for $D$\Hyph module}%
   {associative law, D module}%
\Index
   {associative law for \Drc linear maps of vector bundles}%
   {associative law for drc linear maps of vector bundles}%
\Index
   {associative law for $\mathcal D\star$\Hyph vector fields}%
   {associative law, Dstar vector fields}%
\Index
   {associative law for $D\star$\Hyph vector space}%
   {associative law, Dstar vector space}%
\Index
   {associative law for \rcd linear maps of vector spaces}%
   {associative law for rcd linear maps of vector spaces}%
\Index
   {associative law for $\star A$\Hyph module}%
   {associative law, starA module over algebra}%
\Index
   {associative law for twin representations}%
   {associative law for twin representations}%
\Index
   {associative law of composition of fibered correspondences}%
   {associative law, composition of fibered correspondences}%
\Index
   {associator of $R$\Hyph algebra}%
   {associator of algebra}%
\Index
   {$A\star$\Hyph antilinear mapping of algebra with conjugation}%
   {antilinear mapping of algebra with conjugation}%
\Index
   {$A\star$\Hyph linear map of vector spaces}%
   {Astar linear map of vector spaces}%
\Index
   {$A\star$\Hyph vector space}%
   {Astar vector space}%
\Index
   {$A\star$\Hyph module}%
   {Astar-module}%
\Index
   {$A\star$\Hyph product of vector over scalar}%
   {Astar product of vector over scalar, Astar module}%
\Index
   {$A\star$\hyph product of vector over scalar}%
   {Astar product of vector over scalar, vector space}%
\Index
   {auto parallel line}%
   {auto parallel line}%
\Index
   {automorphism of representation of $\Omega$\Hyph algebra}%
   {automorphism of representation}%
\Index
   {automorphism of tower of representations}%
   {automorphism of tower of representations}%
\Index
   {norm of quaternion}%
   {norm of quaternion}%
\SetIndexSpace%
\Index
   {Banach $D$\Hyph algebra}%
   {Banach algebra}%
\Index
   {base of fibered correspondence}%
   {base of fibered correspondence}%
\Index
   {base of mapping}%
   {base of map}%
\Index
   {basis manifold of affine space}%
   {Basis Manifold, Affine Space}%
\Index
   {basis manifold of central affine space}%
   {Basis Manifold, Central Affine Space, division ring}%
\Index
   {basis manifold of central affine space}%
   {Basis Manifold, Central Affine Space}%
\Index
   {basis manifold of Euclid space}%
   {Basis Manifold, Euclid Space}%
\Index
   {basis manifold of Euclid space}%
   {Basis Manifold, Euclid Space, division ring}%
\Index
   {basis manifold of \rcd affine space}%
   {Basis Manifold, rcd Affine Space, division ring}%
\Index
   {basis manifold of \rcd vector space}%
   {basis manifold of rcd vector space}%
\Index
   {basis manifold of representation}%
   {basis manifold representation F algebra}%
\Index
   {basis manifold of tower of representations}%
   {basis manifold tower of representations}%
\Index
   {basis manifold of vector space}%
   {basis manifold of vector space}%
\Index
   {basis of $A$\Hyph module}%
   {basis of A module}%
\Index
   {basis of $A\CRcirc$\Hyph module}%
   {basis of ACRcirc module}%
\Index
   {basis of algebra $\mathcal L(A;A)$}%
   {basis of algebra L(A,A)}%
\Index
   {basis of $\RCcirc A$\Hyph module}%
   {basis of CRcircA module}%
\Index
   {basis of representation}%
   {basis of representation}%
\Index
   {basis of tower of representations}%
   {basis of tower of representations}%
\Index
   {basis of vector space}%
   {Basis}%
\Index
   {basis vector of representation of Lie group}%
   {basis vector of representation of Lie group}%
\Index
   {basis vector of representation of Lie group over algebra $A$}%
   {basis vector of representation of Lie group over algebra A}%
\Index
   {biring}%
   {biring}%
\Index
   {bundle of level $2$}%
   {bundle of level 2}%
\Index
   {bundle of level $n$}%
   {bundle of level n}%
\SetIndexSpace%
\Index
   {\subs row of matrix}%
   {c row}%
\Index
   {\sups rows \rcd vector space}%
   {sups rows rcd vector space}%
\Index
   {$c$\hyph row of matrix}%
   {c-row}%
\Index
   {Cartan connection}%
   {Cartan connection}%
\Index
   {Cartan curvature}%
   {Cartan curvature}%
\Index
   {Cartan derivative}%
   {Cartan derivative}%
\Index
   {Cartan symbol}%
   {Cartan symbol}%
\Index
   {Cartan transport}%
   {Cartan transport}%
\Index
   {Cartesian power $\Bundle A$ of bundle $\Bundle B$}%
   {Cartesian power A of bundle B}%
\Index
   {Cartesian power $A$ of set $B$}%
   {Cartesian power of set}%
\Index
   {Cartesian power $n$ of bundle $\Bundle E$}%
   {Cartesian power n of bundle E}%
\Index
   {Cartesian power $n$ of $\mathfrak{H}$\Hyph algebra}%
   {Cartesian power of algebra}%
\Index
   {category of \drc vector spaces}%
   {category of drc vector spaces}%
\Index
   {category of fibered correspondences over diagonal}%
   {category of fibered correspondences over diagonal}%
\Index
   {category of reduced fibered correspondences}%
   {category of reduced fibered correspondences}%
\Index
   {category of \Ts representations of $\Omega_1$\Hyph algebra $A$}%
   {category of Tstar representations of Omega1 algebra}%
\Index
   {category of \Ts representations of $\Omega_1$\Hyph algebra from category $\mathcal A$}%
   {category of Tstar representations of Omega1 algebra from category}%
\Index
   {Cauchy sequence in normed algebra}%
   {Cauchy sequence, normed algebra}%
\Index
   {Cauchy sequence in normed ring}%
   {Cauchy sequence, normed ring}%
\Index
   {Cauchy sequence in valued division ring}%
   {Cauchy sequence, valued division ring}%
\Index
   {Cauchy sequence in valued ring}%
   {Cauchy sequence, valued ring}%
\Index
   {center of an $R$\Hyph algebra $A$}%
   {center of algebra}%
\Index
   {center of ring $D$}%
   {center of ring}%
\Index
   {central affine basis}%
   {Central Affine Basis, division ring}%
\Index
   {central affine basis}%
   {Central Affine Basis}%
\Index
   {column determinant}%
   {column determinant}%
\Index
   {column vector}%
   {column vector}%
\Index
   {commutative $D$\Hyph algebra}%
   {commutative D algebra}%
\Index
   {commutative diagram of correspondences}%
   {commutative diagram of correspondences}%
\Index
   {commutator of $R$\Hyph algebra}%
   {commutator of algebra}%
\Index
   {compact\hyph open topology}%
   {compact open topology}%
\Index
   {complete division ring}%
   {complete division ring}%
\Index
   {complete ring}%
   {complete ring}%
\Index
   {complete system of linear partial differential equations}%
   {Complete System of Linear Partial Differential Equations}%
\Index
   {completely integrable system}%
   {completely integrable system}%
\Index
   {component of the G\^ateaux derivative of map $\Vector f(\Vector x)$}%
   {component of Gateaux derivative of map, D vector space}%
\Index
   {component of the G\^ateaux derivative of map $f(x)$}%
   {component of Gateaux derivative of map, division ring}%
\Index
   {component of the G\^ateaux derivative of second order of map  of division ring}%
   {component of Gateaux derivative of Second Order, division ring}%
\Index
   {component of the G\^ateaux derivative of second order of map $\Vector f(\Vector x)$}%
   {component of Gateaux derivative of Second Order, D vector space}%
\Index
   {component of linear map $f$ of division ring}%
   {component of linear map, division ring}%
\Index
   {component of linear map of $D$\Hyph vector space}%
   {component of linear map, D vector space}%
\Index
   {component of polylinear map into associative algebra}%
   {component of polylinear map, associative algebra}%
\Index
   {component of polylinear map of division ring}%
   {component of polylinear map, division ring}%
\Index
   {component of polylinear mapping $\Vector A$}%
   {component of polyadditive map, D vector space}%
\Index
   {component of the G\^ateaux derivative of map $f(x)$ of algebra}%
   {component of Gateaux derivative of map, algebra}%
\Index
   {component of the G\^ateaux derivative of second order of map $f(x)$ of algebra}%
   {component of Gateaux derivative of Second Order, algebra}%
\Index
   {composition of fibered correspondences}%
   {composition of fibered correspondences}%
\Index
   {composition of reduced fibered correspondences}%
   {composition of reduced fibered correspondences}%
\Index
   {condition of reducibility of products}%
   {condition of reducibility of products}%
\Index
   {conjugate of quaternion $x$}%
   {conjugate of quaternion}%
\Index
   {conjugation in algebra}%
   {conjugation in algebra}%
\Index
   {conjugation in ring}%
   {conjugation in ring}%
\Index
   {connection coefficients in $D$\Hyph affine space}%
   {connection coefficients, D affine space}%
\Index
   {continuous correspondence}%
   {continuous correspondence}%
\Index
   {continuous function of division ring}%
   {continuous function, division ring}%
\Index
   {continuous function over $D$\Hyph algebra}%
   {continuous function, algebra}%
\Index
   {contravariant \sT representation of fibered group}%
   {contravariant starT representation of fibered group}%
\Index
   {contravariant \sT representation of group}%
   {contravariant starT representation of group}%
\Index
   {contravariant \Ts representation of fibered group}%
   {contravariant Tstar representation of fibered group}%
\Index
   {contravariant \Ts representation of group}%
   {contravariant Tstar representation of group}%
\Index
   {coordinate \Drc vector bundle}%
   {coordinate drc vector bundle}%
\Index
   {coordinate isomorphism}%
   {coordinate isomorphism}%
\Index
   {coordinate matrix of set of vectors in \dcr rows vector space}%
   {coordinate matrix of set of vectors, dcr vector space}%
\Index
   {coordinate matrix of set of vectors in \rcd rows vector space}%
   {coordinate matrix of set of vectors, rcd vector space}%
\Index
   {coordinate matrix of vector field in \rcD basis}%
   {coordinate matrix of vector field in drc basis}%
\Index
   {coordinate matrix of vector in basis}%
   {coordinate matrix of vector in basis, module}%
\Index
   {coordinate matrix of vector in \drc basis}%
   {coordinate matrix of vector in rcd basis}%
\Index
   {coordinate \rcd isomorphism}%
   {coordinate rcd isomorphism}%
\Index
   {coordinate \rcd vector space}%
   {coordinate rcd vector space}%
\Index
   {coordinate reference frame}%
   {coordinate reference frame}%
\Index
   {coordinate representation in $\Omega_2$\Hyph algebra}%
   {coordinate representation, Omega_2 algebra}%
\Index
   {coordinate representation in \rcd vector space}%
   {coordinate representation, rcd vector space}%
\Index
   {coordinate representation in tuple of $\VX\Omega$\Hyph algebras}%
   {coordinate tower of representations, Omega algebra}%
\Index
   {coordinate representation of group in vector space}%
   {coordinate representation, vector space}%
\Index
   {coordinate vector space}%
   {coordinate vector space}%
\Index
   {coordinates of a geometric object in $\Omega_2$\Hyph algebra $M$}%
   {coordinates of geometric object, representation g}%
\Index
   {coordinates of a geometric object in tuple of $\VX\Omega$\Hyph algebras}%
   {coordinates of geometric object, tower of representations g}%
\Index
   {coordinates of basis of representation}%
   {coordinates of basis relative to basis, representation}%
\Index
   {coordinates of element $m$ of representation $f$ relative to set $X$}%
   {coordinates of element relative to set, representation}%
\Index
   {coordinates of endomorphism of representation}%
   {coordinates of endomorphism, representation}%
\Index
   {coordinates of endomorphism of tower of representations}%
   {coordinates of endomorphism, tower of representations}%
\Index
   {coordinates of geometric object}%
   {coordinates of geometric object, vector space}%
\Index
   {coordinates of geometric object in coordinate \rcd vector space}%
   {coordinates of geometric object, coordinate rcd vector space}%
\Index
   {coordinates of geometric object in coordinate representation}%
   {coordinates of geometric object, coordinate vector space}%
\Index
   {coordinates of geometric object in coordinate space of representation}%
   {coordinates of geometric object, coordinate representation}%
\Index
   {coordinates of geometric object in coordinate space of tower of representations}%
   {coordinates of geometric object, coordinate tower of representations}%
\Index
   {coordinates of geometric object in \rcd vector space}%
   {coordinates of geometric object, rcd vector space}%
\Index
   {coordinates of point $A$ of affine space $\overset{\circ}{A}$ relative to basis $(O,\Basis e)$}%
   {coordinates in affine space}%
\Index
   {coordinates of point of \rcd affine space relative to basis}%
   {coordinates in rcd affine space}%
\Index
   {coordinates of representation}%
   {coordinates of representation, drc vector space}%
\Index
   {coordinates of representation}%
   {coordinates of representation}%
\Index
   {coordinates of set of vectors in \dcr vector space}%
   {coordinates of set of vectors, dcr vector space}%
\Index
   {coordinates of set of vectors in \rcd vector space}%
   {coordinates of set of vectors, rcd vector space}%
\Index
   {coordinates of vector field in \Drc basis}%
   {coordinates of vector field in drc basis}%
\Index
   {coordinates of vector in basis}%
   {coordinates of vector in basis, module}%
\Index
   {coordinates of vector in \rcd basis}%
   {coordinates of vector in rcd basis}%
\Index
   {coordinates of vector relative to Schauder basis}%
   {coordinates of vector, Schauder basis}%
\Index
   {correspondence continuous on the set}%
   {correspondence continuous on the set}%
\Index
   {correspondence of homomorphism}%
   {correspondence of homomorphism}%
\Index
   {covariant \sT representation of fibered group}%
   {covariant starT representation of fibered group}%
\Index
   {covariant \sT representation of group}%
   {covariant starT representation of group}%
\Index
   {covariant \Ts representation of fibered group}%
   {covariant Tstar representation of fibered group}%
\Index
   {covariant \Ts representation of group}%
   {covariant Tstar representation of group}%
\Index
   {\CR inverse element of biring}%
   {cr-inverse element}%
\Index
   {\CR matrix group}%
   {cr-matrix group}%
\Index
   {\CR power}%
   {cr power}%
\Index
   {\CR product of matrices}%
   {cr-product of matrices}%
\Index
   {$\CRcirc$\Hyph product of matrices of mappings}%
   {cr product of matrices of mappings}%
\Index
   {\crd vector space}%
   {crd vector space}%
\Index
   {$C^*$\Hyph algebra}%
   {Cstar-algebra}%
\Index
   {curvilinear coordinates of point in affine space}%
   {curvilinear coordinates of point in affine space}%
\Index
   {\subs rows \drc vector space}%
   {subs rows drc vector space}%
\SetIndexSpace%
\Index
   {$D$\Hyph affine space}%
   {d affine space}%
\Index
   {$D$\Hyph linearly dependent vectors of $D$\Hyph module}%
   {D linearly dependent, module}%
\Index
   {$D$\Hyph linearly independent vectors of $D$\Hyph module}%
   {D linearly independent, module}%
\Index
   {$D$\Hyph affine connection on manifold with affine connections}%
   {D affine connection, affine manifold}%
\Index
   {$D$\Hyph basis for module}%
   {D basis, module}%
\Index
   {$D$\Hyph module}%
   {D-module}%
\Index
   {$D$\Hyph valued variable}%
   {D valued variable}%
\Index
   {$D$\Hyph vector function}%
   {d vector function}%
\Index
   {$D$\Hyph affine connection coefficients on manifold}%
   {D affine connection coefficients, manifold}%
\Index
   {$D$\hyph vector space}%
   {D vector space}%
\Index
   {\dcr vector space}%
   {dcr vector space}%
\Index
   {determinant of matrix}%
   {determinant}%
\Index
   {deviation of trajectories}%
   {deviation of trajectories}%
\Index
   {diagonal in bundle}%
   {diagonal in bundle}%
\Index
   {diagram of correspondences}%
   {diagram of correspondences}%
\Index
   {diagram of representations}%
   {diagram of representations}%
\Index
   {dimension of \rcd vector space}%
   {dimension of vector space}%
\Index
   {direct product of bundles}%
   {Cartesian product of bundles}%
\Index
   {direct product of $D$\Hyph vector spaces}%
   {direct product of D vector spaces}%
\Index
   {direct product of division rings}%
   {direct product of division rings}%
\Index
   {direct product of \rcd vector spaces}%
   {direct product, rcd vector space}%
\Index
   {direct product of representations of fibered group}%
   {direct product of representations of fibered group}%
\Index
   {direct product of representations of group}%
   {direct product of representations of group}%
\Index
   {direct product of total spaces}%
   {Cartesian product of total spaces}%
\Index
   {direct product of \Ts representations of group}%
   {direct product of Tstar representations of group}%
\Index
   {direct sum of representations}%
   {direct sum of representations}%
\Index
   {distributive law for $A\star$\Hyph module}%
   {distributive law, Astar module over algebra}%
\Index
   {distributive law for $A\star$\Hyph vector space}%
   {distributive law, Astar vector space}%
\Index
   {distributive law for $D$\Hyph module}%
   {distributive law, D module}%
\Index
   {distributive law for $\mathcal D\star$\Hyph vector fields}%
   {distributive law, Dstar vector fields}%
\Index
   {distributive law for $D\star$\Hyph vector space}%
   {distributive law, Dstar vector space}%
\Index
   {distributive law for $\star A$\Hyph module}%
   {distributive law, starA module over algebra}%
\Index
   {double determinant}%
   {double determinant}%
\Index
   {\Drc basis for vector  bundle}%
   {drc basis, vector bundle}%
\Index
   {\drc basis of \subs rows vector space}%
   {drc basis, r rows vector space}%
\Index
   {\Drc linear map of vector bundles}%
   {drc linear map of vector bundles}%
\Index
   {\Drc linearly dependent vector fields}%
   {linearly dependent vector fields}%
\Index
   {\Drc linearly independent vector fields}%
   {linearly independent vector fields}%
\Index
   {\drc representation of group}%
   {drc linear representation of group}%
\Index
   {\drc vector}%
   {drc vector}%
\Index
   {\drc vector function}%
   {drc vector function}%
\Index
   {\drc vector space}%
   {drc vector space}%
\Index
   {$D\star$\Hyph affine space}%
   {Dstar affine space}%
\Index
   {$D\star$\Hyph antilinear homomorphism}%
   {Dstar antilinear homomorphism}%
\Index
   {$D\star$\Hyph antilinear mapping of ring with conjugation}%
   {antilinear mapping of ring with conjugation}%
\Index
   {\Ds component of coordinates of vector $\Vector r$}%
   {Dstar component of coordinates of vector, D vector space}%
\Index
   {$D\star$\Hyph linear homomorphism}%
   {Dstar linear homomorphism}%
\Index
   {$\mathcal D\star$\Hyph vector bundle}%
   {Dstar vector bundle}%
\Index
   {$\mathcal D\star$\Hyph vector field}%
   {Dstar vector field}%
\Index
   {$D\star$\hyph  vector space}%
   {Dstar vector space}%
\Index
   {$\mathcal D\star$\hyph linear composition of vector fields}%
   {linear composition of vector fields}%
\Index
   {$D\star$\Hyph module}%
   {Dstar-module}%
\Index
   {$\mathcal D\star$\hyph product of vector field over scalar}%
   {Dstar product of vector field over scalar, vector space}%
\Index
   {$D\star$\hyph product of vector over scalar}%
   {Dstar product of vector over scalar, vector space}%
\Index
   {dual space of \rcd vector space}%
   {dual space of rcd vector space}%
\Index
   {duality principle for biring}%
   {duality principle for biring}%
\Index
   {duality principle for biring of matrices}%
   {duality principle for biring of matrices}%
\Index
   {\rcd basis for \sups rows vector space}%
   {rcd basis, c rows vector space}%
\Index
   {\rcd linear span in vector space}%
   {linear span, vector space}%
\SetIndexSpace%
\Index
   {effective representation of division ring}%
   {effective representation of division ring}%
\Index
   {effective representation of fibered $\Omega$\Hyph algebra}%
   {effective representation of fibered Omega-algebra}%
\Index
   {effective representation of group}%
   {effective representation of group}%
\Index
   {effective representation of $\Omega$\Hyph algebra $A$}%
   {effective representation of algebra}%
\Index
   {effective representation of ring}%
   {effective representation of ring}%
\Index
   {effective \Ts representation of fibered division ring}%
   {effective representation of fibered division ring}%
\Index
   {effective \Ts representation of fibered group}%
   {effective representation of fibered group}%
\Index
   {effective \Ts representation of group}%
   {effective Tstar representation of group}%
\Index
   {endomorphism of representation of $\Omega$\Hyph algebra}%
   {endomorphism of representation}%
\Index
   {endomorphism of representation regular on generating set $X$}%
   {endomorphism of representation, regular on set}%
\Index
   {endomorphism of representation singular on generating set $X$}%
   {endomorphism of representation, singular on set}%
\Index
   {endomorphism of tower of representations}%
   {endomorphism of tower of representations}%
\Index
   {endomorphism of tower of representations regular on tuple of generating sets}%
   {endomorphism of representation, regular on tuple}%
\Index
   {endomorphism of tower of representations singular on tuple of generating sets}%
   {endomorphism of representation, singular on tuple}%
\Index
   {enhanced Lie group}%
   {enhanced Lie group}%
\Index
   {equivalence generated by representation $f$}%
   {equivalence of representation}%
\Index
   {essential parameters in a set of functions}%
   {essential parameters}%
\Index
   {Euclidean metric on division ring}%
   {Euclidean metric on division ring}%
\Index
   {Euclidean scalar product in $D$\Hyph vector space}%
   {Euclidean scalar product, vector space}%
\Index
   {Euclidean scalar product on division ring}%
   {Euclidean scalar product on division ring}%
\Index
   {extended matrix of \drc linear equations}%
   {extended matrix, system of drc linear equations}%
\Index
   {extended matrix of \rcd linear equations}%
   {extended matrix, system of rcd linear equations}%
\Index
   {extension of correspondence}%
   {extension of correspondence}%
\Index
   {extreme line}%
   {extreme line}%
\SetIndexSpace%
\Index
   {fibered coordinate \Drc isomorphism}%
   {fibered coordinate drc isomorphism}%
\Index
   {fibered correspondence from $\Bundle A$ to $\Bundle B$}%
   {fibered correspondence from A to B}%
\Index
   {fibered correspondence in $\Bundle{A}$}%
   {fibered correspondence in A}%
\Index
   {fibered correspondence of homomorphism}%
   {fibered correspondence of homomorphism}%
\Index
   {fibered equivalence}%
   {fibered equivalence}%
\Index
   {fibered group}%
   {fibered group}%
\Index
   {fibered identification morphism}%
   {fibered identification morphism}%
\Index
   {fibered little group}%
   {fibered little group}%
\Index
   {fibered morphism from bundle $\Bundle A$ into $\Bundle B$}%
   {fibered morphism from A into B}%
\Index
   {fibered natural morphism}%
   {fibered natural morphism}%
\Index
   {fibered $\Omega$\Hyph algebra}%
   {fibered Omega-algebra}%
\Index
   {fibered $\Omega$\Hyph subalgebra}%
   {fibered Omega-subalgebra}%
\Index
   {fibered ordering}%
   {fibered ordering}%
\Index
   {fibered preordering}%
   {fibered preordering}%
\Index
   {fibered ring}%
   {fibered ring}%
\Index
   {fibered stability group}%
   {fibered stability group}%
\Index
   {fibered subset}%
   {fibered subset}%
\Index
   {field-strength tensor}%
   {field-strength tensor}%
\Index
   {filter $\mathfrak{F}$ converges to $A$}%
   {filter converges}%
\Index
   {Finsler metric}%
   {Finsler metric}%
\Index
   {Finsler space}%
   {Finsler space}%
\Index
   {Finsler structure}%
   {Finsler structure}%
\Index
   {first Newton law}%
   {First Newton law}%
\Index
   {free $A$\Hyph module}%
   {free A module}%
\Index
   {free algebra over ring}%
   {free algebra over ring}%
\Index
   {free module over ring}%
   {free module over ring}%
\Index
   {free \Ts representation of fibered group}%
   {free representation of fibered group}%
\Index
   {free \Ts representation of group}%
   {free representation of group}%
\Index
   {Frenet transport}%
   {Frenet transport}%
\Index
   {function continuous with respect to set of arguments}%
   {function continuous with respect to set of arguments}%
\Index
   {function homogeneous of degree $k$}%
   {function homogeneous}%
\Index
   {function of $\gi n$ $D$\Hyph valued variables}%
   {function of n D valued variables}%
\Index
   {function of algebra differentiable in the G\^ateaux sense}%
   {function differentiable in Gateaux sense, algebra}%
\Index
   {function of $D$\Hyph vector space $\Vector{V}$ to $D$|Hyph vector space $\Vector W$ differentiable in the G\^ateaux sense}%
   {function differentiable in Gateaux sense, D vector space}%
\Index
   {function of division ring differentiable in the G\^ateaux sense}%
   {function differentiable in Gateaux sense, division ring}%
\Index
   {function of division ring \Ds differentiable in the Fr\'echet sense}%
   {function Dstar differentiable in Frechet sense, division ring}%
\Index
   {fundamental sequence in normed algebra}%
   {fundamental sequence, normed algebra}%
\Index
   {fundamental sequence in normed ring}%
   {fundamental sequence, normed ring}%
\Index
   {fundamental sequence in valued division ring}%
   {fundamental sequence, valued division ring}%
\Index
   {fundamental sequence in valued ring}%
   {fundamental sequence, valued ring}%
\SetIndexSpace%
\Index
   {$G$\Hyph reference frame}%
   {G reference frame}%
\Index
   {$G$\Hyph basis of vector space}%
   {G-basis}%
\Index
   {$G$\Hyph coordinates of basis}%
   {G-coordinates}%
\Index
   {$G$\Hyph space}%
   {GSpace}%
\Index
   {the G\^ateaux \crd derivative of map $\Vector f$ of $D$\hyph vector space $\Vector V$ to $D$\hyph vector space $\Vector W$}%
   {Gateaux crd derivative of map, D vector space}%
\Index
   {the G\^ateaux derivative of map $f$}%
   {Gateaux derivative of map, division ring}%
\Index
   {the G\^ateaux derivative of map $\Vector f$ of normed $D$\Hyph vector space $\Vector{V}$ to normed $D$\Hyph vector space $\Vector{W}$}%
   {Gateaux derivative of map, D vector space}%
\Index
   {the G\^ateaux derivative of mapping of algebra}%
   {Gateaux derivative of map, algebra}%
\Index
   {the G\^ateaux derivative of order $n$ of map $\Vector f$}%
   {Gateaux derivative of Order n, D vector space}%
\Index
   {the G\^ateaux derivative of order $n$ of map $f$ of division ring}%
   {Gateaux derivative of Order n, division ring}%
\Index
   {the G\^ateaux derivative of order $n$ of map $f$ of algebra}%
   {Gateaux derivative of Order n, algebra}%
\Index
   {the G\^ateaux derivative of second order of mapping of algebra}%
   {Gateaux derivative of Second Order, algebra}%
\Index
   {the G\^ateaux derivative of second order of map of division ring}%
   {Gateaux derivative of Second Order, division ring}%
\Index
   {the G\^ateaux differential of map $\Vector f$ of normed $D$\Hyph vector space $\Vector{V}$ to normed $D$\Hyph vector space $\Vector{W}$}%
   {Gateaux differential of map, D vector space}%
\Index
   {the G\^ateaux differential of map $f$}%
   {Gateaux differential of map, division ring}%
\Index
   {the G\^ateaux differential of mapping $f$ of algebra}%
   {Gateaux differential of map, algebra}%
\Index
   {the G\^ateaux differential of second order of map of division ring}%
   {Gateaux differential of Second Order, division ring}%
\Index
   {the G\^ateaux differential of second order of mapping $\Vector f$}%
   {Gateaux differential of Second Order, D vector space}%
\Index
   {the G\^ateaux differential of second order of mapping $f$ of algebra}%
   {Gateaux differential of Second Order, algebra}%
\Index
   {the G\^ateaux \drc derivative of map $\Vector f$ of $D$\Hyph vector space $\Vector V$ to $D$\Hyph vector space $\Vector W$}%
   {Gateaux drc derivative of map, D vector space}%
\Index
   {the G\^ateaux \Ds derivative of map $f$ of division ring $D$}%
   {Gateaux Dstar derivative of map, division ring}%
\Index
   {the G\^ateaux mixed partial derivative of map $f^j$ with respect to variables $v^i$, $v^j$}%
   {Gateaux partial derivative of Second Order, D vector space}%
\Index
   {the G\^ateaux partial derivative of map $f^j$ with respect to variable $v^i$}%
   {Gateaux partial derivative, D vector space}%
\Index
   {the G\^ateaux partial \drc derivative of map $f^b$ with respect to variable $x^a$}%
   {Gateaux partial drc derivative of map with respect to variable, D vector space}%
\Index
   {the G\^ateaux partial \drc derivative of map $f^b$ with respect to variable $x^a$}%
   {Gateaux partial crd derivative of map with respect to variable, D vector space}%
\Index
   {the G\^ateaux \sD derivative of map $f$ of division ring $D$}%
   {Gateaux starD derivative of map, division ring}%
\Index
   {generating set of representation}%
   {generating set of representation}%
\Index
   {generating set of subrepresentation}%
   {generating set of subrepresentation}%
\Index
   {generator of linear mapping}%
   {generator of linear map, division ring}%
\Index
   {geometric object defined in $\Omega_2$\Hyph algebra $M$}%
   {geometric object, representation g}%
\Index
   {geometric object defined in \rcd vector space}%
   {geometric object, rcd vector space}%
\Index
   {geometric object defined in tuple of $\VX\Omega$\Hyph algebras $\VX A$}%
   {geometric object, tower of representations g}%
\Index
   {geometric object in coordinate representation}%
   {geometric object, coordinate vector space}%
\Index
   {geometric object in coordinate representation defined in $\Omega_2$\Hyph algebra $M$}%
   {geometric object, coordinate representation g}%
\Index
   {geometric object in coordinate representation defined in \rcd vector space}%
   {geometric object, coordinate rcd vector space}%
\Index
   {geometric object in coordinate representation defined in tuple of $\VX\Omega$\Hyph algebras $\VX A$}%
   {geometric object, coordinate tower of representations g}%
\Index
   {geometric object in vector space}%
   {geometric object, vector space}%
\Index
   {geometric object of type $H$}%
   {geometric object of type H, representation g}%
\Index
   {geometric object of type $A$ in vector space}%
   {geometric object of type A, vector space}%
\Index
   {group algebra}%
   {group algebra}%
\Index
   {group of automorphisms of representation}%
   {group of automorphisms of representation}%
\SetIndexSpace%
\Index
   {Hadamard inverse of matrix}%
   {Hadamard inverse of matrix}%
\Index
   {Hamel basis}%
   {Hamel basis}%
\Index
   {hermitian conjugated vector}%
   {hermitian conjugated vector}%
\Index
   {hermitian conjugation in division ring}%
   {hermitian conjugation, division ring}%
\Index
   {hermitian matrix}%
   {hermitian matrix}%
\Index
   {hermitian metric on division ring}%
   {hermitian metric on division ring}%
\Index
   {hermitian scalar product in $D$\Hyph vector space}%
   {hermitian scalar product, vector space}%
\Index
   {hermitian scalar product on division ring}%
   {hermitian scalar product on division ring}%
\Index
   {holonomic coordinates of connection}%
   {holonomic coordinates of connection}%
\Index
   {holonomic coordinates of vector}%
   {vector holonomic coordinates}%
\Index
   {homogeneous bundle of fibered group}%
   {homogeneous bundle of fibered group}%
\Index
   {homogeneous linear geometric object}%
   {homogeneous linear geometric object}%
\Index
   {homogeneous linear representation of Lie group}%
   {homogeneous Linear Representation of Lie Group}%
\Index
   {homogeneous map of degree $k$ over field $F$}%
   {homogeneous map of degree over field, D vector space}%
\Index
   {homogeneous space of group}%
   {homogeneous space of group}%
\Index
   {homomorphism of fibered groups}%
   {homomorphism of fibered groups}%
\Index
   {homomorphism of fibered universal algebras}%
   {homomorphism of fibered universal algebras}%
\SetIndexSpace%
\Index
   {infinitesimal generator}%
   {infinitesimal generator}%
\Index
   {infinitesimal generators of group Lie}%
   {infinitesimal generators of group Lie}%
\Index
   {invariance principle in \drc vector space}%
   {invariance principle}%
\Index
   {invariance principle in representation of universal algebra}%
   {invariance principle, representation g}%
\Index
   {invariance principle in tower of representations of universal algebras}%
   {invariance principle, tower of representations g}%
\Index
   {invariance principle in vector space}%
   {invariance principle, vector space}%
\Index
   {inverse fibered correspondence}%
   {inverse fibered correspondence}%
\Index
   {inverse reduced fibered correspondence}%
   {inverse reduced fibered correspondence}%
\Index
   {involution in quaternion algebra}%
   {involution, quaternion algebra}%
\Index
   {isomorphism of fibered $\Omega$\Hyph algebras}%
   {isomorphism of fibered Omega-algebras}%
\Index
   {isomorphism of repesentations of $\Omega$\Hyph algebra}%
   {isomorphism of repesentations of Omega algebra}%
\Index
   {isotropic vector}%
   {isotropic vector}%
\SetIndexSpace%
\Index
   {$(^j_i)$\hyph $\RCcirc$\Hyph quasideterminant}%
   {j i RCcirc-quasideterminant}%
\Index
   {Jacobian complete system of differential equations}%
   {Jacobian complete system of differential equations}%
\Index
   {Jacobian complete system of \drv differential equations}%
   {Jacobian complete system of drc differential equations}%
\Index
   {$(ji)$\hyph quasideterminant}%
   {j i quasideterminant}%
\Index
   {the Jacobi\Hyph G\^ateaux matrix of map of $D$\Hyph vector space}%
   {Jacobi Gateaux matrix of map, D vector space}%
\SetIndexSpace%
\Index
   {kernel of inefficiency of representation of fibered group}%
   {kernel of inefficiency of representation of fibered group}%
\Index
   {kernel of inefficiency of representation of group}%
   {kernel of inefficiency of representation of group}%
\Index
   {kernel of inefficiency of \Ts representation of group $G$}%
   {kernel of inefficiency of Tstar representation of group}%
\Index
   {kernel of linear mapping of $D$\Hyph vector space}%
   {kernel of linear map, D vector space}%
\Index
   {kernel of linear mapping of division ring}%
   {kernel of linear map, division ring}%
\Index
   {Killing equation}%
   {Killing equation}%
\Index
   {Killing equation of second type}%
   {Killing equation second type}%
\Index
   {Killing vector of second type}%
   {Killing vector second type}%
\Index
   {Kronecker symbol}%
   {Kronecker symbol}%
\SetIndexSpace%
\Index
   {left cofactor of entry of matrix}%
   {left cofactor, matrix}%
\Index
   {left defined Lie algebra of Lie group}%
   {left defined Lie algebra}%
\Index
   {left double cofactor of entry of matrix}%
   {left double cofactor}%
\Index
   {left invariant vector field}%
   {left invariant vector}%
\Index
   {left module over $D$\Hyph algebra $A$}%
   {left module over algebra}%
\Index
   {left module over a ring $D$}%
   {left module over ring}%
\Index
   {left shift of $R$\Hyph module}%
   {left shift of module}%
\Index
   {left shift on fibered group}%
   {Tstar shift, fibered group}%
\Index
   {left shift on group}%
   {left shift}%
\Index
   {left shift on group}%
   {left shift, group}%
\Index
   {left structural constant of Lie algebra}%
   {left structural constant of Lie algebra}%
\Index
   {left vector space}%
   {left vector space}%
\Index
   {left-ordered cycle notation of permutation}%
   {left-ordered cycle notation of permutation}%
\Index
   {left-side contravariant representation of group}%
   {left-side contravariant representation of group}%
\Index
   {left-side covariant representation of group}%
   {left-side covariant representation of group}%
\Index
   {left-side representation of fibered $\Omega$\Hyph algebra}%
   {left-side representation of fibered Omega-algebra}%
\Index
   {left-side representation of $\Omega_1$\Hyph algebra $A$ in $\Omega_2$\Hyph algebra $M$}%
   {left-side representation of algebra}%
\Index
   {left-side transformation}%
   {left-side transformation}%
\Index
   {left-side transformation on bundle}%
   {left-side transformation of bundle}%
\Index
   {Lie algebra of Lie group}%
   {algebra Lie group Lie}%
\Index
   {Lie derivative}%
   {Lie derivative}%
\Index
   {Lie derivative of connection}%
   {Lie derivative of connection}%
\Index
   {Lie derivative of metric}%
   {Lie derivative of metric}%
\Index
   {Lie group basic operators}%
   {Lie group basic operators}%
\Index
   {lift of correspondence}%
   {lift of correspondence}%
\Index
   {lift of mapping}%
   {lift of map}%
\Index
   {limit of correspondence with respect to the filter}%
   {limit of correspondence with respect to the filter}%
\Index
   {limit of filter}%
   {limit of filter}%
\Index
   {limit of sequence in normed algebra}%
   {limit of sequence, normed algebra}%
\Index
   {limit of sequence in normed ring}%
   {limit of sequence, normed ring}%
\Index
   {limit of sequence in ring}%
   {limit of sequence, valued ring}%
\Index
   {limit of sequence in valued division ring}%
   {limit of sequence, valued division ring}%
\Index
   {limit set of filter}%
   {limit set of filter}%
\Index
   {linear combination of vectors of $A$\Hyph module}%
   {linear combination in A module}%
\Index
   {linear dependent vectors of $A$\Hyph module}%
   {linear dependent vectors, module}%
\Index
   {linear geometric object}%
   {linear geometric object}%
\Index
   {linear independent vectors of $A$\Hyph module}%
   {linear independent vectors, module}%
\Index
   {linear map of $D$\Hyph vector spaces over field $F$}%
   {linear map over field, vector space}%
\Index
   {linear mapping of $R$\Hyph algebra $A_1$ into $R$\Hyph algebra $A_2$}%
   {linear mapping of R algebras}%
\Index
   {linear mapping of $R$\Hyph module $A_1$ into $R$\Hyph module $A_2$}%
   {linear mapping of R modules}%
\Index
   {linear mapping of division ring}%
   {linear mapping of division ring}%
\Index
   {linear mapping of division ring generated by mapping $G$}%
   {linear map generated by map, division ring}%
\Index
   {linear mapping of $R_1$\Hyph module $A_1$ into $R_2$\Hyph module $A_2$}%
   {linear mapping of R_1 module into R_2 module}%
\Index
   {linear representation of group}%
   {linear representation of group}%
\Index
   {linear representation of Lie group}%
   {Linear Representation of Lie Group}%
\Index
   {linear transformation of \rcd affine space}%
   {linear transformation, rcd affine space}%
\Index
   {little group}%
   {little group}%
\Index
   {local reference frame}%
   {local reference frame}%
\Index
   {locally compact at point $p$ space}%
   {locally compact at point space}%
\Index
   {locally compact space}%
   {locally compact space}%
\Index
   {Lorentz transformation}%
   {Lorentz transformation}%
\SetIndexSpace%
\Index
   {$m$\Hyph dimensional parallelepiped}%
   {m dimensional parallelepiped}%
\Index
   {$m$\Hyph vector}%
   {m-vector}%
\Index
   {manifold with $D$\Hyph affine connections}%
   {manifold with D- affine connections}%
\Index
   {map of type $G$ on manifold}%
   {map of type G on manifold}%
\Index
   {map polylinear over finite dimensional algebras}%
   {map polylinear over finite dimensional algebras}%
\Index
   {mapping of rings polylinear over commutative ring}%
   {map polylinear over commutative ring, ring}%
\Index
   {mapping space}%
   {mapping space}%
\Index
   {matrix of \Acr linear mapping of module}%
   {matrix of Acr linear map}%
\Index
   {matrix of antilinear homomorphism}%
   {matrix of antilinear homomorphism}%
\Index
   {matrix of $A\star$\Hyph linear map}%
   {matrix of Astar linear map}%
\Index
   {matrix of bilinear function}%
   {matrix of bilinear function}%
\Index
   {matrix of endomorphisms of $\Omega$\Hyph algebra}%
   {matrix of endomorphisms of Omega algebra}%
\Index
   {matrix of fibered \Drc linear map}%
   {matrix of fibered drc linear map}%
\Index
   {matrix of linear homomorphism}%
   {matrix of linear homomorphism}%
\Index
   {matrix of linear mappings}%
   {matrix of linear mappings}%
\Index
   {matrix of mappings}%
   {matrix of mappings}%
\Index
   {matrix of quadratic map}%
   {matrix of quadratic map, division ring}%
\Index
   {matrix of \rcd linear mapping}%
   {matrix of rcd linear map}%
\Index
   {metric tensor in Minkowski space}%
   {metric tensor, Minkowski space}%
\Index
   {metric-affine manifold}%
   {metric-affine manifold}%
\Index
   {Minkowski space}%
   {Minkowski space, Finsler}%
\Index
   {module over ring}%
   {module over ring}%
\Index
   {morphism from tower of \Ts representations into tower of \Ts representations}%
   {morphism from tower of representations into tower of representations}%
\Index
   {morphism of fibered \Ts representations from $\Bundle F$ into $\Bundle G$}%
   {morphism of fibered representations from f into g}%
\Index
   {morphism of representation $f$}%
   {morphism of representation f}%
\Index
   {morphism of representations from $f$ into $g$}%
   {morphism of representations from f into g}%
\Index
   {morphism of representations of $\Omega_1$\Hyph algebra in $\Omega_2$\Hyph algebra}%
   {morphism of representations of Omega1 algebra in Omega2 algebra}%
\Index
   {morphism of \Ts representations of fibered $\Omega$\Hyph algebra}%
   {morphism of representations of fibered Omega algebra}%
\Index
   {motion of Minkowski space}%
   {motion, Minkowski space}%
\Index
   {movement on basis manifold}%
   {movement transformation}%
\SetIndexSpace%
\Index
   {$n$\Hyph algebra over the ring}%
   {n algebra over ring}%
\Index
   {$n$\Hyph ary fibered relation}%
   {fibered relation}%
\Index
   {nonmetricity}%
   {nonmetricity}%
\Index
   {nonsingular bilinear function}%
   {nonsingular bilinear function}%
\Index
   {nonsingular system of \rcd linear equations}%
   {nonsingular system of linear equations}%
\Index
   {nonsingular tensor}%
   {nonsingular tensor, algebra}%
\Index
   {nonsingular transformation}%
   {nonsingular transformation}%
\Index
   {norm in quaternion algebra}%
   {norm, quaternion algebra}%
\Index
   {norm of map $\Vector A$ of normed $D$\hyph vector space}%
   {norm of map, D vector space}%
\Index
   {norm of mapping into $D$\Hyph algebra}%
   {norm of map, algebra}%
\Index
   {norm of mapping of division ring}%
   {norm of map, division ring}%
\Index
   {norm on $D$\Hyph algebra}%
   {norm on D algebra}%
\Index
   {norm on $D$\Hyph vector space}%
   {norm on D vector space}%
\Index
   {norm on ring}%
   {norm on ring}%
\Index
   {normed $D$\Hyph algebra}%
   {normed D algebra}%
\Index
   {normed $D$\Hyph vector space}%
   {normed D vector space}%
\Index
   {normed ring}%
   {normed ring}%
\Index
   {not complete group}%
   {not complete group}%
\Index
   {not complete $\Omega$\Hyph algebra}%
   {not complete Omega algebra}%
\Index
   {nucleus of $R$\Hyph algebra $A$}%
   {nucleus of algebra}%
\SetIndexSpace%
\Index
   {octonion algebra}%
   {octonion algebra}%
\Index
   {$\Omega$\Hyph group}%
   {Omega group}%
\Index
   {$\Omega$\Hyph linear mapping}%
   {Omega linear map}%
\Index
   {$\Omega_2$\Hyph word of element of representation relative to generating set}%
   {word of element relative to generating set, representation}%
\Index
   {operation on bundle}%
   {operation on bundle}%
\Index
   {opposite algebra to algebra $P$}%
   {opposite algebra}%
\Index
   {opposite fibered preordering}%
   {opposite fibered preordering}%
\Index
   {orbit of linear mapping}%
   {orbit of linear mapping}%
\Index
   {orbit of representation of fibered group}%
   {orbit of representation of fibered group}%
\Index
   {orbit of representation of group}%
   {orbit of representation of group}%
\Index
   {orbit of \Ts representation of group}%
   {orbit of Tstar  representation of group}%
\Index
   {origin of coordinate system of affine space}%
   {origin of coordinate system of affine space}%
\Index
   {origin of coordinate system of $\star D$\Hyph affine space}%
   {origin of coordinate system of starD affine space}%
\Index
   {orthogonal basis in Minkowski space}%
   {orthogonal basis, Minkowski space}%
\Index
   {orthogonality in Minkowski space}%
   {Minkowski orthogonality}%
\Index
   {orthonormal basis in Minkowski space}%
   {orthonormal basis, Minkowski space}%
\Index
   {orthonornal basis}%
   {Orthonornal Basis}%
\Index
   {orthonornal basis}%
   {Orthonornal Basis, division ring}%
\SetIndexSpace%
\Index
   {passive representation of group $G(f)$ in basis manifold of representation}%
   {passive representation in basis manifold}%
\Index
   {parallel shift of \rcd affine space}%
   {parallel shift, rcd affine space}%
\Index
   {parallelogram}%
   {parallelogram}%
\Index
   {partial linear mapping of variable $v^i$}%
   {partial linear map of variable}%
\Index
   {passive representation}%
   {passive representation}%
\Index
   {passive representation of group $G(\Vector f)$ in basis manifold of tower of representations}%
   {passive representation in basis manifold, tower of representations}%
\Index
   {passive \sT representation}%
   {passive starT representation}%
\Index
   {passive transformation of the basis manifold of representation}%
   {passive transformation of basis, representation}%
\Index
   {passive transformation of the basis manifold of tower of representations}%
   {passive transformation of basis, tower of representations}%
\Index
   {passive transformation on basis manifold}%
   {passive transformation}%
\Index
   {passive transformation on the set of \rcd bases}%
   {passive transformation, vector space}%
\Index
   {permutability property of trace}%
   {permutability property of trace}%
\Index
   {pfaffian derivative}%
   {pfaffian derivative}%
\Index
   {polylinear mapping of $(n)$\hyph $D$\hyph vector spaces}%
   {polylinear map of D vector spaces}%
\Index
   {polylinear mapping of algebras}%
   {polylinear map of algebras}%
\Index
   {polylinear mapping of modules}%
   {polylinear map of modules}%
\Index
   {polylinear skew symmetric map}%
   {polylinear map skew symmetric, division ring}%
\Index
   {polylinear symmetric map}%
   {polylinear map symmetric, division ring}%
\Index
   {polymorphism of representations}%
   {polymorphism of representations}%
\Index
   {polyvector}%
   {polyvector}%
\Index
   {potential energy}%
   {potential energy}%
\Index
   {product of geometric object and constant}%
   {product of geometric object and constant}%
\Index
   {product of geometric object and constant in vector space}%
   {product of geometric object and constant, vector space}%
\Index
   {product of groups}%
   {product of groups}%
\Index
   {product of morphisms of representations of universal algebra}%
   {product of morphisms of representations of universal algebra}%
\Index
   {product of morphisms of tower of representations}%
   {product of morphisms of tower of representations}%
\Index
   {product of morphisms of \Ts representations of fibered $\Omega$\Hyph algebra}%
   {product of morphisms of representations of fibered Omega algebra}%
\Index
   {product of objects in category}%
   {product of objects in category}%
\Index
   {projection of bundle $\Bundle E$ along fiber $E$}%
   {projection of bundle along fiber}%
\Index
   {pseudo\Hyph Euclidean metric on division ring}%
   {pseudo-Euclidean metric on division ring}%
\Index
   {pseudo\Hyph Euclidean scalar product in $D$\Hyph vector space}%
   {pseudo-Euclidean scalar product, vector space}%
\Index
   {pseudo-Euclidean scalar product on division ring}%
   {pseudo-Euclidean scalar product on division ring}%
\SetIndexSpace%
\Index
   {quadratic form in division ring}%
   {quadratic form, division ring}%
\Index
   {quadratic map of division ring}%
   {Quadratic Map of Division Ring}%
\Index
   {quasi affine transformation on basis manifold}%
   {quasi affine transformation}%
\Index
   {quasi affine transformation on basis manifold}%
   {quasi affine drc transformation}%
\Index
   {quasi movement on basis manifold}%
   {quasi movement, division ring}%
\Index
   {quasi movement on basis manifold}%
   {quasi movement}%
\Index
   {quasiclosed ring of mappings}%
   {quasiclosed ring of mappings}%
\Index
   {quasideterminant}%
   {quasideterminant definition}%
\Index
   {quasimotion of Minkowski space}%
   {Quasimotion, Minkowski space}%
\Index
   {quaternion algebra}%
   {quaternion algebra}%
\Index
   {quaternion algebra $E$ over the field $F$}%
   {quaternion algebra over the field}%
\Index
   {quotient bundle}%
   {quotient bundle}%
\SetIndexSpace%
\Index
   {\drc basis dual to \rcd basis of vector space}%
   {basis dual to basis, rcd vector space}%
\Index
   {$(^j_i)$\hyph \RC quasideterminant}%
   {j i RC-quasideterminant}%
\Index
   {\sups row of matrix}%
   {r row}%
\Index
   {$R$\Hyph module}%
   {R- module}%
\Index
   {$r$\hyph row of matrix}%
   {r-row}%
\Index
   {rank of Hermitian matrix by principal minors}%
   {rank of Hermitian matrix by principal minors}%
\Index
   {rank of quadratic map of division ring}%
   {rank of quadratic map, division ring}%
\Index
   {\RC inverse element of biring}%
   {rc-inverse element}%
\Index
   {\RC major minor}%
   {RC-major minor}%
\Index
   {\RC matrix group}%
   {rc-matrix group}%
\Index
   {\RC nonsingular matrix}%
   {RC nonsingular matrix}%
\Index
   {\RC power}%
   {rc power}%
\Index
   {\RC product of matrices}%
   {rc-product of matrices}%
\Index
   {$\RCcirc$\Hyph product of matrices of mappings}%
   {rc product of matrices of mappings}%
\Index
   {\RC quasideterminant}%
   {RC-quasideterminant}%
\Index
   {\RC rank of matrix}%
   {rc-rank of matrix}%
\Index
   {\RC singular matrix}%
   {RC singular matrix}%
\Index
   {$\RCcirc A$\Hyph basis for module}%
   {RCcirc A basis, module over algebra}%
\Index
   {$\RCcirc A$\Hyph linearly dependent set of vectors}%
   {RCcirc linearly dependent, starA module over algebra}%
\Index
   {$\RCcirc A$\Hyph linearly independent set of vectors}%
   {RCcirc linearly independent, starA module over algebra}%
\Index
   {$\RCcirc$\Hyph nonsingular matrix of $\mathcal A(A)$\Hyph mappings}%
   {RCcirc nonsingular matrix of A(A) mappings}%
\Index
   {$\RCcirc$\Hyph nonsingular matrix of endomorphisms}%
   {RCcirc nonsingular matrix of endomorphisms}%
\Index
   {$\RCcirc$\Hyph nonsingular system of additive equations}%
   {RCcirc nonsingular system of additive equations}%
\Index
   {$\RCcirc$\Hyph quasideterminant}%
   {RCcirc-quasideterminant definition}%
\Index
   {$\RCcirc$\Hyph singular matrix of $\mathcal A(A)$\Hyph mappings}%
   {RCcirc singular matrix of A(A) mappings}%
\Index
   {$\RCcirc$\Hyph singular matrix of endomorphisms}%
   {RCcirc singular matrix of endomorphisms}%
\Index
   {$\RCcirc A$\Hyph linear combination}%
   {RCcircA linear combination}%
\Index
   {\rcd affine basis}%
   {rcd affine basis, division ring}%
\Index
   {\rcd affine plane}%
   {rcd affine plane}%
\Index
   {\rcd affine space}%
   {rcd affine space}%
\Index
   {\rcd affine transformation}%
   {rcd affine transformation}%
\Index
   {\rcd automorphism of vector space}%
   {automorphism of vector space}%
\Index
   {\rcd basis for vector space}%
   {rcd basis, vector space}%
\Index
   {\rcd isomorphism of vector spaces}%
   {isomorphism of vector spaces}%
\Index
   {\rcd linear map of vector spaces}%
   {rcd linear map of vector spaces}%
\Index
   {\rcd linear \Ts representation of group}%
   {rcd linear Tstar representation of group}%
\Index
   {\rcd linearly dependent vectors}%
   {linearly dependent, vector space}%
\Index
   {\rcd linearly independent vectors}%
   {linearly independent, vector space}%
\Index
   {\rcd vector}%
   {rcd vector}%
\Index
   {\rcd vector space}%
   {rcd vector space}%
\Index
   {reduced Cartesian product of bundles}%
   {reduced Cartesian product of bundles}%
\Index
   {reduced Cartesian product of total spaces}%
   {reduced Cartesian product of total spaces}%
\Index
   {reduced fibered correspondence from $\Bundle{A}$ to $\Bundle B$}%
   {reduced fibered correspondence from A to B}%
\Index
   {reduced fibered correspondence in $\Bundle{A}$}%
   {reduced fibered correspondence in A}%
\Index
   {reduced polymorphism of representations}%
   {reduced polymorphism of representations}%
\Index
   {reducible biring}%
   {reducible biring}%
\Index
   {reference frame in event space}%
   {reference frame in event space}%
\Index
   {reference frame manifold}%
   {reference frame manifold}%
\Index
   {reflexive $2$\Hyph ary fibered relation}%
   {reflexive 2 ary fibered relation}%
\Index
   {regular endomorphism of representation}%
   {regular endomorphism of representation}%
\Index
   {regular endomorphism of tower of representations}%
   {regular endomorphism of tower of representations}%
\Index
   {regular quadratic map in division ring}%
   {regular quadratic map, division ring}%
\Index
   {representation of group}%
   {representation of group}%
\Index
   {representation of $\Omega$\Hyph algebra in representation}%
   {representation of Omega algebra in representation}%
\Index
   {representation of $\Omega$\Hyph algebra in tower of representations}%
   {representation of Omega algebra in tower of representations}%
\Index
   {representation of $\Omega$\Hyph algebra $A$ in category $\mathcal B$}%
   {representation of Omega algebra in category}%
\Index
   {representation of $\Omega_1$\Hyph algebra $A$ in $\Omega_2$\Hyph algebra $M$}%
   {representation of algebra}%
\Index
   {representative of geometric object in \drc vector space}%
   {representative of geometric object, drc vector space}%
\Index
   {representative of geometric object in $\Omega_2$\Hyph algebra}%
   {representative of geometric object, representation g}%
\Index
   {representative of geometric object in \rcd vector space}%
   {representative of geometric object, rcd vector space}%
\Index
   {representative of geometric object in tuple of $\VX\Omega$\Hyph algebras}%
   {representative of geometric object, tower of representations g}%
\Index
   {representative of geometric object in vector space}%
   {representative of geometric object, vector space}%
\Index
   {restriction of correspondence $\Phi$ to set $C$}%
   {restriction of correspondence}%
\Index
   {right cofactor of entry of matrix}%
   {right cofactor, matrix}%
\Index
   {right defined Lie algebra of Lie group}%
   {right defined Lie algebra}%
\Index
   {right double cofactor of entry of matrix}%
   {right double cofactor}%
\Index
   {right invariant vector field}%
   {right invariant vector}%
\Index
   {right module over $D$\Hyph algebra $A$}%
   {right module over algebra}%
\Index
   {right module over a ring $D$}%
   {right module over ring}%
\Index
   {right shift on group}%
   {right shift}%
\Index
   {right shift on group}%
   {right shift, group}%
\Index
   {right structural constant of Lie algebra}%
   {right structural constant of Lie algebra}%
\Index
   {right vector space}%
   {right vector space}%
\Index
   {right-ordered cycle notation of permutation}%
   {right-ordered cycle notation of permutation}%
\Index
   {right-side contravariant representation of group}%
   {right-side contravariant representation of group}%
\Index
   {right-side covariant representation of group}%
   {right-side covariant representation of group}%
\Index
   {right-side representation of fibered $\Omega$\Hyph algebra}%
   {right-side representation of fibered Omega-algebra}%
\Index
   {right-side representation of $\Omega_1$\Hyph algebra $A$ in $\Omega_2$\Hyph algebra $M$}%
   {right-side representation of algebra}%
\Index
   {right-side transformation}%
   {right-side transformation}%
\Index
   {ring has characteristic $0$}%
   {ring has characteristic 0}%
\Index
   {ring has characteristic $p$}%
   {ring has characteristic p}%
\Index
   {ring with conjugation}%
   {ring with conjugation}%
\Index
   {row determinant}%
   {row determinant}%
\Index
   {row vector}%
   {row vector}%
\SetIndexSpace%
\Index
   {$\star A$\Hyph module}%
   {starA-module}%
\Index
   {$D\star$\Hyph product of \rcd linear map $A$ over scalar}%
   {Dstar product of rcd linear map over scalar}%
\Index
   {$(\RCstar S,\RCstar T)$\Hyph linear map of vector spaces}%
   {rcs rct linear map of vector spaces}%
\Index
   {scalar algebra of algebra}%
   {scalar algebra of algebra}%
\Index
   {scalar algebra of ring}%
   {scalar algebra of ring}%
\Index
   {scalar of element of algebra}%
   {scalar of algebra}%
\Index
   {scalar of element of ring}%
   {scalar of ring}%
\Index
   {scalar of mapping}%
   {scalar of mapping}%
\Index
   {scalar potential}%
   {scalar potential}%
\Index
   {Schauder basis}%
   {Schauder basis}%
\Index
   {second Newton law}%
   {Second Newton law}%
\Index
   {section of bundle}%
   {section of bundle}%
\Index
   {set of coordinates of representation}%
   {coordinate set of representation}%
\Index
   {set of $\Omega_2$\Hyph words of representation}%
   {word set of representation}%
\Index
   {set of tuples of coordinates of tower of representations}%
   {coordinate set of tower of representations}%
\Index
   {set of tuples of $\Vector\Omega$\Hyph words of tower of representations}%
   {word set of tower of representations}%
\Index
   {simple polyvector}%
   {simple polyvector}%
\Index
   {single transitive representation of fibered $\Omega$\Hyph algebra}%
   {single transitive representation of fibered Omega-algebra}%
\Index
   {single transitive representation of group}%
   {single transitive representation of group}%
\Index
   {single transitive representation of $\Omega$\Hyph algebra $A$}%
   {single transitive representation of algebra}%
\Index
   {singular linear mapping of $D$\Hyph vector space}%
   {singular linear map, D vector space}%
\Index
   {singular linear mapping of division ring}%
   {singular linear map, division ring}%
\Index
   {skew product of vectors}%
   {skew product of vectors}%
\Index
   {skew symmetric polylinear mapping into associative algebra}%
   {polylinear map skew symmetric, associative algebra}%
\Index
   {space of orbits of \Ts representation}%
   {space of orbits of Ts representation}%
\Index
   {spacelike vector}%
   {spacelike vector}%
\Index
   {speed of deviation}%
   {speed of deviation}%
\Index
   {$(\mathcal S\RCstar,\mathcal T\RCstar)$\Hyph linear map of vector bundles}%
   {src trc linear map of vector bundles}%
\Index
   {($S\star$, $\star T$)\hyph bimodule}%
   {(Sstar,starT)-bimodule}%
\Index
   {stability group}%
   {stability group}%
\Index
   {stable set of representation}%
   {stable set of representation}%
\Index
   {standard component of the G\^ateaux differential of map $f$}%
   {standard component of Gateaux differential, division ring}%
\Index
   {standard component of polylinear map $f$ of division ring}%
   {standard component of polylinear map, division ring}%
\Index
   {standard component of polylinear mapping into associative algebra}%
   {standard component of polylinear map, associative algebra}%
\Index
   {standard component of quadratic map $f$ over field $F$}%
   {standard component of quadratic map, division ring}%
\Index
   {standard component of tensor}%
   {standard component of tensor, division ring}%
\Index
   {standard component of tensor in tensor product of algebras}%
   {standard component of tensor, algebra}%
\Index
   {standard component of the G\^ateaux derivative of mapping $f$}%
   {standard component of Gateaux derivative, algebra}%
\Index
   {standard component over field $F$ of bilitnear map $f$}%
   {standard component of bilinear map, division ring}%
\Index
   {standard coordinates of basis}%
   {standard coordinates of basis}%
\Index
   {standard coordinates of \rcd basis}%
   {standard coordinates of rcd basis}%
\Index
   {standard $F$\Hyph component of linear mapping $f$}%
   {standard component of linear map, division ring}%
\Index
   {standard $F$\Hyph representation of linear mapping of division ring}%
   {linear map, standard representation, division ring}%
\Index
   {standard representation of the G\^ateaux differential of map of division ring over field $F$}%
   {Gateaux differential, standard representation, division ring}%
\Index
   {standard representation of matrix}%
   {Standard representation}%
\Index
   {standard representation of polylinear map of division ring}%
   {polylinear map, standard representation, division ring}%
\Index
   {standard representation of polylinear mapping into associative algebra}%
   {polylinear map, standard representation, associative algebra}%
\Index
   {standard representation of quadratic map of division ring over field $F$}%
   {quadratic map, standard representation, division ring}%
\Index
   {standard representation over field $F$ of bilinear map of division ring}%
   {bilinear map, standard representation, division ring}%
\Index
   {$\star A$\Hyph product of \Acr linear map $\Vector f$ over scalar}%
   {starA product of Acr linear map over scalar}%
\Index
   {$\star A$\Hyph product of $A\star$\Hyph linear mapping over scalar}%
   {starA product of Astar linear map over scalar}%
\Index
   {$\star D$\hyph vector space}%
   {starD-vector space}%
\Index
   {$\star R$\hyph module}%
   {starR-module}%
\Index
   {$\star A$\Hyph product of vector over scalar}%
   {starA product of vector over scalar, starA module}%
\Index
   {$\star D$\Hyph affine space}%
   {starD affine space}%
\Index
   {\sD component of coordinates of vector $\Vector r$}%
   {starD component of coordinates of vector, D vector space}%
\Index
   {\sT representation of fibered group}%
   {starT representation of fibered group}%
\Index
   {\sT representation of fibered $\Omega$\Hyph algebra}%
   {starT representation of fibered Omega-algebra}%
\Index
   {\sT representation of $\Omega_1$\Hyph algebra $A$ in $\Omega_2$\Hyph algebra $M$}%
   {starT representation of algebra}%
\Index
   {\sT shift}%
   {starT shift}%
\Index
   {\sT shift on fibered group}%
   {starT shift, fibered group}%
\Index
   {\sT transformation}%
   {starT transformation}%
\Index
   {\sT transformation on bundle}%
   {starT transformation of bundle}%
\Index
   {structural constants of algebra $P$ over ring $D$}%
   {structural constants of algebra}%
\Index
   {structural constants of division ring $D$ over field $F$}%
   {structural constants of division ring over field}%
\Index
   {subbundle}%
   {subbundle}%
\Index
   {subbundle of $\mathcal D\star$\hyph vector space}%
   {subbundle of Dstar vector bundle}%
\Index
   {subrepresentation generated by set $X$}%
   {subrepresentation generated by set}%
\Index
   {subrepresentation of representation}%
   {subrepresentation of representation}%
\Index
   {sum of \Acr linear mappings of module}%
   {sum of Acr linear maps, module}%
\Index
   {sum of geometric objects in vector space}%
   {sum of geometric objects, vector space}%
\Index
   {sum of geometric objects}%
   {sum of geometric objects}%
\Index
   {sum of \rcd linear maps}%
   {sum of rcd linear maps, rcd vector spaces}%
\Index
   {superposition of coordinates of the representation $f$ and the element $m$}%
   {superposition of coordinates, representation}%
\Index
   {superposition of coordinates of the tower of representations $\Vector f$ and the element $\VX a$}%
   {superposition of coordinates, tower of representations}%
\Index
   {symmetric $2$\Hyph ary fibered relation}%
   {symmetric 2 ary fibered relation}%
\Index
   {symmetric bilinear map of $D$\Hyph vector space to division ring}%
   {symmetric bilinear map, vector space to division ring}%
\Index
   {symmetric polylinear mapping into associative algebra}%
   {polylinear map symmetric, associative algebra}%
\Index
   {symmetry group}%
   {symmetry group}%
\Index
   {symmetry group}%
   {SymmetryGroup}%
\Index
   {synchronization of reference frame}%
   {synchronization of reference frame}%
\Index
   {system of additive equations}%
   {system of additive equations}%
\Index
   {system of \drc linear equations}%
   {system of drc linear equations}%
\Index
   {system of \rcd linear equations}%
   {system of rcd linear equations}%
\Index
   {tandard representation of the G\^ateaux derivative of mapping over algebra}%
   {Gateaux derivative, standard representation, algebra}%
\SetIndexSpace%
\Index
   {Taylor polynomial}%
   {Taylor polynomial, division ring}%
\Index
   {Taylor series}%
   {Taylor series, division ring}%
\Index
   {tensor power of algebra}%
   {tensor power of algebra}%
\Index
   {tensor power of representation}%
   {tensor power of representation}%
\Index
   {tensor product of algebras}%
   {tensor product of algebras}%
\Index
   {tensor product of $D$\Hyph vector spaces}%
   {tensor product of D vector spaces}%
\Index
   {tensor product of division rings}%
   {tensor product of division rings}%
\Index
   {tensor product of \Ds vector spaces}%
   {tensor product of Dstar vector spaces}%
\Index
   {tensor product of representations}%
   {tensor product of representations}%
\Index
   {tensor product of representations}%
   {tensor product of representations}%
\Index
   {tensor product of rings over commutative ring}%
   {tensor product of rings}%
\Index
   {the Fr\'echet \Ds derivative of map $f$ of division ring $D$ at point $x$}%
   {Frechet Dstar derivative of map, division ring}%
\Index
   {timelike vector}%
   {timelike vector}%
\Index
   {topological $D$\Hyph vector space}%
   {topological D vector space}%
\Index
   {topological $D$\Hyph algebra}%
   {topological D algebra}%
\Index
   {topological division ring}%
   {topological division ring}%
\Index
   {topological \drc vector space}%
   {topological drc vector space}%
\Index
   {topological ring}%
   {topological ring}%
\Index
   {torsion form}%
   {torsion form}%
\Index
   {torsion tensor}%
   {torsion tensor}%
\Index
   {tower of bundles}%
   {tower of bundles}%
\Index
   {tower of effective representations}%
   {tower of effective representations}%
\Index
   {tower of representations of $\Vector{\Omega}$\Hyph algebras}%
   {tower of representations of algebras}%
\Index
   {tower of subrepresentations}%
   {tower of subrepresentations}%
\Index
   {tower of subrepresentations of tower of representations $\Vector f$ generated by tuple of sets $\VX X$}%
   {subrepresentation generated by tuple of sets}%
\Index
   {trace of quaternion}%
   {trace, quaternion algebra}%
\Index
   {transformation coordinated with equivalence}%
   {transformation coordinated with equivalence}%
\Index
   {transformation of universal algebra}%
   {transformation of universal algebra}%
\Index
   {transformation on bundle}%
   {transformation of bundle}%
\Index
   {transitive $2$\Hyph ary fibered relation}%
   {transitive 2 ary fibered relation}%
\Index
   {transitive representation of fibered $\Omega$\Hyph algebra}%
   {transitive representation of fibered Omega-algebra}%
\Index
   {transitive representation of group}%
   {transitive representation of group}%
\Index
   {transitive representation of $\Omega$\Hyph algebra $A$}%
   {transitive representation of algebra}%
\Index
   {\Ts linear composition of  vectors}%
   {linear composition of  vectors}%
\Index
   {\Ts matrices vector space}%
   {matrices vector space}%
\Index
   {\Ts representation of fibered $\Omega$\Hyph algebra}%
   {Tstar representation of fibered Omega-algebra}%
\Index
   {\Ts representation of $\Omega_1$\Hyph algebra $A$ in $\Omega_2$\Hyph algebra $M$}%
   {Tstar representation of algebra}%
\Index
   {\Ts shift}%
   {Tstar shift}%
\Index
   {\Ts transformation}%
   {Tstar transformation}%
\Index
   {\Ts transformation on bundle}%
   {Tstar transformation of bundle}%
\Index
   {tuple of coordinates of element $\Vector a$ relative to tuple of sets $\VX X$}%
   {coordinates of element, tower of representations}%
\Index
   {tuple of equivalence generated by tower of representations $\Vector f$}%
   {tuple of equivalence of tower of representations}%
\Index
   {tuple of generating sets of tower of representations}%
   {tuple of generating sets of tower of representations}%
\Index
   {tuple of generating sets of tower subrepresentations}%
   {tuple of generating sets of tower subrepresentations}%
\Index
   {tuple of $\Vector{\Omega}$\Hyph words of element of tower of representations relative to tuple of generating sets}%
   {tuple of words relative to tuple of generating sets, tower of representations}%
\Index
   {tuple of stable sets of tower of representation}%
   {tuple of stable sets of tower of representations}%
\Index
   {twin representations of associative algebra}%
   {twin representations of associative algebra}%
\Index
   {twin representations of $D$\Hyph algebra}%
   {twin representations of D algebra}%
\Index
   {twin representations of division ring}%
   {twin representations of division ring}%
\Index
   {twin representations of fibered group}%
   {twin representations of fibered group}%
\Index
   {twin representations of group}%
   {twin representations of group}%
\SetIndexSpace%
\Index
   {unit sphere in $D$\Hyph algebra}%
   {unit sphere in algebra}%
\Index
   {unit sphere in division ring}%
   {unit sphere in division ring}%
\Index
   {unit vector}%
   {unit vector}%
\Index
   {unitarity law for  $\mathcal D\star$\Hyph vector fields}%
   {unitarity law, Dstar vector fields}%
\Index
   {unitarity law for $A\star$\Hyph module}%
   {unitarity law, Astar module over algebra}%
\Index
   {unitarity law for $A\star$\Hyph vector space}%
   {unitarity law, Astar vector space}%
\Index
   {unitarity law for $D$\Hyph module}%
   {unitarity law, D module}%
\Index
   {unitarity law for $D\star$\Hyph vector space}%
   {unitarity law, Dstar vector space}%
\Index
   {unitarity law for $\star A$\Hyph module}%
   {unitarity law, starA module over algebra}%
\SetIndexSpace%
\Index
   {valued division ring}%
   {valued division ring}%
\Index
   {valued ring}%
   {valued ring}%
\Index
   {vector bundle}%
   {vector bundle}%
\Index
   {vector module of algebra}%
   {vector module of algebra}%
\Index
   {vector module of ring}%
   {vector module of ring}%
\Index
   {vector of element of algebra}%
   {vector of algebra}%
\Index
   {vector of element of ring}%
   {vector of ring}%
\Index
   {vector of mapping}%
   {vector of mapping}%
\Index
   {vector potential}%
   {vector potential}%
\Index
   {vector space over field}%
   {vector space over field}%
\Index
   {vector space type}%
   {vector space type}%

\CloseIndex

%% file: Symbol.English.tex
\def\indexname{Special Symbols and Notations}
\OpenIndex

\SetIndexSpace
\Symb
   {$(^a_b)$\hyph\CR quasideterminant}%
   {a b CR quasideterminant definition}%
\Symb
   {minor}%
   {A from b a}%
\Symb
   {minor}%
   {A from columns T}%
\Symb
   {minor}%
   {A from rows S}%
\Symb
   {minor}%
   {A without column a}%
\Symb
   {minor}%
   {A without columns T}%
\Symb
   {minor}%
   {A without row b}%
\Symb
   {minor}%
   {A without rows S}%
\Symb
   {$A\CRcirc$\Hyph linear combination}%
   {ACRcirc linear combination 1}%
\Symb
   {$A\CRcirc$\Hyph linear combination}%
   {ACRcirc linear combination 2}%
\Symb
   {active representation of group $G(f)$ in basis manifold $\mathcal B(f)$}%
   {active representation in basis manifold}%
\Symb
   {active representation of group $G(\Vector f)$ in basis manifold $\mathcal B(\Vector f)$}%
   {active representation in basis manifold, tower of representations}%
\Symb
   {affine space}%
   {affine space, division ring}%
\Symb
   {affine space}%
   {An}%
\Symb
   {associator of $R$\Hyph algebra}%
   {associator of algebra}%
\Symb
   {\subs row ($c$\hyph row) of matrix}%
   {c row}%
\Symb
   {commutator of $R$\Hyph algebra}%
   {commutator of algebra}%
\Symb
   {component of linear map $\Vector{A}$ of $D$\Hyph vector space}%
   {component of linear map, D vector space}%
\Symb
   {component $p$ of polylinear mapping $\Vector A$}%
   {component of polyadditive map, D vector space}%
\Symb
   {linear combination of vectors of $A$\Hyph module}%
   {CR linear combination in A module}%
\Symb
   {\CR power of element $A$ of biring}%
   {cr power}%
\Symb
   {\CR inverse element of biring}%
   {cr-inverse element}%
\Symb
   {\CR product of matrices}%
   {cr-product of matrices}%
\Symb
   {derivative of left shift}%
   {derivative of left shift}%
\Symb
   {derivative of left shift in $1$\Hyph parameter Lie group}%
   {derivative of left shift, 1-Parameter Group}%
\Symb
   {derivative of left shift in $1$\Hyph parameter Lie D group}%
   {derivative of left shift, 1-Parameter Group, algebra}%
\Symb
   {}%
   {derivative of right shift}%
\Symb
   {}%
   {derivative of right shift}%
\Symb
   {derivative of right shift in $1$\Hyph parameter Lie group}%
   {derivative of right shift, 1-Parameter Group}%
\Symb
   {derivative of right shift in $1$\Hyph parameter Lie D group}%
   {derivative of right shift, 1-Parameter Group, algebra}%
\Symb
   {derivative of left shift}%
   {derivative of Tstar shift}%
\Symb
   {\drc vector}%
   {drc vector}%
\Symb
   {hermitian conjugation in division ring}%
   {hermitian conjugation, division ring}%
\Symb
   {$(ji)$\hyph quasideterminant of matrix $\bfA$}%
   {j i quasideterminant definition}%
\Symb
   {$(^j_i)$\hyph \RC quasideterminant}%
   {j i RC-quasideterminant definition}%
\Symb
   {$(^j_i)$\hyph $\RCcirc$\Hyph quasideterminant}%
   {j i RCcirc-quasideterminant definition}%
\Symb
   {left shift in $D$\Hyph algebra}%
   {left shift, D algebra}%
\Symb
   {linear combination of vectors of $A$\Hyph module}%
   {linear combination in A module}%
\Symb
   {transformation of matrix}%
   {matrix, replacing its column}%
\Symb
   {transformation of matrix}%
   {matrix, replacing its row}%
\Symb
   {norm of map $\Vector A$ of normed $D$\hyph vector space}%
   {norm of map, D vector space}%
\Symb
   {opposite algebra to algebra $A$}%
   {opposite algebra}%
\Symb
   {orbit of linear mapping}%
   {orbit of linear mapping}%
\Symb
   {derivative}%
   {overline nabla_l, definition 2}%
\Symb
   {partial linear map of variable $v^i$}%
   {partial linear map of variable}%
\Symb
   {quasideterminant of matrix $\bfA$}%
   {quasideterminant definition}%
\Symb
   {\sups row ($r$\hyph row) of matrix}%
   {r row}%
\Symb
   {\RC power of element $A$ of biring}%
   {rc power}%
\Symb
   {\RC inverse element of biring}%
   {rc-inverse element}%
\Symb
   {\RC product of matrices}%
   {rc-product of matrices}%
\Symb
   {\RC quasideterminant}%
   {RC-quasideterminant definition}%
\Symb
   {$\RCcirc$\Hyph quasideterminant}%
   {RCcirc-quasideterminant definition}%
\Symb
   {\rcd vector}%
   {rcd vector}%
\Symb
   {right shift in $D$\Hyph algebra}%
   {right shift, D algebra}%
\Symb
   {coordinates of vector $\Vector a$ relative to Schauder basis}%
   {Schauder basis, coordinates}%
\Symb
   {set of $A$\Hyph linear mappings of module $\Vector V$ into module $\Vector W$}%
   {set A linear maps, module}%
\Symb
   {set of linear maps of division ring $D_1$ into division ring $D_2$}%
   {set linear maps, division ring}%
\Symb
   {set of polylinear mappings of rings $R_1$, ..., $R_n$ into module $S$}%
   {set polylinear maps, ring}%
\Symb
   {skew product of vectors $\Vector a_1$, ..., $\Vector a_m$}%
   {skew product of vectors}%
\Symb
   {standard component of tensor in tensor product of algebras}%
   {standard component of tensor, algebra}%
\Symb
   {right shift}%
   {starT shift}%
\Symb
   {\sT shift}%
   {starT shift, fibered group}%
\Symb
   {tensor power of algebra $A$}%
   {tensor power of algebra}%
\Symb
   {tensor product of algebras}%
   {tensor product of algebras}%
\Symb
   {left shift}%
   {Tstar shift}%
\Symb
   {\Ts shift}%
   {Tstar shift, fibered group}%
\Symb
   {anholonomic coordinates of vector}%
   {vector anholonomic coordinates}%
\Symb
   {holonomic coordinates of vector}%
   {vector holonomic coordinates}%

\SetIndexSpace
\Symb
   {basis manifold of \rcd vector space $\Vector V$}%
   {basis manifold of rcd vector space}%
\Symb
   {basis manifold of vector space}%
   {basis manifold of vector space}%
\Symb
   {basis manifold of representation $f$}%
   {basis manifold representation F algebra}%
\Symb
   {basis manifold of tower of representations $\Vector f$}%
   {basis manifold tower of representations}%
\Symb
   {basis manifold of affine space}%
   {Basis Manifold, Affine Space}%
\Symb
   {basis manifold of \rcd affine space}%
   {Basis Manifold, rcd Affine Space, division ring}%
\Symb
   {basis manifold of central affine space}%
   {BCAn}%
\Symb
   {basis manifold of Euclid space}%
   {BEn}%
\Symb
   {Cartesian power $\Bundle A$ of bundle $\Bundle B$}%
   {Cartesian power A of bundle B}%
\Symb
   {Cartesian power $A$ of set $B$}%
   {Cartesian power of set}%
\Symb
   {basis manifold of central affine space}%
   {FCAn}%
\Symb
   {basis manifold of Euclid space}%
   {FEn}%
\Symb
   {lattice of subrepresentations of representation $f$}%
   {lattice of subrepresentations}%
\Symb
   {lattice of towers of subrepresentations of tower of representations $\Vector f$}%
   {lattice of subrepresentations, tower of representations}%
\Symb
   {product of objects $B_1$, ..., $B_n$ in category $\mathcal A$}%
   {product of objects in category, 1 n}%
\Symb
   {structural constants of division ring $D$ over field $F$}%
   {structural constants of division ring over field}%
\Symb
   {tensor power of representation}%
   {tensor power of representation}%
\Symb
   {tensor product of representations}%
   {tensor product of representations}%

\SetIndexSpace
\Symb
   {central affine space}%
   {CAn}%
\Symb
   {central affine space}%
   {central affine space}%
\Symb
   {$j$th column determinant of matrix $\bfA$}%
   {column determinant}%
\Symb
   {$\CRcirc$\Hyph product of matrices of mappings}%
   {cr product of matrices of mappings}%
\Symb
   {left structural constant of Lie algebra}%
   {left structural constant of Lie algebra}%
\Symb
   {right structural constant of Lie algebra}%
   {right structural constant of Lie algebra}%
\Symb
   {structural constants of algebra $A$ over ring $D$}%
   {structural constants of algebra}%

\SetIndexSpace
\Symb
   {basis vector of representation of Lie group}%
   {basis vector of representation of Lie group}%
\Symb
   {basis vector of representation of Lie group over algebra $A$}%
   {basis vector of representation of Lie group over algebra A}%
\Symb
   {coordinates of basis vector of representation of Lie group over algebra $A$}%
   {basis vector of representation of Lie group over algebra A, coordinates}%
\Symb
   {coordinates of basis vector of representation of Lie group}%
   {basis vector of representation of Lie group, coordinates}%
\Symb
   {component of the G\^ateaux derivative of map $f(x)$ of algebra}%
   {component of Gateaux derivative of map, algebra}%
\Symb
   {component of the G\^ateaux derivative of map $\Vector f(\Vector x)$}%
   {component of Gateaux derivative of map, D vector space}%
\Symb
   {component of the G\^ateaux derivative of map $\Vector f(\Vector x)$}%
   {component of Gateaux derivative of map, D vector space, short}%
\Symb
   {component of the G\^ateaux differential of map $f(x)$}%
   {component of Gateaux derivative of map, division ring}%
\Symb
   {component of the G\^ateaux derivative of second order of map $f(x)$ of algebra}%
   {component of Gateaux derivative of Second Order, algebra}%
\Symb
   {component of the G\^ateaux derivative of second order of map $\Vector f(\Vector x)$}%
   {component of Gateaux derivative of Second Order, D vector space}%
\Symb
   {component of the G\^ateaux derivative of second order of map $f(x)$ of division ring}%
   {component of Gateaux derivative of Second Order, division ring}%
\Symb
   {conjugation in algebra}%
   {conjugation in algebra}%
\Symb
   {conjugation in ring}%
   {conjugation in ring}%
\Symb
   {coordinate \Drc vector bundle}%
   {coordinate drc vector bundle}%
\Symb
   {coordinate \rcd vector space}%
   {coordinate rcd vector space}%
\Symb
   {coordinate reference frame}%
   {coordinate reference frame, extensive definition}%
\Symb
   {diagonal in bundle $\Bundle A$}%
   {diagonal in bundle, 1}%
\Symb
   {direct product of division rings $D_1$, ..., $D_n$}%
   {direct product of division rings, 1 n}%
\Symb
   {double determinant of matrix $\bfA$}%
   {double determinant}%
\Symb
   {the Fr\'echet \Ds derivative of map $f$ of division ring}%
   {Frechet Dstar derivative of map, division ring}%
\Symb
   {the G\^ateaux \crd derivative of map $\Vector f$ of $D$\hyph vector space $\Vector V$ to $D$\hyph vector space $\Vector W$}%
   {Gateaux crd derivative of map, D vector space}%
\Symb
   {the G\^ateaux derivative of map $f$ of algebra}%
   {Gateaux derivative of map, algebra}%
\Symb
   {the G\^ateaux derivative of map $\Vector f$ of normed $D$\Hyph vector space $\Vector{V}$ to normed $D$\Hyph vector space $\Vector{W}$}%
   {Gateaux derivative of map, D vector space}%
\Symb
   {the G\^ateaux derivative of map $f$}%
   {Gateaux derivative of map, division ring}%
\Symb
   {the G\^ateaux derivative of map $f$ of algebra}%
   {Gateaux derivative of map, fraction, algebra}%
\Symb
   {the G\^ateaux derivative of map $f$}%
   {Gateaux derivative of map, fraction, division ring}%
\Symb
   {the G\^ateaux derivative of order $n$ of map $f$ of algebra}%
   {Gateaux derivative of Order n, algebra}%
\Symb
   {the G\^ateaux derivative of order $n$ of map $\Vector f$}%
   {Gateaux derivative of Order n, D vector space}%
\Symb
   {the G\^ateaux derivative of order $n$ of map $f$ of division ring}%
   {Gateaux derivative of Order n, division ring}%
\Symb
   {the G\^ateaux derivative of order $n$ of map $f$ of algebra}%
   {Gateaux derivative of Order n, fraction, algebra}%
\Symb
   {the G\^ateaux derivative of order $n$ of map $f$ of division ring}%
   {Gateaux derivative of Order n, fraction, division ring}%
\Symb
   {the G\^ateaux derivative of second order of mapping $f$ of algebra}%
   {Gateaux derivative of Second Order, algebra}%
\Symb
   {the G\^ateaux derivative of second order of map $\Vector f$}%
   {Gateaux derivative of Second Order, D vector space}%
\Symb
   {the G\^ateaux derivative of second order of map $f$ of division ring}%
   {Gateaux derivative of Second Order, division ring}%
\Symb
   {the G\^ateaux derivative of second order of mapping $f$ of algebra}%
   {Gateaux derivative of Second Order, fraction, algebra}%
\Symb
   {the G\^ateaux derivative of second order of map $f$ of division ring}%
   {Gateaux derivative of Second Order, fraction, division ring}%
\Symb
   {the G\^ateaux differential of mapping $f$ of algebra}%
   {Gateaux differential of map, algebra}%
\Symb
   {the G\^ateaux differential of map $\Vector f$ of normed $D$\Hyph vector space $\Vector{V}$ to normed $D$\Hyph vector space $\Vector{W}$}%
   {Gateaux differential of map, D vector space}%
\Symb
   {the G\^ateaux differential of map $f$}%
   {Gateaux differential of map, division ring}%
\Symb
   {the G\^ateaux differential of second order of mapping $f$ of algebra}%
   {Gateaux differential of Second Order, algebra}%
\Symb
   {the G\^ateaux differential of second order of mapping $\Vector f$}%
   {Gateaux differential of Second Order, D vector space}%
\Symb
   {the G\^ateaux differential of second order of mapping $f$ of division ring}%
   {Gateaux differential of Second Order, division ring}%
\Symb
   {the G\^ateaux \drc derivative of map $\Vector f$ of $D$\Hyph vector space $\Vector V$ to $D$\Hyph vector space $\Vector W$}%
   {Gateaux drc derivative of map, D vector space}%
\Symb
   {the G\^ateaux \Ds derivative of map $f$ of division ring $D$}%
   {Gateaux Dstar derivative of map, division ring}%
\Symb
   {the G\^ateaux Jacobian of map of $D$\Hyph vector space}%
   {Gateaux Jacobian of map, D vector space}%
\Symb
   {the G\^ateaux partial \drc derivative of map $f^b$ with respect to variable $v^a$}%
   {Gateaux partial crd derivative of map, 1, D vector space}%
\Symb
   {the G\^ateaux partial \drc derivative of map $f^b$ with respect to variable $v^a$}%
   {Gateaux partial crd derivative of map, 2, D vector space}%
\Symb
   {the G\^ateaux partial \drc derivative of map $f^b$ with respect to variable $v^a$}%
   {Gateaux partial crd derivative of map, 3, D vector space}%
\Symb
   {the G\^ateaux mixed partial derivative of map $f^j$ with respect to variables $v^i$, $v^j$}%
   {Gateaux partial derivative of Second Order, D vector space}%
\Symb
   {the G\^ateaux partial derivative of map $f^j$ with respect to variable $v^i$}%
   {Gateaux partial derivative, D vector space}%
\Symb
   {the G\^ateaux partial \drc derivative of map $f^b$ with respect to variable $v^a$}%
   {Gateaux partial drc derivative of map, 1, D vector space}%
\Symb
   {the G\^ateaux partial \drc derivative of map $f^b$ with respect to variable $v^a$}%
   {Gateaux partial drc derivative of map, 2, D vector space}%
\Symb
   {the G\^ateaux partial \drc derivative of map $f^b$ with respect to variable $v^a$}%
   {Gateaux partial drc derivative of map, 3, D vector space}%
\Symb
   {the G\^ateaux \sD derivative of map $f$ of division ring $D$}%
   {Gateaux starD derivative of map, division ring}%
\Symb
   {matrices vector space}%
   {matrices vector space}%
\Symb
   {Cartan derivative}%
   {overbrace D}%
\Symb
   {derivative}%
   {overline D}%
\Symb
   {derivative $e_{(k)}$}%
   {partial(k)}%
\Symb
   {\subs rows \drc vector space}%
   {r rows drc vector space}%
\Symb
   {speed of deviation}%
   {speed of deviation}%
\Symb
   {standard component of the G\^ateaux derivative of mapping $f$}%
   {standard component of Gateaux derivative, algebra}%
\Symb
   {standard component of the G\^ateaux differential of map $f$}%
   {standard component of Gateaux differential, division ring}%
\Symb
   {tensor product of division rings}%
   {tensor product of division rings}%
\Symb
   {tensor product of rings}%
   {tensor product of rings}%
\Symb
   {vector space type}%
   {vector space type}%

\SetIndexSpace
\Symb
   {$A\CRcirc$\Hyph basis for module}%
   {A CRcirc basis, module over algebra}%
\Symb
   {affine basis}%
   {Affine Basis}%
\Symb
   {basis of vector space}%
   {Basis e}%
\Symb
   {basis in vector space $\Vector V$}%
   {basis in V}%
\Symb
   {basis of vector space}%
   {basis, vector space}%
\Symb
   {basis of $(n)$\hyph vector space}%
   {basis,n vector space}%
\Symb
   {Cartesian power of total spaces}%
   {Cartesian power of total spaces}%
\Symb
   {Cartesian product of total spaces}%
   {Cartesian product of total spaces, definition 1}%
\Symb
   {central affine basis}%
   {Central Affine Basis}%
\Symb
   {basis for \Drc vector bundle}%
   {drc basis, vector bundle}%
\Symb
   {form of reference frame}%
   {dual forms, reference frame}%
\Symb
   {Euclid space}%
   {Euclid space}%
\Symb
   {Euclid space}%
   {Euclid space, division ring}%
\Symb
   {identical transformation of bundle}%
   {identical transformation of bundle}%
\Symb
   {linear automorphism of quaternioin algebra}%
   {mapping E, quaternion}%
\Symb
   {linear automorphism of quaternioin algebra}%
   {mapping E_1, quaternion}%
\Symb
   {linear automorphism of quaternioin algebra}%
   {mapping E_2, quaternion}%
\Symb
   {orthonornal basis}%
   {Orthonornal Basis}%
\Symb
   {pseudo Euclid space}%
   {pseudo Euclid space}%
\Symb
   {pseudo Euclid space}%
   {pseudo Euclid space, division ring}%
\Symb
   {quaternion algebra over the field $F$}%
   {quaternion algebra over the field}%
\Symb
   {quaternion division algebra over the field}%
   {quaternion division algebra over the fieldL}%
\Symb
   {$\RCcirc A$\Hyph linear combination}%
   {RCcircA linear combination 1}%
\Symb
   {$\RCcirc A$\Hyph linear combination}%
   {RCcircA linear combination 2}%
\Symb
   {\rcd affine basis}%
   {rcd affine basis, division ring}%
\Symb
   {reduced Cartesian product of total spaces}%
   {reduced Cartesian product of total spaces, definition 1}%
\Symb
   {Schauder basis}%
   {Schauder basis}%
\Symb
   {set of nonsingular \sT transformations of bundle $\Bundle E$}%
   {set of starT nonsingular transformations of bundle}%
\Symb
   {set of nonsingular \Ts transformations of bundle $\Bundle E$}%
   {set of Tstar nonsingular transformations of bundle}%
\Symb
   {standard coordinates of basis}%
   {standard coordinates of basis}%
\Symb
   {standard coordinates of reference frame}%
   {standard coordinates of reference frame}%
\Symb
   {vector field of reference frame}%
   {vector field of reference frame}%
\Symb
   {vector of basis}%
   {vector of basis}%

\SetIndexSpace
\Symb
   {coordinates of basis in \sups rows \rcd vector space}%
   {basis coordinates, c rows rcd vector space}%
\Symb
   {coordinates of basis in \subs rows \drc vector space}%
   {basis coordinates, r rows drc vector space}%
\Symb
   {basis for \sups rows \rcd vector space}%
   {basis, c rows rcd vector space}%
\Symb
   {basis for \subs rows \drc vector space}%
   {basis, r rows drc vector space}%
\Symb
   {central affine basis}%
   {Central Affine Basis, division ring}%
\Symb
   {component of linear map $f$ of division ring}%
   {component of linear map, division ring}%
\Symb
   {component of polylinear map into associative algebra}%
   {component of polylinear map, associative algebra}%
\Symb
   {component of polylinear map of division ring}%
   {component of polylinear map, division ring}%
\Symb
   {fibered morphism from bundle $\Bundle A$ into $\Bundle B$}%
   {fibered morphism from A into B}%
\Symb
   {filter $\mathfrak{F}$ converges to set $A$}%
   {filter converges}%
\Symb
   {homomorphism of fibered universal algebras}%
   {homomorphism of fibered universal algebras}%
\Symb
   {inverse fibered correspondence}%
   {inverse fibered correspondence, 1}%
\Symb
   {inverse reduced fibered correspondence}%
   {inverse reduced fibered correspondence, 1}%
\Symb
   {map to Cartesian product}%
   {map to Cartesian product}%
\Symb
   {norm of mapping into $D$\Hyph algebra}%
   {norm of map, algebra}%
\Symb
   {norm of map $f$ of division ring}%
   {norm of map, division ring}%
\Symb
   {representation orbit of group $G$}%
   {orbit of Tstar representation of group}%
\Symb
   {orthonornal basis}%
   {Orthonornal Basis, division ring}%
\Symb
   {quaternion algebra  over field ${\rm {\mathbb{F}}}$}%
   {quaternion algebra F a b}%
\Symb
   {reference frame}%
   {reference frame}%
\Symb
   {reference frame, extensive definition}%
   {reference frame, extensive definition}%
\Symb
   {standard component of biadditive map $f$ over field $F$}%
   {standard component of biadditive map, division ring}%
\Symb
   {standard $F$\Hyph component of linear mapping $f$}%
   {standard component of linear map, division ring}%
\Symb
   {standard component of polylinear mapping into associative algebra}%
   {standard component of polylinear map, associative algebra}%
\Symb
   {standard component of polylinear map $f$ of division ring}%
   {standard component of polylinear map, division ring}%
\Symb
   {standard component of quadratic map $f$ over field $F$}%
   {standard component of quadratic map, division ring}%
\Symb
   {standard component of tensor}%
   {standard component of tensor, division ring}%

\SetIndexSpace
\Symb
   {affine transformation group}%
   {affine transformation group}%
\Symb
   {\CR matrix group}%
   {cr-matrix group}%
\Symb
   {affine transformation group}%
   {drc affine transformation group}%
\Symb
   {fibered little group of section $h$}%
   {fibered little group}%
\Symb
   {fibered stability group of section $h$}%
   {fibered stability group}%
\Symb
   {algebra Lie of group Lie}%
   {g}%
\Symb
   {left defined algebra Lie of group Lie}%
   {gl}%
\Symb
   {right defined algebra Lie of group Lie}%
   {gr}%
\Symb
   {group of automorphisms of representation $f$}%
   {group of automorphisms of representation}%
\Symb
   {group of homomorphisms of vector space $\Vector V$}%
   {GV}%
\Symb
   {little group of $x$}%
   {little group}%
\Symb
   {orbit of effective covariant \sT representation of fibered group}%
   {orbit of effective covariant starT representation of fibered group}%
\Symb
   {orbit of effective covariant \sT representation of group}%
   {orbit of effective covariant starT representation of group}%
\Symb
   {orbit of effective covariant \Ts representation of fibered group}%
   {orbit of effective covariant Tstar representation of fibered group}%
\Symb
   {orbit of effective covariant \Ts representation of group}%
   {orbit of effective covariant Tstar representation of group}%
\Symb
   {product of groups $G_1$, ..., $G_n$}%
   {product of groups, 1 n}%
\Symb
   {\RC matrix group}%
   {rc-matrix group}%
\Symb
   {stability group of $x$}%
   {stability group}%

\SetIndexSpace
\Symb
   {Hadamard inverse of matrix}%
   {Hadamard inverse of matrix}%
\Symb
   {quaternion algebra}%
   {quaternion algebra H a b}%
\Symb
   {quaternion algebra over real field}%
   {quaternion algebra over real field}%

\SetIndexSpace
\Symb
   {infinitesimal generator of representation}%
   {infinitesimal generator I of representation}%
\Symb
   {infinitesimal generator of representation}%
   {infinitesimal generator i of representation}%
\Symb
   {Lie group infinitesimal generators}%
   {Lie group infinitesimal generators}%
\Symb
   {vector module of algebra $A$}%
   {vector module of algebra}%
\Symb
   {vector module of ring $D$}%
   {vector module of ring}%
\Symb
   {vector of element $d$ of algebra}%
   {vector of algebra}%
\Symb
   {vector of mapping $f$}%
   {vector of mapping}%
\Symb
   {vector of element $d$ of ring}%
   {vector of ring}%

\SetIndexSpace
\Symb
   {Jacobian matrix of left shift}%
   {aE, quaternion, Jacobian matrix}%
\Symb
   {closure operator of representation $f$}%
   {closure operator, representation}%
\Symb
   {closure operator of tower of representations $\Vector f$}%
   {closure operator, tower of representations}%
\Symb
   {Jacobian matrix of right shift}%
   {Ea, quaternion, Jacobian matrix}%
\Symb
   {tower of subrepresentations of tower of representations $\Vector f$ generated by tuple of sets $\VX X$}%
   {subrepresentation generated by tuple of sets}%

\SetIndexSpace
\Symb
   {kernel of linear mapping of $D$\Hyph vector space}%
   {kernel of linear map, D vector space}%
\Symb
   {kernel of linear mapping of division ring}%
   {kernel of linear map, division ring}%

\SetIndexSpace
\Symb
   {left $ij$th cofactor of entry of matrix}%
   {left cofactor, matrix}%
\Symb
   {left double $ij$th cofactor of entry of matrix}%
   {left double cofactor}%
\Symb
   {left shift}%
   {left shift}%
\Symb
   {Lie derivative of connection}%
   {Lie derivative of connection}%
\Symb
   {Lie derivative of metric}%
   {Lie derivative of metric}%
\Symb
   {limit of correspondence $\Phi$ with respect to the filter $\mathfrak{F}$}%
   {limit of correspondence with respect to the filter}%
\Symb
   {limit of sequence in normed ring}%
   {limit of sequence, normed ring}%
\Symb
   {limit of sequence in valued division ring}%
   {limit of sequence, valued division ring}%
\Symb
   {limit of sequence in valued ring}%
   {limit of sequence, valued ring}%
\Symb
   {passive transformation}%
   {passive transformation}%
\Symb
   {set of \Acr linear mappings of module $\Vector V$ into module $\Vector W$}%
   {set Acr linear maps, module}%
\Symb
   {\rcd vector space of \drc linear maps of \drc vector space $\Vector V$ into \drc vector space $\Vector W$}%
   {set drc linear maps, drc vector space}%
\Symb
   {set of linear mappings of algebra $A_1$ into algebra $A_2$}%
   {set linear mappings, algebra}%
\Symb
   {set of linear mappings of $D$\Hyph vector space $\Vector{V}$ into $D$\Hyph vector space $\Vector{W}$}%
   {set linear maps, D vector space}%
\Symb
   {set of left-side nonsingular transformations of set $M$}%
   {set of left-side nonsingular transformations}%
\Symb
   {set of $n$\hyph linear mappings of algebra $A$ into module $S$}%
   {set polylinear mappings An, algebra}%
\Symb
   {set of polylinear maps of algebras $A_1$, ..., $A_n$ into module $S$}%
   {set polylinear mappings, algebra}%
\Symb
   {set of polylinear mappings}%
   {set polylinear maps, D vector space}%
\Symb
   {set of polylinear maps of algebras $A_1$, ..., $A_n$ into algebra $A$}%
   {set polylinear maps, finite dimensional algebra}%
\Symb
   {\drc vector space of \rcd linear maps of \rcd vector space $\Vector{V}$ into \rcd vector space $\Vector{W}$}%
   {set rcd linear maps, rcd vector space}%
\Symb
   {set of \sT representations of division ring $S$ in additive group of division ring $R$}%
   {set sT representations, division ring}%
\Symb
   {set of \Ts representations of division ring $S$ in additive group of division ring $R$}%
   {set Ts representations, division ring}%

\SetIndexSpace
\Symb
   {set of \sT transformations of set $M$}%
   {set of starT transformations}%
\Symb
   {set of transformations of set $M$}%
   {set of transformations}%
\Symb
   {set of \Ts transformations of set $M$}%
   {set of Tstar transformations}%
\Symb
   {space of orbits of effective \sT covariant representation of the group}%
   {space of orbits of effective sT representation}%
\Symb
   {space of orbits of effective \Ts covariant representation of the group}%
   {space of orbits of effective Ts representation}%
\Symb
   {space of orbits of \Ts representation $f$ of group $G$ in set $M$}%
   {space of orbits of Ts representation}%

\SetIndexSpace
\Symb
   {norm of quaternion $x$}%
   {norm, quaternion algebra}%
\Symb
   {nucleus of $R$\Hyph algebra $A$}%
   {nucleus of algebra}%

\SetIndexSpace
\Symb
   {geometric object in coordinate representation defined in \rcd vector space}%
   {geometric object, coordinate rcd vector space}%
\Symb
   {geometric object in coordinate representation}%
   {geometric object, coordinate vector space}%
\Symb
   {geometric object defined in \rcd vector space}%
   {geometric object, rcd vector space}%
\Symb
   {octonion algebra}%
   {octonion algebra}%
\Symb
   {orbit of representation of fibered group $\Bundle G$}%
   {orbit of representation of fibered group}%
\Symb
   {orbit of representation of the group $G$}%
   {orbit of representation of group}%

\SetIndexSpace
\Symb
   {bundle}%
   {bundle}%
\Symb
   {bundle of level $2$}%
   {bundle of level 2}%
\Symb
   {bundle of level $n$}%
   {bundle of level n}%
\Symb
   {Cartesian power $n$ of bundle $\bundle{}{p}{E}{}$}%
   {Cartesian power of bundle}%
\Symb
   {Cartesian product of bundles}%
   {Cartesian product of bundles, definition 1}%
\Symb
   {passive representation of group $G(f)$ in basis manifold $\mathcal B(f)$}%
   {passive representation in basis manifold}%
\Symb
   {passive representation of group $G(\Vector f)$ in basis manifold $\mathcal B(\Vector f)$}%
   {passive representation in basis manifold, tower of representations}%
\Symb
   {reduced Cartesian product of bundles}%
   {reduced Cartesian product of bundles, definition 1}%
\Symb
   {set of nonsingular \sT transformations of bundle $\bundle{}pE{}$}%
   {set of starT nonsingular transformations of bundle, projection}%
\Symb
   {set of nonsingular \Ts transformations of bundle $\bundle{}pE{}$}%
   {set of Tstar nonsingular transformations of bundle, projection}%

\SetIndexSpace
\Symb
   {active transformation}%
   {active transformation}%
\Symb
   {\sups rows \rcd vector space}%
   {c rows rcd vector space}%
\Symb
   {Cartan curvature}%
   {Cartan curvature}%
\Symb
   {\CR rank of matrix}%
   {cr-rank of matrix}%
\Symb
   {diagonal in bundle  $\bundle{}pA{}$}%
   {diagonal in bundle, 2}%
\Symb
   {diagonal in bundle $\Bundle A$}%
   {diagonal in reduced bundle, 2}%
\Symb
   {\Ds component of coordinates of vector $\Vector r$}%
   {Dstar component of coordinates of vector, D vector space}%
\Symb
   {image of $m$ under endomorphism $R$ of effective representation}%
   {endomorphism image, effective representation}%
\Symb
   {image of tuple $\VX a$ under endomorphism $\VX r$ of tower of effective representations}%
   {endomorphism image, tower of effective representations}%
\Symb
   {curvature}%
   {GLn curvature_overline}%
\Symb
   {$\RCcirc$\Hyph product of matrices of mappings}%
   {rc product of matrices of mappings}%
\Symb
   {\RC rank of matrix}%
   {rc-rank of matrix}%
\Symb
   {right $ij$th cofactor of entry of matrix}%
   {right cofactor, matrix}%
\Symb
   {right double $ij$th cofactor of entry of matrix}%
   {right double cofactor}%
\Symb
   {right shift}%
   {right shift}%
\Symb
   {$i$th row determinant of matrix $\bfA$}%
   {row determinant}%
\Symb
   {scalar algebra of algebra $A$}%
   {scalar algebra of algebra}%
\Symb
   {scalar algebra of ring $D$}%
   {scalar algebra of ring}%
\Symb
   {scalar of element $d$ of algebra}%
   {scalar of algebra}%
\Symb
   {scalar of mapping $f$}%
   {scalar of mapping}%
\Symb
   {scalar of element $d$ of ring}%
   {scalar of ring}%
\Symb
   {set of right-side nonsingular transformations of set $M$}%
   {set of right-side nonsingular transformations}%
\Symb
   {\sD component of coordinates of vector $\Vector r$}%
   {starD component of coordinates of vector, D vector space}%

\SetIndexSpace
\Symb
   {composition of fibered correspondences}%
   {composition of fibered correspondences}%
\Symb
   {inverse fibered correspondence}%
   {inverse fibered correspondence, 2}%
\Symb
   {inverse reduced fibered correspondence}%
   {inverse reduced fibered correspondence, 2}%
\Symb
   {linear span in vector space}%
   {linear span, vector space}%
\Symb
   {image of basis $X$ under passive transformation $S$}%
   {passive transformation of basis, representation}%
\Symb
   {image of basis $\VX  X$ under passive transformation $\VX s$}%
   {passive transformation of basis, tower of representations}%
\Symb
   {symmetric group}%
   {symmetric group}%

\SetIndexSpace
\Symb
   {category of \Ts representations of $\Omega_1$\Hyph algebra $A$}%
   {category of Tstar representations of Omega1 algebra}%
\Symb
   {category of \Ts representations of $\Omega_1$\Hyph algebra from category $\mathcal A$}%
   {category of Tstar representations of Omega1 algebra from category}%
\Symb
   {tangent plane to group $G$}%
   {TaG}%
\Symb
   {trace of quaternion $x$}%
   {trace, quaternion algebra}%

\SetIndexSpace
\Symb
   {coordinate vector space}%
   {coordinate vector space}%
\Symb
   {coordinates in vector space}%
   {coordinates in vector space}%
\Symb
   {direct product of $\RCstar D_i$\hyph vector spaces $\Vector V_1$, ..., $\Vector V_n$}%
   {direct product, rcd vector space, 1 n}%
\Symb
   {dual space of \rcd vector space $\Vector V$}%
   {dual space of rcd vector space}%
\Symb
   {hermitian conjugated vector}%
   {hermitian conjugated vector}%
\Symb
   {\dcr vector space}%
   {left CR vector space}%
\Symb
   {\drc vector space}%
   {left RC vector space}%
\Symb
   {\crd vector space}%
   {right CR vector space}%
\Symb
   {\rcd vector space}%
   {right RC vector space}%
\Symb
   {tensor product of $D$\Hyph vector spaces}%
   {tensor product of D vector spaces}%
\Symb
   {tensor product of \Ds vector spaces}%
   {tensor product of Dstar vector spaces}%
\Symb
   {vector space}%
   {V}%

\SetIndexSpace
\Symb
   {set of coordinates of representation $J(f,X)$}%
   {coordinate set of representation}%
\Symb
   {set of tuples of coordinates of tower of representations $\Vector J(\Vector f,\VX X)$}%
   {coordinate set of tower of representations}%
\Symb
   {coordinates of basis $X'$ relative to basis $X$ of representation}%
   {coordinates of basis relative to basis, representation}%
\Symb
   {coordinates of element $m$ of representation $f$ relative to set $X$}%
   {coordinates of element relative to generating set, representation}%
\Symb
   {coordinates of element $m$ relative to set $X$}%
   {coordinates of element relative to set, representation}%
\Symb
   {tuple of coordinates of element $\Vector a*$ relative to tuple of sets $\VX X$}%
   {coordinates of element, tower of representations}%
\Symb
   {geometric object in coordinate representation defined in $\Omega_2$\Hyph algebra $M$}%
   {geometric object, coordinate representation g}%
\Symb
   {geometric object in coordinate representation defined in tuple of $\VX\Omega$\Hyph algebras $\VX A$}%
   {geometric object, coordinate tower of representations g}%
\Symb
   {geometric object defined in $\Omega_2$\Hyph algebra $M$}%
   {geometric object, representation g}%
\Symb
   {geometric object defined in tuple of $\VX\Omega$\Hyph algebras $\VX A$}%
   {geometric object, tower of representations g}%
\Symb
   {geometric object in vector space}%
   {geometric object, vector space}%
\Symb
   {set of coordinates of set $B\subset J(f,X)$}%
   {subset of coordinates of representation}%
\Symb
   {coordinates of tuple of sets $\VX B$ relative to tuple of sets $\VX X$}%
   {subset of coordinates of tower of representations}%
\Symb
   {coordinates of set $B_k$ relative to tuple of sets $\VX X$}%
   {subset of coordinates of tower of representations, k}%
\Symb
   {set of $\Omega_2$\Hyph words representing set $B\subset J(f,X)$}%
   {subset of words of representation}%
\Symb
   {superposition of coordinates of the representation $f$ and the element $m$}%
   {superposition of coordinates, representation}%
\Symb
   {superposition of coordinates of the tower of representations $\Vector f$ and the element $\VX a$}%
   {superposition of coordinates, tower of representations}%
\Symb
   {$\Omega_2$\Hyph word representing element $m\in J(f,X)$}%
   {word of element relative to generating set, representation}%
\Symb
   {set of $\Omega_2$\Hyph words of representation $J(f,X)$}%
   {word set of representation}%
\Symb
   {set of tuples of $\VX{\Omega}$\Hyph words of tower of representations $\Vector J(\Vector f,\VX X)$}%
   {word set of tower of representations}%
\Symb
   {tuple of words of element $\Vector a*$ relative to tuple of sets $\VX X$}%
   {words of element, tower of representations}%

\SetIndexSpace
\Symb
   {conjugate of quaternion $x$}%
   {conjugate of quaternion}%
\Symb
   {local basis of affine space}%
   {local basis of affine space}%
\Symb
   {anholonomic coordinate}%
   {x(k)}%

\SetIndexSpace
\Symb
   {center of an $R$\Hyph algebra $A$}%
   {center of algebra}%
\Symb
   {center of ring $D$}%
   {center of ring}%

\SetIndexSpace
\Symb
   {deviation of trajectories}%
   {deviation of trajectories}%
\Symb
   {identical transformation}%
   {identical transformation}%
\Symb
   {image of vector $\Vector e_k\in\Basis e$ under isomorphism to coordinate vector space}%
   {image of vector e_k, coordinate vector space}%
\Symb
   {Kronecker symbol}%
   {Kronecker symbol}%

\SetIndexSpace
\Symb
   {anholonomic coordinates of connection}%
   {anholonomic coordinates of connection}%
\Symb
   {Cartan symbol}%
   {Cartan symbol}%
\Symb
   {connection}%
   {conection overline}%
\Symb
   {connection coefficients in $D$\Hyph affine space}%
   {connection coefficients, D affine space}%
\Symb
   {connection in $D$\Hyph affine manifold}%
   {connection, affine manifold}%
\Symb
   {$D$\Hyph affine connection coefficients on manifold}%
   {D affine connection coefficients, manifold}%
\Symb
   {holonomic coordinates of connection}%
   {holonomic coordinates of connection}%
\Symb
   {Cartan connection}%
   {overbrace Gamma i kl}%
\Symb
   {set of sections of bundle}%
   {set of sections of bundle}%

\SetIndexSpace
\Symb
   {inverse operator to operator $\psi_l$}%
   {inverse operator to operator psi l}%
\Symb
   {inverse operator to operator $\psi_r$}%
   {inverse operator to operator psi r}%

\SetIndexSpace
\Symb
   {anholonomity object}%
   {anholonomity object}%

\SetIndexSpace
\Symb
   {left basic operator of Lie group over algebra $A$}%
   {L basic operator of Lie group over algebra A}%
\Symb
   {left basic operator of group Lie}%
   {Lie Basic Operator L}%
\Symb
   {left basic operator of Lie 1-parameter group}%
   {Lie Basic Operator L, 1-Parameter Group}%
\Symb
   {left basic operator of Lie 1-parameter group over algebra $A$}%
   {Lie Basic Operator L, 1-Parameter Group, algebra}%
\Symb
   {right basic operator of group Lie}%
   {Lie Basic Operator R}%
\Symb
   {right basic operator of Lie 1-parameter group}%
   {Lie Basic Operator R, 1-Parameter Group}%
\Symb
   {right basic operator of Lie 1-parameter group over algebra $A$}%
   {Lie Basic Operator R, 1-Parameter Group, algebra}%
\Symb
   {right basic operator of Lie group over algebra $A$}%
   {R basic operator of Lie group over algebra A}%

\SetIndexSpace
\Symb
   {Lie group composition law}%
   {Lie group composition law}%

\SetIndexSpace
\Symb
   {Cartan derivative}%
   {overbrace nabla_l}%
\Symb
   {derivative}%
   {overline nabla_l, definition 1}%

\SetIndexSpace
\Symb
   {restriction of correspondence $\Phi$ to set $C$}%
   {restriction of correspondence}%

\SetIndexSpace
\Symb
   {Cartesian product of bundles}%
   {Cartesian product of bundles, definition 2}%
\Symb
   {Cartesian product of total spaces}%
   {Cartesian product of total spaces, definition 2}%
\Symb
   {direct product of division rings $D_i$, $i\in I$}%
   {direct product of division rings}%
\Symb
   {direct product of division rings $D_1$, ..., $D_n$}%
   {direct product of division rings, i 1 n}%
\Symb
   {direct product of $\RCstar D_i$\hyph vector spaces $\Vector V_i$, $i\in I$}%
   {direct product, rcd vector space}%
\Symb
   {direct product of $\RCstar D_i$\hyph vector spaces}%
   {direct product, rcd vector space, i 1 n}%
\Symb
   {product of groups $G_i$, $i\in I$}%
   {product of groups}%
\Symb
   {product of groups $G_1$, ..., $G_n$}%
   {product of groups, i 1 n}%
\Symb
   {product of objects $\{B_i,i\in I\}$ in category $\mathcal A$}%
   {product of objects in category}%
\Symb
   {product of objects $B_1$, ..., $B_n$ in category $\mathcal A$}%
   {product of objects in category, i 1 n}%
\Symb
   {reduced Cartesian product of bundles}%
   {reduced Cartesian product of bundles, definition 2}%
\Symb
   {reduced Cartesian product of total spaces}%
   {reduced Cartesian product of total spaces, definition 2}%

\SetIndexSpace
\Symb
   {fibered subset}%
   {fibered subset}%
\Symb
   {subbundle}%
   {subbundle}%

\CloseIndex